%% file: spincentralizerqmath.tex
\documentclass[a4paper,12pt]{amsart}
\usepackage{amsmath}
\usepackage{amssymb}
\usepackage{amsthm}
\usepackage{verbatim}

\numberwithin{equation}{section}

\newtheorem{thm}{Theorem}[section]
\newtheorem{cor}[thm]{Corollary}
\newtheorem{lem}[thm]{Lemma}

\newtheorem{defn}[thm]{Definition}
\newtheorem{exmp}[thm]{Example}
\newtheorem{rem}[thm]{Remark}
\newcommand{\ds}{\displaystyle}

\newcommand{\SOe}{SO(2n)}

\newcommand{\Syp}{Sp(2n)}

\newcommand{\al}{\alpha}

\newcommand{\lam}{\lambda}

\newcommand{\ep}{\epsilon}
\newcommand{\be}{{\bold{e}}}
\newcommand{\bx}{\bold{x}}

\newcommand{\diag}{\operatorname{diag}}

\newcommand{\ad}{\operatorname{ad}}
\newcommand{\id}{\operatorname{id}}

\newcommand{\sign}{\operatorname{\ep}}
\newcommand{\Hom}{\operatorname{Hom}}
\newcommand{\End}{\operatorname{End}}

\newcommand{\ua}[1]{ \ensuremath{u_{#1}}}
\newcommand{\ub}[1]{ \ensuremath{u_{\overline{#1}}}}
\newcommand{\uo}{\ensuremath{u_{0'}}}

\newcommand{\bfv}[1]{\ensuremath{\mathbf{v}_{#1}}}
\newcommand{\bfx}[1]{\ensuremath{\mathbf{x}_{#1}}}
\newcommand{\bfy}[1]{\ensuremath{\mathbf{y}_{#1}}}

\newcommand{\C}{\ensuremath{\mathbb{C}}}
\newcommand{\R}{\ensuremath{\mathbb{R}}}
\newcommand{\G}{\ensuremath{\frak{so}(N)}}
\newcommand{\GC}{\ensuremath{\frak{so}(N,\C)}}

\newcommand{\ur}[1]{\ensuremath{{\mathtt{#1}}}}

\newcommand{\ulb}[1]{\ensuremath{\overline{{\mathtt{#1}}}}}
\newcommand{\dphi}{\ensuremath{\operatorname{d\phi}}}
\newcommand{\tW}{\ensuremath{\mathtt{W}}}
\newcommand{\tS}{\ensuremath{\mathtt{S}}}
\newcommand{\tI}{\ensuremath{\mathtt{I}}}
\newcommand{\tJ}{\ensuremath{\mathtt{J}}}
\newcommand{\tK}{\ensuremath{\mathtt{K}}}
\newcommand{\tWc}{\ensuremath{\mathtt{W^c}}}
\newcommand{\pr}[1]{\ensuremath{\operatorname{pr}_{#1}}}
\newcommand{\inj}[1]{\ensuremath{\operatorname{inj}_{#1}}}
\newcommand{\tT}{\ensuremath{\mathtt{T}}}
\newcommand{\Alt}[1]{\ensuremath{\operatorname{Alt_{#1}}}}
\newcommand{\cont}[1]{\ensuremath{\operatorname{C}_{#1}}}

\begin{document}

\title{Spin representations and centralizer algebras for the Spinor groups}

\date{\today}
\author{Kazuhiko Koike}
\thanks{Supported by the Grant-in-Aid for Scientific Research}
\address{
Department of Mathematics, 
Aoyama Gakuin University,
5-10-1, Sagamihara-shi, fuchinobe, 229-8558,Japan}
\email{koike@gem.aoyama.ac.jp}

\maketitle

\markboth{KAZUHIKO KOIKE}{CENTRALIZER ALGEBRAS FOR SPINOR GROUPS}

\begin{abstract}
We pursue an analogy of the 
Schur-Weyl reciprocity for the spinor groups and pick up the irreducible spin representations in the tensor space $\Delta \textstyle{\bigotimes \bigotimes^k V}$. Here $\Delta$ is the fundamental representation of $Pin(N)$ and $V$ is the natural (vector) representation of the orthogonal group $O(N)$.
 We consider the centralizer algebra $\mathbf{CP_k} = \End_{Pin(N)}(\Delta 
\textstyle{\bigotimes \bigotimes^k V})$
 for $Pin(N)$, the double covering group of $O(N)$ and define two kinds of linear basis in $\mathbf{CP_k}$ (one comes from invariant theory and the other from representation theory), both of which are parameterized by the 'generalized Brauer diagrams'. We develop analogous argument to the original Brauer centralizer algebra for $O(N)$ and determine the transformation matrices between the above two basis and give the multiplication rules of those basis. Finally we define the subspaces in $\Delta \textstyle{\bigotimes \bigotimes^k V}$, on which the symmetric group $\frak{S}_k$ and $Pin(N)$ or $Spin(N)$ act as a dual pair.
\end{abstract}

\section{Introduction}
In this paper we consider analogs of the Brauer centralizer algebras for the spinor groups and define the subspace $T_{k,s}^{0}$ of the tensor space  $\Delta \bigotimes \bigotimes^k V$, on which the symmetric group $\frak{S}_k$ and $Pin(N)$ act as the dual pair. (See Definition \ref{def T0} and Theorem \ref{thm dualpin}.)
Here $V \cong \R^{N}$ is the vector (natural) representation of $O(N)$ and $\Delta$ is the fundamental representation of $Pin(N)$.

For the irreducible representations of the orthogonal groups, together with the other classical groups,  this procedure is carried out by H. Weyl in his book {\it The Classical Groups}. 
Namely the centralizer algebra of the orthogonal group in $\bigotimes^k V$ is generated by the symmetric group $\frak{S}_k$ and the contractions and the immersions of the invariant form and is called the Brauer centralizer algebras.
  The basis of this algebra is parameterized by the diagrams, which is called the Brauer diagrams and the multiplication rules of the basis elements are given in \cite{br}.
 H. Weyl defined the subspace consisting of the traceless tensors in $\bigotimes^k V$ and he showed that the symmetric group $\frak{S}_k$ and the orthogonal group act on this space as the dual pair and he picked up the irreducible representations of $O(N)$ as the images of the Young symmetrizers corresponding to the highest weights of the irreducible representations.

\begin{figure}[htbp]
\begin{minipage}{\linewidth}
\begin{center}
\input{spincentfig1.tex}
\end{center}
\caption{The generalized Brauer diagrams of $k=2$ and $l=2$}
\end{minipage}
\label{spincentfig1}
\end{figure}
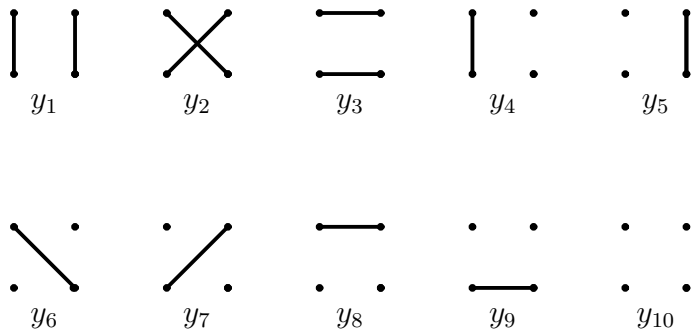

So we follow this line and first determine the basis of the centralizer algebra
$$
\mathbf{CP_k} = \Hom_{Pin(N)}(\Delta \textstyle{\bigotimes \bigotimes\limits^k V, \Delta \bigotimes \bigotimes\limits^k V}).
$$
In this case, we can define two kinds of basis of the centralizer algebra $\mathbf{CP_k}$ and 
the basis elements of both bases are parameterized by the same diagrams, which we call \lq the generalized Brauer diagrams'. (See figure \ref{spincentfig1}.)

The generalized Brauer diagrams are, by definition, 
 the graphs of two lines of vertices with $k$ vertices in the upper row and $l$ vertices in the lower row, in which vertices are connected with each other as in the usual Brauer diagrams except for admitting isolated vertices. Namely they are graphs with no loops in which the number of edges connected to each vertex is either $0$ or $1$.

Let us denote the set of the generalized Brauer diagrams with $k$ vertices in the upper row and $l$ vertices in the lower row by $\mathbf{GB^k_l}$.

Precisely speaking, for each element of $\mathbf{GB^k_l}$, we define two elements in $\mathbf{CP^k_l} = \Hom_{Pin(N)}(\Delta \bigotimes \bigotimes^k V, \Delta \bigotimes \bigotimes^l V)$, one of which comes from the invariant theory (we call this element an invariant theoretic element.) and the other of which comes from the representation theory (we call this element a representation theoretic element.).
The edges in the generalized Brauer diagram convey the same meaning as in the usual Brauer diagrams. As for the exact definitions of the above elements, see the part after Lemma \ref{lem indep} in the section \ref{sect inv} for the invariant theoretic elements and the part after Definition \ref{def injind} in the section \ref{sect reppar} for the representation theoretic elements.

For any element $d \in \mathbf{GB^k_l}$, let us denote the corresponding invariant theoretic element by $d_{inv}$ and the corresponding representation theoretic element by $d_{rt}$ respectively.  Then we can show that both of the sets $\{d_{inv} \vert d \in \mathbf{GB^k_l}\}$
 and $\{d_{rt} \vert d \in \mathbf{GB^k_l}\}$ are bases of the centralizer space  $\mathbf{CP^k_l} = \Hom_{Pin(N)}(\Delta \bigotimes \bigotimes^k V, \Delta \bigotimes \bigotimes^l V)$ if $k \leqq n$ and $l \leqq n$, here $N = 2n$ or $N = 2n +1$.
If $k = l$ and $k \leqq n$, we obtain two kinds of basis of the centralizer algebra $\mathbf{CP_k}$.

In Theorem \ref{thm psitorep} and Theorem \ref{thm psitorep2}, 
we give the transformation matrices between the above two bases explicitly.

If $k > n$ or $l > n$, the sets $\{d_{inv} \vert d \in \mathbf{GB^k_l}\}$
 and $\{d_{rt} \vert d \in \mathbf{GB^k_l}\}$  are not basis of  $\mathbf{CP^k_l}$ at all, but the invariant theory tells us that the elements of the set $\{d_{inv} \vert d \in \mathbf{GB^k_l}\}$ still span the whole space $\mathbf{CP^k_l}$. (Lemma \ref{lem indep}.)

We also give a series of the formulas to calculate the products of the basis elements in Theorem \ref{thm reprel}, from which all the products of the basis elements can be calculated.

Those formulas tell us that we can define the \lq generic' centralizer algebra of $\mathbf{CP_k}$ if we put $N=X$ ($X$ is an indeterminate.) as in the usual Brauer centralizer algebras.

From the definition, $\mathbf{CP_k}$ naturally contains the ordinary Brauer centralizer algebra and the symmetric group of degree $k$.  

Since we already know the decomposition rules (see the formula \eqref{eqn evendelnat}, \eqref{eqn evendelnnat}, \eqref{eqn odddelnat}.) of the tensor products, all the irreducible spin representations (that is, the representations not coming from $O(N)$ or $SO(N)$ ) of $Pin(N)$ or $Spin(N)$ occur in this space $\Delta \bigotimes \bigotimes^k V$ for some $k$. (See Section \ref{sect par}.)

So  we define the subspace $T_{k,s}^{0}$ $( s = 1,2, \ldots, \min(k,n))$
of $\Delta \bigotimes \bigotimes^k V$
 to be the intersection of the kernels of all the projections from 
$\Delta \bigotimes \bigotimes^{k} V$ to $\Delta$ (See Definition \ref{def prind}.) and of the contractions and of the alternating operators of degrees greater than $s$. This space is an analog of the traceless tensors in the classical arguments (\cite{weyl}).

Then on this space $T_{k,s}^{0}$, the symmetric group $\frak{S}_k$ of degree $k$ and $Pin(N)$ act as a dual pair and we obtain an analogy of Schur-Weyl reciprocity for $Pin(N)$. (See Theorem \ref{thm dualpin}.)

For odd spinor groups, the irreducible $Pin(2n+1)$ modules are still irreducible as $Spin(2n+1)$-modules, hence $Spin(2n+1)$ and $\frak{S}_k$  act on this space $T_{k,s}^{0}$ as a dual pair.

For $Spin(2n)$, if we define the \lq associator' $A$ of $\End(\Delta)$ by 
\newline
$A = \id \oplus -\id$ on the direct sum $\Delta = \Delta^+ \bigoplus \Delta^-$ and consider the endomorphism $(A \otimes \id) \in \End(\Delta \bigotimes 
\bigotimes^k V)$, then $(A \otimes \id)$ is an involution 
and the action of $(A \otimes \id )$ commutes with those of $\mathbf{CP_k}$ and $Spin(2n)$. (See Lemma \ref{lem commute}.)

The centralizer algebra $\mathbf{CS_k} =  \End_{Spin(2n)}(\Delta \bigotimes 
\bigotimes^k V)$ for $Spin(2n)$ is generated by $\mathbf{CP_k}$  and $(A \otimes \id)$.

We denote the $\pm 1$ eigenspaces of $(A \otimes \id)$ in the spaces 
$T_{k,s}^{0}$ by $T_{k,s}^{0,\pm}$.
On each of the spaces $T_{k,s}^{0,\pm}$, $Spin(2n)$ and $\frak{S}_k$ act as the dual pair. (See Theorem \ref{thm dualspin}.)

The key theorems, which connect the fundamental pin or spin representations to the tensor calculus, are Theorem \ref{thm phiodd}, \ref{thm phieven} and Theorem \ref{thm phipinodd}, \ref{thm phipin}, which give
 the explicit $Pin(N)$-equivariant (or $Spin(N)$-equivariant) embeddings from $\bigwedge^{\ell} V$ to $\Delta \bigotimes \Delta^*$ under the specified basis of both spaces.
The transformation matrices of the same weight spaces in the direct sum of the exterior products and $\Delta \bigotimes \Delta^*$ are given by the Hadamard matrices.

We note that the construction here goes over a field $\mathbb{Q}(\sqrt{2})$.

\section{parametrization of the Irreducible Representations}\label{sect par}
We take a symmetric anti-diagonal matrix $S = (\delta_{i,N+1-i})$ of size $N$, 
where 
$\delta_{i,j}$ is Kronecker's delta 
and  define  the orthogonal group as
$$
O(N,\mathbb{R})= \{ g \in GL(N); g S {}^t\! g = S \}.
$$

Then its Lie algebra is given by
$$
\frak{so}(N,\mathbb{R})= \{ X \in M(N,\mathbb{R}); X S + S {}^t\! X = 0 \}.
$$

As a Cartan subalgebra $\frak{H}$, we can take the intersection of the above Lie 
algebra and the diagonal matrices and let $\be_i$ be a homomorphism from Cartan 
subalgebra to $\mathbb{R}$ defined by taking the $(i,i)$-th component of the diagonal 
matrices. 

Let $P^{+}$ be the set of the dominant integral weights for $\frak{so}(N,\mathbb{R})$.   
Namely for $N=2n+1$, $P^{+}$ is given by
$$
P^{+}=\{\lam_1 \be_1+\lam_2 \be_2+\ldots+\lam_n \be_n; \lam_1 \ge \lam_2 \ge \ldots 
\ge \lam_n \ge 0\}
$$
and for $N=2n$,  $P^{+}$ is given by
$$
P^{+}=\{\lam_1 \be_1+\lam_2 \be_2+\ldots+\lam_n \be_n; \lam_1 \ge \lam_2 \ge \ldots 
\ge \vert \lam_n \vert \ge 0\}.
$$

Here all the $\lam_i$'s are integers or half-integers (namely $1/2 +\Bbb{Z}$) 
simultaneously.

We put $\lam=(\lam_1,\lam_2,\ldots,\lam_n)$ and denote the corresponding 
irreducible character by $\lam_{Spin(N)}$.

The set of the dominant integral weights of the group $SO(N)$ is the subset of the above, defined by the condition that  all the $\lam_i$'s are integers.
So in that case 
 we also write $\lam_{SO(N)}$ for $\lam_{Spin(N)}$.
Then the vector representation (the natural representation) is denoted by 
$(1)_{SO(N)}$ and the representation $\bigwedge^i V$ of \G \, ($i=1,2,\ldots,n$) is 
denoted by $e_i$ as usual.
$e_i$ is the irreducible representation $(1^i)_{SO(N)}$ if $i < n$ or if  $i = n$ and $N=2n+1$.

If  $i = n$ and $N=2n$, 
$e_n$ is not irreducible and is decomposed into the sum of two irreducible 
representations $e_n^+ =(1^n)_{SO(2n)}$ and $e_n^- =(1^{n-1},-1)_{SO(2n)}$.
Later we give an explicit description of $e_n^+$ and $e_n^-$ in $\bigwedge^n V$.

We note that $\bigwedge^i V \cong \bigwedge^{N-i} V$ as $SO(N)$ module.

For a partition $\delta=(\delta_1,\delta_2,\ldots, \delta_n)$, 
we denote the irreducible representations of $Spin(N)$ as follows.

For $Spin(2n+1)$, we put 
 $\lam=(1/2+\delta_1,1/2+\delta_2,\ldots,1/2+\delta_n)$ and 
denote this irreducible representation by 
$\lam_{Spin(2n+1)}=[\Delta,\delta]_{Spin(2n+1)}$.
 Then $[\Delta,\delta]_{Spin(2n+1)}$ can be considered as an irreducible character of $Pin(2n+1)$. 
When we consider this as the representation of $Pin(2n+1)$, we write $[\Delta, \delta]_{Pin(2n+1)}$. (See also Remark \ref{rem dphi}.)

For $Spin(2n)$, we put 
$$
(1/2 + \delta)^+_{Spin(2n)} = 
(1/2+\delta_1,1/2+\delta_2,\ldots,1/2+\delta_n)_{Spin(2n)}
$$
and 
$$
(1/2 + \delta)^-_{Spin(2n)} = (1/2+\delta_1,1/2+\delta_2,\ldots,1/2+\delta_{n-1}, 
-1/2-\delta_n)_{Spin(2n)}.
$$

We denote simply by $\Delta^+$ the irreducible representation 
$(1/2 + \emptyset)^+_{Spin(2n)}$ and by $\Delta^-$ the irreducible representation 
$(1/2+\emptyset)^-_{Spin(2n)}$.

Then we have 
$$
\Delta = \Delta^+ + \Delta^-=\prod_{i=1}^{n} ({x_i}^{1/2} + {x_i}^{-1/2})$$
and
$$
\Delta^+ - \Delta^-=\prod_{i=1}^{n} ({x_i}^{1/2} - {x_i}^{-1/2}).
$$
Here $x_i$ denotes the $i$-th diagonal component of the maximal torus (the diagonal matrices = $\{\diag(x_1,x_2, \ldots,x_n, x_n^{-1}, \ldots,x_1^{-1}) \}$ ) of the group $SO(2n)$.
 We put $\Delta'=\Delta^+ - \Delta^-$.
Moreover 
we introduce the sum character and the difference character for $Spin(2n)$ and denote 
them by
$$
[\Delta, \delta]_{Spin(2n)} =\lam_{Spin(2n)}^{(+)}= (1/2 + \delta )^+_{Spin(2n)} + 
(1/2 + \delta )^-_{Spin(2n)}
$$
 and
$$
[\Delta', \delta]_{Spin(2n)} =\lam_{Spin(2n)}^{(-)}=  (1/2 + \delta )^+_{Spin(2n)} 
- (1/2 + \delta )^-_{Spin(2n)}.
$$

Then $[\Delta, \delta]_{Spin(2n)}$ becomes the irreducible character of $Pin(2n)$ and when we consider this as the representation of $Pin(2n)$, we write $[\Delta, \delta]_{Pin(2n)}$. 
(See the section \ref{sect clif}.)

For a partition $\mu=(\mu_1,\mu_2,\ldots,\mu_n)$, we also denote the sum character 
and the difference character for $\SOe$ by 

$$
\mu_{SO(2n)}^{(+)}= (\mu_1,\mu_2,\ldots,\mu_n)_{SO(2n)} + (\mu_1,\ldots,\mu_{n-
1},-\mu_n)_{SO(2n)}
$$
and 
$$
\mu_{SO(2n)}^{(-)}= (\mu_1,\mu_2,\ldots,\mu_n)_{SO(2n)} - (\mu_1,\ldots,\mu_{n-
1},-\mu_n)_{SO(2n)}.
$$
Then if $\ell(\mu) <n$, $\mu_{SO(2n)}^{(+)}= 2 \mu_{SO(2n)}$ and 
$\mu_{SO(2n)}^{(-)}= 0$. For a partition $\mu$ with $\ell(\mu) <n$, 
$\mu_{SO(2n)}$ becomes an irreducible character of $O(2n)$, so we sometimes write 
$\mu_{O(2n)}$ for $\mu_{SO(2n)}$ in this case.

 If $\ell(\mu) = n$, $\mu_{SO(2n)}^{(+)}$ is the irreducible character of $O(2n)$ and 
we sometimes write $\mu_{O(2n)}$ for $\mu_{SO(2n)}^{(+)}$ in this case.

First we give complementary formulas which are not dealt with in the previous paper \cite{k3}.
We fix some notations. 
For a partition $\lam$ with its length at most $n$, we denote the irreducible 
character of the Symplectic group $\Syp$ by 
 $\lam_{\Syp}$. Then
 ${}^n\! C^{\lam}_{\mu,\nu}$ denotes the multiplicity 
of $\lam_{\Syp}$ in the irreducible decomposition of the tensor product 
$\mu_{\Syp}\nu_{\Syp}$. Namely 
$$
\mu_{\Syp}\nu_{\Syp} = \sum_{\lam}
 {}^n\! C^{\lam}_{\mu,\nu} \lam_{\Syp}.
$$

We note that we can calculate these structure constants  ${}^n\! C^{\lam}_{\mu,\nu}$ explicitly by using the Littlewood-Richardson rules. $($ See 
Theorem 3.1 in \cite{k1} and also \cite{k3},\cite{k4}, \cite{kt1}.$)$

\begin{thm}[Decomposition of Tensor Products for Spin(2n)]\label{thm dtp}

\begin{enumerate}

\item
If $\ell(\mu) = n$, we have 
\begin{equation}
[\Delta,\delta]_{Spin(2n)} \mu_{SO(2n)}^{(-)} 
= \sum_{\substack{\nu, \lam \\ \delta/\nu: vertical\,strip}}
(-1)^{\vert \delta/\nu \vert}
 \; {}^n\! C^{\lam}_{\nu, (\mu-(1^n))}
\sum_{\substack{\kappa \supseteq \lam, \ell(\kappa) \leqq n \\ \kappa/\lam: 
vertical\,strip}}
[\Delta',\kappa]_{Spin(2n)}.
\tag{\ref{thm dtp}.1}
\label{eqn spin+so-}\end{equation}

Here $(\mu-(1^n))$ denotes the partition $(\mu_1-1,\mu_2-1,\ldots,\mu_n-1)$ and $\delta/\nu$ denotes the skew Young diagram and $\vert \delta/\nu \vert$ denotes the size of the skew Young diagram.

\item
If $\ell(\mu) = n$, we have 
\begin{equation}
[\Delta',\delta]_{Spin(2n)} \mu_{SO(2n)}^{(-)} 
= \sum_{\substack{\nu, \lam \\ \delta/\nu: vertical\,strip}}
(-1)^{\vert \mu \vert + \vert \nu \vert}
\; {}^n\! C^{\lam}_{\nu, (\mu-(1^n))}
\sum_{\substack{\kappa \supseteq \lam, \ell(\kappa) \leqq n \\ \kappa/\lam: 
vertical\,strip}} (-1)^{\vert \kappa \vert}
[\Delta,\kappa]_{Spin(2n)}.
\tag{\ref{thm dtp}.2}
\label{eqn spin-so-}\end{equation}

\end{enumerate}

\end{thm}

\begin{proof}

The proof is almost same as in the proof of Theorem 1.2 of \cite{k3}.

From Theorem 10.1 in  \cite{k3}, we have 
\begin{equation*}
\begin{split}
[\Delta,\delta]_{Spin(2n)}(x_1,x_2,\ldots,x_n)
=&\prod_{i=1}^{n} ({x_i}^{1/2} + {x_i}^{-1/2}) 
\sum_{ \delta/\mu: vertical\,strip } (-1)^{\vert \delta/\mu \vert} 
\mu_{Sp(2n)}(x_1,x_2,\ldots,x_n) \\
=&\prod_{i=1}^{n} ({x_i}^{1/2} + {x_i}^{-1/2}) 
(-1)^{\vert \delta \vert} \delta_{SO(2n+1)}(-x_1,-x_2,\ldots,-x_n)
\end{split}
\end{equation*}
and
\begin{equation*}
\begin{split}
[\Delta',\delta]_{Spin(2n)}(x_1,x_2,\ldots,x_n) 
=&\prod_{i=1}^{n} ({x_i}^{1/2} - {x_i}^{-1/2}) 
\sum_{ \delta/\mu: vertical\,strip } \mu_{Sp(2n)}(x_1,x_2,\ldots,x_n) \\
=&\prod_{i=1}^{n} ({x_i}^{1/2} - {x_i}^{-1/2}) 
\delta_{SO(2n+1)}(x_1,x_2,\ldots,x_n)
\end{split}
\end{equation*}
and 
for a partition $\mu$ with  $\ell(\mu)=n$, we have the following.
\begin{equation*}
\mu_{SO(2n)}^{(-)}(x_1,x_2,\ldots,x_n)
=\prod_{i=1}^{n} (x_i - {x_i}^{-1}) (\mu - (1^n))_{Sp(2n)}(x_1,x_2,\ldots,x_n).
\end{equation*}

Here $x_i$ denotes the $(i,i)$-th component of the diagonal matrices which we take 
as a maximal torus for each of the classical groups.

Moreover from Lemma 10.2 in \cite{k3}, we have

\begin{equation*}
\prod_{i=1}^{n} ({x_i}^{1/2} + {x_i}^{-1/2})^2 
\delta_{Sp(2n)}(x_1,x_2,\ldots,x_n)
 = \sum_{\substack{
\kappa/\delta: vertical\,strip \\
 \ell(\kappa) \leqq n }}
\kappa_{SO(2n+1)}(x_1,x_2,\ldots,x_n).
\end{equation*}

If we write the above in the form 

\begin{equation*}
\dfrac{\prod_{i=1}^{n} (1 + {x_i})^2}{x_1 x_2 \cdots x_n} 
\delta_{Sp(2n)}(x_1,x_2,\ldots,x_n)
 = \sum_{\substack{
\kappa/\delta: vertical\,strip \\
 \ell(\kappa) \leqq n }}
\kappa_{SO(2n+1)}(x_1,x_2,\ldots,x_n).
\end{equation*}
 and substitute $-x_i$ for $x_i$, we have 
\begin{equation*}
(-1)^{n + \vert \delta \vert} \prod_{i=1}^{n} ({x_i}^{1/2} - {x_i}^{-1/2})^2 
\delta_{Sp(2n)}(x_1,x_2,\ldots,x_n)
 = \sum_{\substack{
\kappa/\delta: vertical\,strip \\
 \ell(\kappa) \leqq n }}
\kappa_{SO(2n+1)}(-x_1,-x_2,\ldots,-x_n).
\end{equation*}

We note that $\delta_{Sp(2n)}(-x_1,-x_2,\ldots,-x_n) = (-1)^{\vert \delta \vert} 
\delta_{Sp(2n)}(x_1,x_2,\ldots,x_n)$. 

So we have ( for abbreviation, we write $\bx$ for  $x_1,x_2,\ldots,x_n$.)
\begin{equation*}
\begin{split}
& [\Delta,\delta]_{Spin(2n)} \mu_{SO(2n)}^{(-)} \\
& \qquad =\prod_{i=1}^{n} \{({x_i}^{1/2} + {x_i}^{-1/2})(x_i - {x_i}^{-1})\} 
(\mu - (1^n))_{Sp(2n)}(\bx)
\sum_{ \delta/\nu: vertical\,strip } (-1)^{\vert \delta/\nu \vert} 
\nu_{Sp(2n)}(\bx)  \\
& \qquad =\prod_{i=1}^{n} \{({x_i}^{1/2} + {x_i}^{-1/2})(x_i - {x_i}^{-1})\} 
\sum_{ \delta/\nu: vertical\,strip } (-1)^{\vert \delta/\nu \vert}
\:{}^n\! C^{\lam}_{\nu, (\mu-(1^n))}\lam_{Sp(2n)}(\bx)
\end{split}
\end{equation*}

Since ${\ds \prod_{i=1}^{n} ({x_i}^{1/2} + {x_i}^{-1/2})^2 \lam_{Sp(2n)}(\bx) = 
\sum_{\substack{
\kappa/\delta: vertical\,strip \\
 \ell(\kappa) \leqq n }}
\kappa_{SO(2n+1)}(\bx)}$, we have the formula (\ref{thm dtp}.1).

Also we have 
\begin{equation*}
\begin{split}
& [\Delta',\delta]_{Spin(2n)} \mu_{SO(2n)}^{(-)} \\
& \qquad =\prod_{i=1}^{n} \{({x_i}^{1/2} - {x_i}^{-1/2})(x_i - {x_i}^{-1})\} 
(\mu - (1^n))_{Sp(2n)}(\bx)
\sum_{ \delta/\nu: vertical\,strip } \nu_{Sp(2n)}(\bx)  \\
& \qquad =\prod_{i=1}^{n} \{({x_i}^{1/2} + {x_i}^{-1/2})({x_i}^{1/2} - {x_i}^{-
1/2})^2\} \sum_{ \delta/\nu: vertical\,strip } 
\:{}^n\! C^{\lam}_{\nu, (\mu-(1^n))}\lam_{Sp(2n)}(\bx)
\end{split}
\end{equation*}

Since ${\ds (-1)^{n + \vert \lam \vert} \prod_{i=1}^{n} ({x_i}^{1/2} - {x_i}^{-
1/2})^2 \lam_{Sp(2n)}(\bx) = \sum_{\substack{
\kappa/\delta: vertical\,strip \\
 \ell(\kappa) \leqq n }} \kappa_{SO(2n+1)}(-x_1,-x_2,\ldots,-x_n)}$, 
and $\vert \lam \vert  \equiv \vert \nu \vert + \vert (\mu-(1^n)) \vert \mod{2}$,  
we have the formula (\ref{thm dtp}.2).

\end{proof}

From the above and  Theorem 1.2 of \cite{k3}, we can deduce the followings.
For abbreviation, we introduce the following convention.
For $\ep \in \{ \pm 1\}$, if $\ep = 1$, we consider $\Delta^{\ep} = \Delta^+$, $e_n^{\ep} = e_n^+$, 
$(1/2 + \delta)^{\ep}_{Spin(2n)} = (1/2 + \delta)^+_{Spin(2n)}$, etc.  and  if $\ep = -1$, 
$\Delta^{\ep} = \Delta^-$, $e_n^{\ep} = e_n^-$, $(1/2 + \delta)^{\ep}_{Spin(2n)} = (1/2 + \delta)^-_{Spin(2n)}$, etc.
\begin{thm}\label{thm repeven}
For $\ep, \ep_1, \ep_2 \in \{ \pm 1\}$, we have the following formulas.
\begin{enumerate}

\item
\begin{equation}
(\Delta^{\ep})^2 = e_n^{\ep} + \sum_{\substack{0 \leqq 2i \leqq n }} e_{n-2i}.
\tag{\ref{thm repeven}.1}\label{eqn evendeldel}
\end{equation}

\begin{equation}
\Delta^+ \Delta^- =  \sum_{\substack{0 \leqq 2i \leqq n-1 }} e_{n-1-2i}.
\tag{\ref{thm repeven}.2}\label{eqn evendel+del-}
\end{equation}

\item
For a partition  $\delta$ with its length $\ell(\delta) < n$, we have 
\begin{equation}
\begin{split}
&(1/2 + \delta)^{\ep}_{Spin(2n)} (1)_{SO(2n)} \\
& = (1/2 + \delta)^{-\ep}_{Spin(2n)} + 
\sum_{\substack{ \mu \supset \delta \\
 \vert \mu/\delta \vert = 1}} (1/2 + \mu)^{\ep}_{Spin(2n)}
 + \sum_{\substack{ \delta \supset \mu \\
 \vert \delta/\mu \vert = 1}} (1/2 + \mu)^{\ep}_{Spin(2n)}.
\end{split}
\tag{\ref{thm repeven}.3}\label{eqn evendelnat}
\end{equation}

Here $\mu/\delta$ denotes the skew diagram.

If $\ell(\delta) =n$, we have 
\begin{equation}
\begin{split}
&(1/2 + \delta)^{\ep}_{Spin(2n)} (1)_{SO(2n)} 
 = \sum_{\substack{ \mu \supset \delta \\
 \vert \mu/\delta \vert = 1,\\ \ell(\mu) \leqq n}} (1/2 + \mu)^{\ep}_{Spin(2n)}
 + \sum_{\substack{ \delta \supset \mu \\
 \vert \delta/\mu \vert = 1}} (1/2 + \mu)^{\ep}_{Spin(2n)}.
\end{split}
\tag{\ref{thm repeven}.4}\label{eqn evendelnnat}
\end{equation}

\item 
For a partition $\mu$ with $\ell(\mu) < n$, we have 

\begin{equation}
\Delta^{\ep} \mu_{SO(2n)} =   \sum_{\substack{ \mu \supseteq \nu \\
\mu/\nu \, \text{vertical strip} \\ \vert \mu/\nu \vert \equiv 0 \mod{2} }} (1/2 + 
\nu)^{\ep}_{Spin(2n)} + \sum_{\substack{ \mu \supseteq \nu \\
\mu/\nu \, \text{vertical strip} \\ \vert \mu/\nu \vert \equiv 1 \mod{2} }} (1/2 + 
\nu)^{-\ep}_{Spin(2n)}.
\tag{\ref{thm repeven}.5}\label{eqn evendelmu}
\end{equation}

\item

For a partition $\mu=(\mu_1,\mu_2,\ldots,\mu_n)$ with $\ell(\mu) =n$, 
we write $\mu^+_{SO(2n)}=\mu_{SO(2n)}$ and $\mu^-_{SO(2n)}=(\mu_1,\ldots,\mu_{n-
1},-\mu_n)_{SO(2n)}$. Then we have

\begin{equation}
\Delta^{\ep_1} \mu^{\ep_2}_{SO(2n)} =  \sum_{\substack{ \mu \supseteq \nu \\
\mu/\nu \, \text{vertical strip} \\ (-1)^{\vert \mu/\nu \vert} =\ep_1 \ep_2 }}  (1/2 
+ \nu)^{\ep_2}_{Spin(2n)}.
\tag{\ref{thm repeven}.6}\label{eqn evendelmun}
\end{equation}

Particularly, for the exterior products $e_i$'s 
if $i < n$,  we have 
\begin{equation*}
\Delta^{\ep} (1^i)_{SO(2n)} = \sum_{0 \leqq 2s \leqq i} (1/2 + (1^{i-
2s}))^{\ep}_{Spin(2n)} + \sum_{0 \leqq 2s \leqq i-1} (1/2 + (1^{i-1-2s}))^{-
\ep}_{Spin(2n)}\end{equation*}
 and for $e_n^{\ep}$, we have 
\[
\Delta^{\ep}  e_n^{\ep} = \sum_{0 \leqq 2s \leqq n} (1/2 + (1^{n-
2s}))^{\ep}_{Spin(2n)}
\]
and 
\[
\Delta^{-\ep}  e_n^{\ep} = \sum_{0 \leqq 2s \leqq n-1} (1/2 + (1^{n-1-
2s}))^{\ep}_{Spin(2n)}.
\]

\item
$$
(e_n^+)^2 = \sum_{\substack{s \equiv n \mod{2} \\ t \equiv 0 \mod{2} \\ 0 \leqq s 
+ t \leqq n}}
 (2^s, 1^t)_{SO(2n)}, \quad
(e_n^+)(e_n^-) = \sum_{\substack{s \equiv n-1 \mod{2} \\ t \equiv 0 \mod{2} \\ s + 
t \leqq n-1}}
 (2^s,1^t)_{SO(2n)}.
$$

$$
(e_n^-)^2 = (2^{n-1}, -2)_{SO(2n)} + \sum_{\substack{s \equiv n \mod{2} \\ s \leqq 
n-2}} (2^s, 1^{n-s-1},-1)_{SO(2n)} + 
 \sum_{\substack{s \equiv n \mod{2} \\ t \equiv 0 \mod{2} \\ s + t \leqq n-2}}
 (2^s, 1^t)_{SO(2n)}.
$$

\end{enumerate}
\end{thm}

\begin{proof} 
All these formulas can be deduced easily from Theorem \ref{thm dtp} and Theorem 1.2 of \cite{k3}.
We note that the diagram obtained by attaching a vertical strip to $(\mu - (1^n))$ is equal to the diagram obtained by removing the vertical strip from 
$\mu$.
\end{proof}

\begin{rem}
In Example 10.3 of \cite{k3}, the formulas \eqref{eqn evendelnat}, \eqref{eqn evendelnnat} 
 are stated in different forms, since  the author considered there that if we allow 
$\mu_n$ to take the value $-1$, the formulas become slightly of unified form. But this  notation is not compatible with the other formulas in \cite{k3}. So it might cause some misunderstandings.  Also the author must apologize for writing down the incorrect formula in Example 10.3 of \cite{k3},
$$
\Delta' \Delta' = e_n - e_{n-1} + e_{n-2} - \cdots + (-1)^{n-1}e_1 + (-1)^n.
$$

The correct formula deduced from the formula (\ref{eqn evendeldel}), (\ref{eqn 
evendel+del-}) in the Theorem {thm repeven} is 
$$
\Delta' \Delta' = e_n - 2e_{n-1} + 2e_{n-2} - \cdots + (-1)^{n-1}2 e_1 + (-1)^n 2.
$$
The rest of the formulas are correct in Example 10.3 of \cite{k3}.
\end{rem}

For $Spin(2n+1)$,
 we quote here several formulas, which we can 
deduce easily from the formulas in the Theorem 1.2 in \cite{k3} and Theorem 3.1 in \cite{k1}.  

\begin{thm}\label{thm repodd}
We have the following formulas.
\begin{enumerate}

\item
\begin{equation}
\Delta^2 = e_0 + e_1  + e_2 + \ldots + e_n.
\tag{\ref{thm repodd}.1}\label{eqn odddel2}
\end{equation}

\item

\begin{equation}
\begin{split}
&[\Delta,\delta]_{Spin(2n+1)} (1)_{SO(2n+1)} \\
& = [\Delta,\delta]_{Spin(2n+1)} + 
\sum_{\substack{ \mu \supset \delta \\
 \vert \mu/\delta \vert = 1,\, \ell(\mu) \leqq n}} [\Delta,\mu]_{Spin(2n+1)}
 + \sum_{\substack{ \delta \supset \mu \\
 \vert \delta/\mu \vert = 1}} [\Delta,\mu]_{Spin(2n+1)}.
\end{split}
\tag{\ref{thm repodd}.2}\label{eqn odddelnat}
\end{equation}

Here $\ell(\mu)$ denotes the length of $\mu$ and $\mu/\delta$ denotes the skew 
diagram.

\item
For any partition $\lam$ with  $\ell(\lam) \leqq n$,  we have 
\begin{equation}
\Delta \textstyle{\bigotimes}
 \lam_{SO(2n+1)} = \sum\limits_{\substack{ \lam \supseteq \mu \\ \lam/\mu : \text{vertical strip}}} [\Delta, \mu]_{Spin(2n+1)}.
\tag{\ref{thm repodd}.3}\label{eqn odddellam}
\end{equation}
Hence 
$\Delta \bigotimes \lam_{SO(2n+1)}$ contains $\Delta$
 if and only if $\lam =(1^k)$ for some $k$ $( 1 \leqq k \leqq n)$.
In this case, we have 
\begin{equation*}
\Delta \textstyle{\bigotimes} (1^k)_{SO(2n+1)} = \sum\limits_{i=0}^{k} [\Delta,(1^i)]_{Spin(2n+1)}.
\end{equation*}
\end{enumerate}
\end{thm}

\section{Clifford Algebras and Spin groups}\label{sect clif}

In this section we give explicit actions of the Lie group $Pin(N)$ and its Lie algebra $\frak{so}(N \R)$ on the spaces $\Delta$ and its dual $\Delta^*$ and $\bigwedge^i \R^N$ with respect to their specified bases, since later in the section 4, we need to show that the linear embedding from $\bigwedge^i \R^N$ to $\Delta \bigotimes \Delta^*$ given explicitly there is equivariant homomorphism by verifying the actions of generaters of $Pin(N)$ or $Pin(N)$ on both sides.
As for the reference, see \cite{ch}.

We consider the spin groups and their Lie algebra mainly over $\R$.

Let $V$ be a vector space of dimension $N$ ($N=2n+1$ or $N=2n$) over  the real number field $\mathbb{R}$ and we take a basis $< \ua{1}, \ua{2}, \ldots, \ua{n}, (\ua{0}), \ub{n}, \ldots, \ub{1} >$ of $V$. Here $\{1,2,\ldots,n,(0), \overline{n}, \ldots, \overline{1}\}$ is an index set and we consider that for the case $N=2n$, we omit the index $0$ and that for the case $N=2n+1$, we add the index $0$. We introduce an order in the index set such that $ 1 < 2 < \ldots < n \;(< 0) <  \overline{n} < \overline{n-1} < \ldots < \overline{1}$.

Let $B(\quad, \quad)$ be a non-degenerate symmetric 
bilinear form on $V$ defined by
\[
 B(\ua{k}, \ua{\ell}) =0, \qquad 
 B(\ub{k}, \ub{\ell}) =0, \qquad 
 B(\ua{k}, \ub{\ell}) = \dfrac{1}{2} \delta_{k,\ell} (k, \ell > 0).
\]
\[
 (B(\ua{0}, \ua{k}) = 0, \qquad 
 B(\ua{0}, \ub{k}) = 0 \quad (k > 0), \qquad 
 B(\ua{0}, \ua{0}) = 1.) 
\]

In the tensor algebra $\frak{T}(V) = {\sum\limits_{i=0}^{\infty} \bigotimes^{i} V}$, let $\frak{I}_{B}$ be an ideal generated by $< v \otimes v - B(v,v) >_{ v \in V}$ ( here $v$ runs over all the elements in $V$.) and we define the Clifford algebra $Cl(V)$ by the quotient ${\ds Cl(V) = \frak{T}(V) / \frak{I}_{B}}$.

We note that 
$< \ua{k} + \ub{k}, \quad \ua{k} - \ub{k} \quad ( k > 0 ),
\quad (\ua{0}) >$ \quad is an orthogonal basis of $V$ and 
$B(\ua{k} + \ub{k}, \ua{k} + \ub{k}) = 1, \quad B(\ua{k} - \ub{k}, \ua{k} - \ub{k}) = -1, \quad (B(\ua{0}, \ua{0}) = 1)$.

Then $Pin(N)$ and $Spin(N)$ are defined by 
\begin{defn}\label{def spin}
\[
Pin(N,B):=\{ \pm v_1 v_2 \ldots v_r ; v_i \in V,  B(v_i, v_i) = \pm 1 \}
\]
and
\[
Spin(N,B):=\{ \pm v_1 v_2 \ldots v_{2r} ; v_i \in V,  B(v_i, v_i) = \pm 1 \}.\]
\end{defn}

\begin{rem}
Since the quadratic form $B$ is the form of maximal index (\cite{ch}), that is, it has the signature $(n,n)$ for $N=2n$ and $(n+1,n)$ for $N=2n+1$ , $Pin(N,B)$ and $Spin(N,B)$ are not compact groups.
\end{rem}

The Lie algebra $L(Pin(N))=L(Spin(N)) (\cong \frak{so}(N,B)) $ in $Cl(V)$ 
is generated by the degree $2$ elements   $< u_{k} u_{\ell} >$, ($k \ne \ell$).

From now on,  we simply write $Pin(N)$ and $Spin(N)$ for $Pin(N,B)$ and $Spin(N,B)$ because we fix the defining quadratic form $B$ for the Clifford algebras.

We define the subspace $V_k= \R \ua{k} + \R \ub{k}$ ( $k = 1,2,\ldots, n$ )  (and $V_0= \R \ua{0}$)
 of $V$ . Then $B \vert_{V_k}$ is non-degenerate and $Cl(V_k)$ is isomorphic to all the $2 \times 2$ matrices $M(2,\R)$. 
This isomorphism $\phi_2$ is given by $\phi_2(\ua{k})= X$ and $\phi_2(\ub{k}) = Y$ and $\phi(\ua{k}\ub{k}-\ub{k}\ua{k}) = H$, here $\{X, Y, H\}$ is the standard triple in $\frak{sl}(2)$ and is given by
\[
X = 
\begin{pmatrix}
0 & 1 \\
0 & 0 \\
\end{pmatrix}, \quad
Y = 
\begin{pmatrix}
0 & 0 \\
1 & 0 \\
\end{pmatrix}, \quad
H = 
\begin{pmatrix}
1 & 0 \\
0 & -1 \\
\end{pmatrix}.
\]

Since $V = V_1 \bigoplus V_2 \bigoplus \cdots \bigoplus V_n (\bigoplus V_0)$ is the 
orthogonal direct sum, we have the following isomorphism in this case.

\[
Cl(V) \cong Cl(V_1) \textstyle{\bigotimes} Cl(V_2) \textstyle{\bigotimes} \cdots 
Cl(V_n) (\textstyle{\bigotimes} Cl(V_0)). 
\]
Here $Cl(V_k)$ is isomorphic to $M(2,\R)$ for $k >0$ and $Cl(V_0)$ is a two dimensional 
semisimple commutative algebra.

\begin{rem}\label{rem clten}
In \cite{ch} (see p38 Theorem 3.3.), the tensor products of the clifford algebras are considered in the category of graderd algebras and there appear the signatures in the products. But here in the above we consider the tensor products in the usual (not graded ) ways. 
So we need to say about the above isomorphism.

\begin{lem}[\cite{tak}]\label{lem:tpca} 
Let $(V,B)$ be a pair of vector space $V$ over $\R$ 
and a non-degenerate symmetric bilinear form $B$ on $V$.
Let $V_1$, $V_2$ be a vector subspace of $V$ such that $(V,B)$ is the orthogonal direct sum of $(V_1,B \vert_{V_1})$
 and $(V_2,B \vert_{V_2})$.)
 Let $\al$ be the involutive automorphism of $Cl(V_1)$ i.e., $\al(v_1) = -v_1$ for any elements $v_1 \in V_1$.
 
 If there exists an involutive element $\theta \in Cl(V_1)$ such that  
 $Ad(\theta) = \al$, (i.e. $\theta^2 = 1$ and $ \al(x) = \theta x \theta^{-1}$ for any $x \in Cl(V_1)$.)
 then we have
 $$
 Cl(V) \cong Cl(V_1) \otimes Cl(V_2).
 $$
 Here the tensor product is an ordinary tensor product (i.e., not graded one) and the isomorphism $\Phi$ from $Cl(V)$ to $Cl(V_1) \bigotimes Cl(V_2)$ is given by 
\begin{equation*}
\Phi(u)= 
\begin{cases}
u \otimes 1 & \quad \text{if $u \in V_1$} \\
\theta \otimes u & \quad \text{if $u \in V_2$}.
\end{cases}
\end{equation*}

\end{lem}

\begin{proof}

Since both the algebras have the same dimension, 
it is enough to show that $\Phi$ is surjective homomorphism.

For any $v=(v_1, v_2) \in V_1 \bigoplus V_2$ ($v_i \in V_i$),
 $\Phi(v) \Phi(v)=(v_1 \otimes 1 + \theta \otimes v_2)^2 =  
{v_1}^2  \otimes 1 + {\theta}^2 \otimes {v_2}^2 + v_1 \theta \otimes v_2
 + \theta  v_1 \otimes v_2 = (B(v_1, v_1) + B(v_2, v_2)) 1 \otimes 1$,
 since $\al(v_1) = -v_1 = \theta v_1 {\theta}^{-1}$.
From the universality,  $\Phi$ becomes a homomorphism and  $Cl(V_1)$ is included in the image of  $\Phi$, so $\theta \otimes 1$ is in the image.
Hence $\Phi$ is surjective.
\end{proof}

For the orthogonal direct sum $V = V_1 \bigoplus V_2 \bigoplus \cdots \bigoplus V_n (\bigoplus V_0)$, 
the above $\theta$ in the lemma for each direct summand  
is given by 
$
\theta = \ua{k}\ub{k}-\ub{k}\ua{k},\quad \phi(\theta) =H.
$
\end{rem}

From the above,  $V_1 \bigotimes V_2 \bigotimes \cdots \bigotimes V_n (\bigotimes V_0) $ becomes a representation space of $Cl(V)$.

The central primitive idempotents in $Cl(V_0)$ are given by
${\ds \frac{1 \pm \ua{0}}{{2}}}$.
Then $\ua{0}$ acts on ${\ds V_{0}^{\pm}= \R \frac{1 \pm \ua{0}}{{2}}}$ by $\pm 1$.

And so if $\dim V = 2n+1$,  for the Clifford algebra $Cl(V)$, we have the two different irreducible representations $\phi_{\pm}$ according to $V_0^{\pm}$.

Since the dimensions of $V_0^{\pm}$ are one, so we omit the last tensor and 
define the representation space of $\phi_{\pm}$ by the same $ \Delta = V_1 \bigotimes V_2 \bigotimes \cdots \bigotimes V_n$.

If $\dim V = 2n$, $Cl(V)$ is a central simple algebra over $\R$ and the above 
$\Delta$ is the unique irreducible representation of $Cl(V)$ and we denote it by $\phi$.

Explicit representation matrices of  $\phi$, $\phi_{\pm}$ are given in the following.

\begin{lem}\label{lem clactdel}
The representation matrices of $\phi$, $\phi_{\pm}$ on $ \Delta = V_1 \bigotimes V_2 \bigotimes \cdots \bigotimes V_n$ are given as follows. 
\begin{align*}
& \ua{k} \longmapsto & H \otimes \cdots H \otimes \overset{k}{X} \otimes 1_2 \cdots \otimes 1_2, \\
& \ub{k} \longmapsto & H \otimes \cdots H \otimes \overset{k}{Y} \otimes 1_2 \cdots \otimes 1_2, \\
(& \ua{0} \longmapsto & \pm H \otimes H \otimes H \otimes \cdots \otimes H,)  
\end{align*}
 ,where $k \in \{1,2,\ldots,n\}$ and $1_2$ denotes the identity matrix of size 2 and the sign of the image of $\ua{0}$ is $+1$ for $\phi_{+}$ and $-1$ for $\phi_{-}$. 
\end{lem}

If $N=2n$, $\Delta$ is not irreducible as $Spin(2n)$ module and is decomposed into two irreducible representations $\Delta^+$ and $\Delta^-$.
We note that $\Delta$ is irreducible as a $Pin(2n)$ module.

If $N=2n+1$ and $\phi_{\pm}$ are restricted to $Spin(2n+1)$, it is known that they give equivalent irreducible representations of $Spin(2n+1)$. 
As $Pin(2n+1)$ modules, they are different in the tensor product with the linear character $\det$. 
By \lq $\det$ ', we denote the linear representation of $Pin(N)$ induced from the natural homomorphism $Pin(N) \longrightarrow O(N) 
\overset{\det}{\longrightarrow} \{\pm 1\}$.

We fix a base of the representation space 
$\Delta = V_1 \bigotimes V_2 \bigotimes \cdots \bigotimes V_n$.
 From the definition, the set of vectors 
 $u_{p_1} \otimes u_{p_2} \otimes \cdots \otimes u_{p_n}$'s  
(here ${p_1} \in \{1, \bar{1} \}$, ${p_2} \in \{2, \bar{2} \}$, $\ldots$, ${p_n} \in \{n, \bar{n} \}$ ) becomes a basis. 
We pick up the indices $\ur{I} =\{ i_j \vert p_{i_j} = \overline{i_j}\}$ and 
denote the above vector $u_{p_1} \otimes \cdots \otimes u_{p_n}$ by 
$[\ur{i_1,i_2,\ldots,i_r}] =[\ur{I}]$,  where  $\ur{I} = \{i_1,i_2,\ldots,i_r\}$ and $i_1 < i_2 < \ldots < i_r$.

From Lemma \ref{lem clactdel}, we can write down the action $\ua{k}$ and $\ub{k}$ on this base. 
Since $H u_k = u_k$, $H u_k = -u_{\bar{k}}$, $X u_k  = 0$, $X u_{\bar{k}} = u_k$, and $Y u_k = u_{\bar{k}}$, $Y u_{\bar{k}} = 0$ ($H$, $X$, $Y$ is the standard triple in $\frak{sl}(2,\R)$),
 we have the following lemma.

\begin{lem}\label{lem clactdelbase}
\[
\ua{k} [\ur{i_1, i_2, i_3, \ldots, i_r}] = 
\begin{cases}
(-1)^{s-1} [\ur{i_1, \ldots, \widehat{i_s}, \ldots, i_r}] \quad & \text{if $k = i_s$  for some $\ur{i_s}$} \\
0 \quad & \text{otherwise}, \\
\end{cases}
\]
 where $\widehat{\ur{i_s}}$ denotes the removal of the index $i_s$  from the index set and

\begin{equation*}
\ub{k} [\ur{i_1, i_2, \ldots, i_r}] = 
 \begin{cases}
(-1)^{s-1} [\ur{i_1,\ldots, i_{s-1}, k, i_s,\ldots, i_r}] \quad & \text{if $k \notin \ur{I}$ and $i_{s-1} < k < i_s$} \\
0 \quad & \text{otherwise}. \\
\end{cases}
\end{equation*}

\[
(\ua{0} [\ur{i_1, i_2, i_3, \ldots, i_r}] = \pm (-1)^r [\ur{i_1, i_2, i_3, \ldots, i_r}] \quad \text{according to $\phi_{\pm}$}.)
\]
\end{lem}

This representation is sometimes called the Fock space representation since it is realized as the action on the exterior algebra of the isotropic subspace of 
$V$ generated by $\ub{1}$, $\ub{2}$,\ldots $\ub{n}$. 
That is, under the correspondence of $[\ur{i_1, i_2, i_3, \ldots, i_r}] 
=u_{\overline{i_1}} \wedge u_{\overline{i_2}} \wedge \ldots \wedge 
u_{\overline{i_r}}$, 
$\ua{k}$ acts on this space by the interior product (here we must adjust the inner product such as $( \ua{k}, \ub{\ell} ) = \delta_{k,\ell}$.) and  $\ub{k}$ acts on this space  by the exterior product $\ub{k} \wedge$ and $\ua{0}$ acts by the degree operator $(-1)^{\deg}$.

Therefore for later convenience, we introduce the following convention.
For the indices in the bracket $[\ur{i_1, i_2, i_3, \ldots, i_r}]$ which are not necessarily in increasing order,  we correspond this element to $0$ if the indices are not distinct, otherwise we correspond this to the signed base element $ \sign(\sigma) [\ur{i_{\sigma(1)}, i_{\sigma(2)}, i_{\sigma(3)}, \ldots, i_{\sigma(r)}}]$, where $\sigma$ is the permutation of $\{1,2, \ldots,r\}$ such that $i_{\sigma(1)} < i_{\sigma(2)} < i_{\sigma(3)} < \ldots < i_{\sigma(r)}$
 and  $\sign(\sigma)$ denotes the signature of $\sigma$.
 This convention is compatible with the base of the Fock space representation.

The explicit action of the Lie algebra $L(Spin(N,B))$ on $V$ is given by
\[ 
\ad(uv) w = [uv, w] = uvw - wuv = 2 B(w,v)u - 2 B(w,u)v,
\]
where $u,v,w \in V$ and the product here is that of the Clifford algebra,
 therefore we can calculate the explicit actions of the standard root vectors of $\G$ with respect to the base $< \ua{1}, \ldots, \ua{n}, (\ua{0}), \ub{n}, 
\ldots, \ub{1} >$ of $V$.
 The standard root vectors of $L(Spin(N))$ for the simple roots are given by $X_k = \ua{k}\ub{k+1}$, $Y_k = \ua{k+1}\ub{k}$ $(1 \leqq k \leqq n-1)$ and $h_i =  \dfrac{\ua{i}\ub{i}-\ub{i}\ua{i}}{2}$ $(1 \leqq i \leqq n)$ and for $L(Spin(2n+1))$, $X_{n} = \ua{n}\ua{0}$, $Y_{n} = \ua{0}\ub{n}$ and for $L(Spin(2n))$, $X_{n} = \ua{n-1}\ua{n}$, $Y_{n} = \ub{n}\ub{n-1}$.
 Therefore the actions of the above root vectors on $V$ and $\Delta$ are given as follows.

\begin{lem}\label{lem spinactvdel}
For $k$ $(1 \leqq k \leqq n-1)$ and $i$ $(1 \leqq i \leqq n)$, we have 
\[
\begin{array}{ll}
\ad(X_k) =  E_{k,k+1}-E_{\overline{k+1},\overline{k}}, & 
\dphi_{\pm} (X_k) =  - 1_2  \otimes \cdots 1_2 \otimes \overset{k}{X} \otimes  Y \otimes 1_2 \cdots \otimes 1_2, \\
\ad(Y_k) = E_{k+1,k}-E_{\overline{k},\overline{k+1}}, & 
\dphi_{\pm} (Y_k) =  - 1_2  \otimes \cdots 1_2 \otimes \overset{k}{Y} \otimes  X \otimes 1_2 \cdots \otimes 1_2, \\
\ad(h_i) = E_{i,i}-E_{\overline{i},\overline{i}}, 
 & \dphi_{\pm} (h_i) = \dfrac{1}{2} 1_2  \otimes 
\cdots 1_2  \otimes \overset{i}{H} \otimes 1_2 \cdots \otimes 1_2. \\ 
\end{array} 
\]

And for $L(Spin(2n+1))$, we have
\[
\begin{array}{ll}
\ad(X_n) =  2 E_{n,0}-E_{0,\overline{n}}, & \dphi_{\pm} (X_n) =  \mp 1_2 \otimes 1_2 \otimes \cdots \otimes 1_2 \otimes X,  \\
\ad(Y_n) = E_{0,n} - 2 E_{\overline{n},0}, & \dphi_{\pm} (Y_n) = \mp 1_2 \otimes 1_2 \otimes \cdots \otimes 1_2 \otimes Y.  \\
\end{array}
\]

And for $L(Spin(2n))$, we have 
\[
\begin{array}{ll}
\ad(X_n) =   E_{n-1,\overline{n}}-E_{n,\overline{n-1}}, & \dphi(X_n) = - 1_2 \otimes \cdots \otimes 1_2 \otimes X \otimes X,  \\
\ad(Y_n) = E_{\overline{n},n-1} -  E_{\overline{n-1},n},  & \dphi (Y_n) = - 1_2 \otimes \cdots \otimes 1_2 \otimes Y \otimes Y.  \\
\end{array}
\]
\end{lem}

Then for any $ 1 \leqq k \leqq n-1$, we have 
$[X_k, Y_k] = h_k - h_{k+1} = H_{\al_k}$ and for $L(Spin(2n+1))$,
$[X_n, Y_n] = 2 h_n = 2 H_{\al_n}$ and for $L(Spin(2n))$, 
$[X_n, Y_n] = h_{n-1} + h_{n}= H_{\al_n}$.

For $Spin(2n)$,
$\ad$ is a Lie algebra isomorphism from $L(Spin(2n))$ to $\frak{so}(V,B)$.
Since $B=\dfrac{1}{2} S$ with respect to the fixed basis  ($S$ is the 
defining symmetric form of $\frak{so}(N,\mathbb{R})$ in the section \ref{sect 
par}.) , $\frak{so}(V,B)$ is identified with the Lie algebra $\frak{so}({\R}^{2n}, S)$.

So from now on, we consider $V$ has the standard basis 
$<\ua{i}, \ub{i} >$ such that the defining symmetric bilinear form is given by $S$. Namely $V$ has the inner product 
$(\quad,\quad)$ such that $(\ua{k}, \ub{\ell}) = \delta_{k,\ell}$ and the values of the inner products between the other base elements are all 0.

For $Spin(2n+1)$, we need to change basis.
Let us put $\uo = \dfrac{1}{\sqrt{2}} \ua{0}$. 
Then $B(\uo, \uo)=  \dfrac{1}{2}$ and 
we take a new basis $< \ua{i}, \uo, \ub{i} >$ of $V$. The representation matrices of $X_n$ and $Y_n$ under \lq $\ad$' change into 
$\ad(\ua{n}\ua{0}) = \ad(X_n) = \sqrt{2}( E_{n,0'}-E_{0',\overline{n}})$ and 
$\ad(\ua{0}\ub{n}) = \ad(Y_n) = \sqrt{2}(E_{0',n} - E_{\overline{n},0'})$. 
Then $\ad$ is a Lie algebra isomorphism from $L(Spin(2n+1))$ to $\frak{so}(V,B)$ and the symmetric matrix corresponding to $B$ on the new base 
is given by $\dfrac{1}{2} S$.  

 So similarly as above,
 we consider $V$ has the standard basis 
$<\ua{i}, \uo, \ub{i} >$ such that the 
defining symmetric bilinear form is given by $S$.  
 Namely $V$ has the inner product $(\quad,\quad)$ such that $(\ua{k}, \ub{\ell}) = \delta_{k,\ell}$, $(\uo \uo) = 1$ and the other inner products between these basis elements are all 0.

From Lemma \ref{lem spinactvdel},
We can write down the action of $X_i$'s and $Y_i$'s explicitly on the 
fundamental spin representation space $\Delta$ (the dual of $\Delta$) of $\frak{so}(N)$ with respect to the basis $[\ur{i_1,i_2, \ldots,i_r}]$'s.

\begin{lem}\label{lem spinactdel}
In this lemma,  $\ur{I}$ denotes the index set $\{i_1,i_2,\ldots, i_r\}$, where  $i_1, i_2, i_3, \ldots, i_r$ are in the increasing order.
If $1 \leqq k \leqq n-1$, we can write down the actions of $X_k$'s and $Y_k$'s simultaneously for $\frak{so}(N)$ with respect to the basis $[\ur{i_1,i_2, \ldots,i_r}]$'s as follows.

\[
X_{k} [\ur{I}] = \begin{cases}
- [\ur{i_1, \ldots, i_{s-1}, k+1, i_{s+1},  \ldots, i_r}] \quad & \text{if $k = i_s$  and 
$ k+1 < i_{s+1}$} \\
0 \quad & \text{otherwise},
\end{cases}
\]
and
\[
Y_{k} [\ur{I}] = 
\begin{cases}
- [\ur{i_1, \ldots, i_{s-1}, k, i_{s+1},  \ldots, i_r}] \quad & \text{if $k +1 = i_s$  and 
$ k > i_{s-1}$} \\
0 \quad & \text{otherwise},
\end{cases}
\]
and 
$
h_{i} [\ur{I}] = - \dfrac{1}{2} [\ur{I}]
$   
if $i \in \ur{I}$ and \vspace*{2mm}
$
h_{i} [\ur{I}] =\dfrac{1}{2} [\ur{I}]
$ 
if $i \notin \ur{I}$.

For $L(Spin(2n+1))$, we have 
\[
\dphi_{\pm}(X_{n}) [\ur{i_1, i_2,  \ldots, i_r}] = 
\begin{cases}
 \mp [\ur{i_1, i_2, \ldots, i_{r-1}}] \quad & \text{if $i_r=n$ } \\
0 \quad & \text{otherwise}, \\
\end{cases}
\]
and
\[
\dphi_{\pm}(Y_{n}) [\ur{i_1, i_2, \ldots, i_r}] = 
\begin{cases}
 \mp [\ur{i_1, i_2, \ldots, i_{r},n}] \quad & \text{if $i_r \ne n$ } \\
0 \quad & \text{otherwise}. \\
\end{cases}
\]

For $Spin(2n)$, we have 
\[
X_{n} [\ur{i_1, i_2, \ldots, i_r}] = 
\begin{cases}
 - [\ur{i_1, i_2, \ldots, i_{r-2}}] \quad & \text{if $i_{r-1}=n-1$  and $i_r=n$} \\
0 \quad & \text{otherwise}, \\
\end{cases}
\]
and 
\[
Y_{n} [\ur{i_1, i_2, \ldots, i_r}] = 
\begin{cases}
 - [\ur{i_1, i_2, \ldots, i_{r},n-1,n}] \quad & \text{if $i_r < n-1$ } \\
0 \quad & \text{otherwise} \\
\end{cases}
\]
\end{lem}

\begin{rem}\label{rem dphi}
As noted before, for $Spin(2n+1)$,  $\phi_{\pm}$ are equivalent if we restrict them to the group $Spin(2n+1)$ and the derivations $\dphi_{\pm}$ only differ in the sign of the representation matrices of $X_n$ and $Y_n$, where $X_n$ and $Y_n$ are fixed short root vectors. 

As $\frak{so}(2n+1)$ modules, if we take a new base $\{ [\ur{i_1,i_2,\ldots, i_r}]_{-} = (-1)^r [\ur{i_1, i_2, \ldots, i_r}]\}$ in the representation space of $\dphi_{-}$, then the representation matrices of $X_i$, $Y_i$, $h_i$ under $\dphi_{-}$ with respect to this new base coincide with those of them under $\dphi_{+}$.

So hereafter, in the representation space of $\phi_{-}$ for $Pin(2n+1)$ (and in  its dual space), we only consider the elements $[\ur{i_1, \ldots, i_r}]_{-}$ instead of $[\ur{i_1, \ldots, i_r}]$. 

And from now on, we change the definition of the notation $[\ur{i_1, i_2, \ldots, i_r}]$ in the representation space $\phi_{-}$ and by $[\ur{i_1, i_2, \ldots, i_r}]$, we denote the element $[\ur{i_1, i_2, \ldots, i_r}]_{-}$ in $\phi_{-}$.

For $\dim V = 2n+1$, the center $Z(Cl(V))$ of the Clifford algebra $Cl(V)$ is given by 
$Z(Cl(V)) = \R \bigoplus \R z$, where  
 $
 z = (\ua{1} - \ub{1})(\ua{1} + \ub{1}) \cdots (\ua{n} - \ub{n})(\ua{n} + \ub{n}) \ua{0}.
$
 We have $z^2 = 1$ and $z$ acts on $V$ by the scalar $-1$ 
 and acts on the representation spaces $\phi_{\pm}$ by the scalar $\pm 1$.

The Clifford theorem tells us that there are exactly two non-equivalent irreducible representations of $Pin(2n+1)$ with the same character and that one is obtained from the other by the tensor product with the linear character $\det$. 

Since $z \notin Spin(2n+1)$, we can distinguish them by the action of $z$.
By $[\Delta,\delta]_{Pin(2n+1), \pm}$, we denote the irreducible representation of $Pin(2n+1)$ with its character $[\Delta,\delta]_{Spin(2n+1)}$,  on which $z$ acts by the scalar $\pm 1$.

For the representation space of $\phi_{\pm}$, we simply denote them by 
$\Delta_{\pm}$ instead of $\Delta_{Pin(2n+1), \pm}$. 
We have $(\Delta_{\pm})^* \cong \Delta_{\pm}$.

In this paper, for abbreviation, $\Delta$ stands for $\Delta_{+}$ or $\Delta_{-}$ with the 
base $\{ [\ur{i_1, i_2, \ldots, i_r}]\}$ noted above, if not specified.
\end{rem}

For $\frak{so}(2n)$,
from Lemma \ref{lem spinactdel}, the parity of the number of the indices in the base element $[\ur{i_1, i_2, i_3, \ldots, i_r}]$ under the above action is preserved.  We define $\Delta^+$ to be the linear space spanned by all the base elements with even parity of the number 
of indices and $\Delta^-$ to be the linear space spanned by all the base elements 
with odd parity of indices.
Then the highest weight vector with the highest weight $\dfrac{1}{2}( \be_1 + \be_2 + \cdots + \be_n)$ of $\Delta^+$ is given by $[\ur{\emptyset}]$ and 
the highest weight vector with the highest weight $\dfrac{1}{2}( \be_1 + \cdots  + \be_{n-1}- \be_n)$ of $\Delta^-$ is given by $[\ur{n}]$.
As for the dual representation  $(\Delta^+)^*$ and $(\Delta^-)^*$, if $ n \equiv 0 \mod{2}$, we have $(\Delta^+)^* \cong \Delta^+ $ and $(\Delta^-)^* \cong \Delta^- $. If $ n \equiv 1 \mod{2}$, we have $(\Delta^+)^* \cong \Delta^- $ and $(\Delta^-)^* \cong \Delta^+ $.

As a $Pin(2n)$ module, the action of the element $\ua{k} \pm \ub{k}$ in $Pin(2n)$ on $\Delta$ is given by
\[
\phi(\ua{k} \pm \ub{k}) = H \otimes \otimes \cdots H \otimes \overset{k}{(X \pm Y)} 
\otimes 1_2 \otimes \cdots \otimes 1_2
\]

The adjoint action of the element $\ua{n}-\ub{n}$ on the Lie algebra $L(Spin(2n))$ in the Clifford algebra $Cl(V)$ gives the standard outer automorphism of $Spin(2n)$ and the sum characters $[\Delta, \delta]_{Spin(2n)}$ are the irreducible characters for $Pin(2n)$.

Also the element $\ua{n}-\ub{n}$ acts on $\Delta$ by  $(\ua{n}-\ub{n}) \Delta^{\pm} = \Delta^{\mp}$ and $(\ua{n}-\ub{n})^2 = -\id$.
  Then the Clifford theorem tells us that as $Pin(2n)$ modules, $(\Delta \bigotimes \det) \cong \Delta$. The associator $A$ 
which realizes this isomorphism is ,by definition, the linear isomorphism on $\Delta =\Delta^{+} \bigoplus \Delta^{-}$  given by
\[
A = \id_{\Delta^{+}} \oplus  - \id_{\Delta^{-}}.
\]

For, it is easy to see that $(\ua{n}-\ub{n}) A = - A (\ua{n}-\ub{n})$ and $A$ commutes with the action of $Spin(2n)$. 
This fact tells us that if two irreducible representations $\lam$ and 
$\widetilde{\lam}$ of $O(2n)$ are associated, i.e., $\lam \cong \widetilde{\lam} \bigotimes \det $, we have the equivalent representation $\Delta \bigotimes \lam \cong \Delta \bigotimes \widetilde{\lam}$ of $Pin(2n)$.

$A$ acts on $\Delta$ by the degree operator 
$A [\ur{I}] = (-1)^{\vert \ur{I} \vert}  [\ur{I}]$.

Also by the same reasoning, all the $[\Delta, \delta]_{Pin(2n)}$ become self-associated since $[\Delta, \delta]_{Pin(2n)}$ is decomposed into two non-equivariant representations under $Spin(2n)$ and the Clifford theorem tells us that it occurs if and only if $[\Delta, \delta]_{Pin(2n)}$ is self-associated. (See p161 in \cite{weyl}.) 
So we have $[\Delta, \delta]_{Pin(2n)} \cong [\Delta, \delta]_{Pin(2n)} \bigotimes \det$. 
This fact means that the characters $[\Delta, \delta]_{Spin(2n)}$ uniquely determine the irreducible representations of $Pin(2n)$.

Since $[\Delta, \delta] = (1/2 + \delta )^+_{Spin(2n)} \bigoplus (1/2 + \delta )^-_{Spin(2n)}$ and $(\ua{n}-\ub{n})(1/2 + \delta )^{\pm}_{Spin(2n)}=(1/2 + \delta )^{\mp}_{Spin(2n)}$, the associator $A$ is given by 
\[
A = \id_{(1/2 + \delta )^+_{Spin(2n)}} \oplus - \id_{(1/2 + \delta )^-_{Spin(2n)}}.
\]
 We write down the above as a lemma and use this fact later.

\begin{lem}\label{lem assoc}
We have the isomorphism  $[\Delta, \delta]_{Pin(2n)} \cong [\Delta, \delta]_{Pin(2n)} \bigotimes \det$. The associator $A$ which realizes this $Pin(2n)$ isomorphism is given by 
\[
A = \id_{(1/2 + \delta )^+_{Spin(2n)}} \oplus  - \id_{(1/2 + 
\delta )^-_{Spin(2n)}}.
\]
Hence the characters $[\Delta, \delta]_{Spin(2n)}$ uniquely determine the irreducible representations of $Pin(2n)$.
Since  $(\Delta \bigotimes \det) \cong \Delta$, we have $\Delta \bigotimes (\lam \bigotimes \det) \cong \Delta \bigotimes \lam$.
\end{lem}

A compact real form $\G_{cpt}$ of $L(Spin(N)) \bigotimes \C \cong \GC$ is generated by the elements  $\sqrt{-1}h_i$ and $\sqrt{-1}(X_i+Y_i)$, $X_i-Y_i$, \quad 
$(i=1,2,\ldots,n)$  as a real Lie algebra.

Then we can introduce the invariant hermitian inner products with respect to 
$\G_{cpt}$ in the space $V$ and $\Delta^{(\pm)}$  respectively. 
Such invariant hermitian forms are given explicitly as follows.
If we define a hermitian inner product in $V$ such that the basis $<\ua{i}, (\uo),  \ub{i} >$ are orthonormal basis, then we can show easily that the hermitian inner product becomes invariant under the action of $\G_{cpt}$. 
Also if we define a hermitian inner product on $\Delta$ such that the basis 
$\{[\ur{i_1, i_2, i_3, \ldots, i_r}]\}$ ($\{i_1, i_2, \ldots, i_r\}$ runs over all the subsets of $[1,n] =\{1, 2, \ldots, n\}$.)  is orthonormal basis,
then we can show easily that the hermitian inner product becomes invariant under the action of $\G_{cpt}$. 
So we consider $V$, $\Delta$, $\Delta^{\pm}$ are the hermitian metric spaces with these invariant hermitian metrics.

Let us consider the dual action of $\Delta$ and $V$.
Let $\Delta^*$ and $V^*$ be the dual vector space of $\Delta$ and $V$ respectively and let $< \ua{i}^*, (\uo^*), \ub{i}^* >$ and $\{[\ur{i_1, \ldots, i_r}]^*\}$ be the dual basis of $< \ua{i}, (\uo), \ub{i} >$ and $\{[\ur{i_1, \ldots, i_r}]\}$ respectively. 
Since the representations matrices of the dual representations on the 
dual basis are given by the minus sign of the transposed matrices, we can write down the dual actions explicitly as in Lemma \ref{lem spinactvdel} and Lemma \ref{lem spinactdel}.

We fix some notations.

In the exterior product $\bigwedge^r V$, we have the natural basis $\{ u_{i_1} \wedge u_{i_2} \wedge \ldots \wedge u_{i_r} \}$, where $i_k \in \{ 1, \ldots, n, (0'), \overline{n}, \ldots,\overline{1} \}$.
We denote these basis elements only by their indices
$$
<\ur{i_1,i_2,\dots,i_r}> =\dfrac{1}{r!}\sum_{\sigma \in \frak{S}_r} \sign(\sigma) 
u_{i_{\sigma^{-1}(1)}} \otimes u_{i_{\sigma^{-1}(2)}} \otimes \ldots \otimes 
u_{i_{\sigma^{-1}(r)}}.
$$

We define an order in the index set by $ 1 < 2 < \ldots< n < (0') < \overline{n} < \ldots < \overline{1} $ and introduce the following convention.
For any $\tI, \tJ, \tW \subseteq [1,n]=\{1,2,\ldots,n\}$,
 $<\ur{I}, \ur{W}, \ulb{W}, \ulb{J}>$ denotes the exterior product $u_{\ur{I}} 
\wedge u_{\ur{W}} \wedge u_{\ulb{W}} \wedge u_{\ulb{J}}$, where for any subset 
$\ur{K}$ of the index set,  $u_{\ur{K}}$ denotes the exterior product 
$<\ur{k_1,k_2,\dots,k_r}>$ for $\mathtt{K}=\{ k_1, k_2, \ldots, k_r \}$  and the indices are arranged in the increasing order defined above, namely $k_1 < k_2 < \ldots < k_r$. 
Also we put $<\ur{I}, \ur{W}, \ur{0'}, \ulb{W}, \ulb{J}> = u_{\ur{I}} \wedge u_{\ur{W}} \wedge u_{0'} \wedge u_{\ulb{W}} \wedge u_{\ulb{J}}$.

For example, for subsets $\tI =\{1, 4\}$ and $\tW =\{2, 5\}$ and $\tJ =\{3, 6\}$ of $[1,7]$,  $<\ur{I}, \ur{W}, \ulb{W}, \ulb{J}>$ denotes the exterior product $u_1 \wedge u_4 \wedge u_2 \wedge u_5 \wedge u_{\overline{5}} \wedge u_{\overline{2}} \wedge u_{\overline{6}} \wedge u_{\overline{3}}$.

Also for the base $\{ [\ur{i_1, i_2, \ldots, i_r}] \}$ of $\Delta$, we introduce a similar notation.
Namely for mutually disjoint subsets $\tI$ and $\tK$ of $[1,n]$, $\{[\ur{I}, 
\ur{K}]$ denotes the element 
$[\ur{i_1, \ldots, i_r, k_1, \ldots, k_s}]$, where 
the indices $i_1$, $i_2$, $\ldots$, $i_r$ are in the increasing order and the same for $\tK$.
We note that the indices in the bracket have the alternating property.
For the dual space $\Delta^*$, 
We introduce the same notation $[\ur{J}, \ur{K}]^*$, conveying the same meaning.

 Then  
$\{[\ur{I}, \ur{K}] \otimes [\ur{J}, \ur{K}]^*\}$ becomes a base of $\Delta 
\bigotimes \Delta^*$ if $\tI$ and $\tK$ and $\tJ$ run over all the mutually disjoint subsets of $[1,n]$.

Since $V$ is a hermitian space, so the exterior products $\bigwedge^k V$ of $V$ are naturally endowed with the hermitian inner product. Namely 
$\{ u_{i_1} \wedge u_{i_2} \wedge \ldots \wedge u_{i_r}\}$ becomes an orthonormal basis of $\bigwedge^k V$.

 Since $V$ has the $O(N)$-invariant bilinear form S, we can identify the vector space $V$ and its dual $V^*$.
 
 Let $\iota$ be the $\R$-linear map from $V$ to $V^*$ defined by
 $\iota (u_k) = {u_{\overline{k}}}^*$.
 Here $k \in \{1,\ldots,n,(0'), \overline{n},\ldots, \overline{1} \}$ and we consider $\overline{\overline{k}} =k$, $\overline{0'}=0'$.
 Then $\iota : V \longrightarrow V^*$ is an $O(N)$-equivariant isomorphism.
We also denote the $O(N)$-equivariant isomorphism from $V^*$ to $V$ by
 the same notation $\iota$, i.e., $\iota (u_k^*) = u_{\overline{k}}$.

The dual base is given by ${<\ur{I}, \ur{W}, (\ur{0'}), \ulb{W}, \ulb{J} >}^* = k! <\ur{I}^*, \ur{W}^*, ({\ur{0'}}^*), \ulb{W}^*, \ulb{J}^* >$   (here $k$ is the degree of the exterior product and $*$ denotes the element obtained by replacing each tensor component by its dual base).

\begin{defn}\label{def iota}
Let $\iota= (-1)^{\binom{k}{2}} \bigwedge^{k} \iota :  \bigwedge^{k} V 
\longrightarrow (\bigwedge^{k} V)^*$ be the $O(N)$-equivariant isomorphism  defined 
by 
$$\iota (<\ur{I}, \ur{W},(\ur{0'}), \ulb{W}, \ulb{J} >) =\dfrac{1}{ k!} {<\ur{J}, \ur{W}, 
({\ur{0'}}), \ulb{W}, \ulb{I} >}^*.$$
\end{defn}

We also denote the inverse of the above by the same letter $\iota$, that is,  
$\iota$ is also the map 
$\iota = (-1)^{\binom{k}{2}} \bigwedge^{k} \iota :  (\bigwedge^{k} V)^* 
\longrightarrow \bigwedge^{k} V$ defined by 
$\iota ({<\ur{I}, \ur{W},(\ur{0'}), \ulb{W}, \ulb{J} >}^*) = k! <\ur{J}, \ur{W}, (\ur{0'}), \ulb{W}, \ulb{I} >$.

To prove the next lemma, we also introduce the linear map $\iota' : \bigwedge^{\ell} V \longrightarrow
\bigwedge^{\ell} V$.

\begin{defn}\label{def iota'ext}
$\iota'(<\ur{J}, \ur{W},(\ur{0'}), \ulb{W}, \ulb{I} >) =  <\ur{I}, \ur{W},(\ur{0'}), \ulb{W}, \ulb{J} >$. 
\end{defn}

Then we can check easily that $ \iota' X_k = -Y_k \iota'$, $ \iota' Y_k = -X_k \iota'$ and $ \iota' h_k = -h_k \iota'$ on $\bigwedge^{\ell} V$.

\begin{lem}\label{lem rl}
The $\frak{so}(2n+1)$-isomorphism $\textstyle{
r_{\ell} : \bigwedge\limits^{\ell} V \longrightarrow \bigwedge\limits^{2n+1-\ell} V}
$ 
is given by follows.

$$
r_{\ell}(<\ur{J}, \ur{W}, \ulb{W}, \ulb{I} >) = (-1)^{\vert \tW \vert} 
<\ur{J}, \ur{W^c}, \ur{0'}, \ulb{W^c}, \ulb{I}>
$$
and 
$$
r_{\ell}(<\ur{J}, \ur{W}, \ur{0'}, \ulb{W}, \ulb{I} >) = (-1)^{\vert \mathtt{W^c} \vert} 
<\ur{J}, \ur{W^c}, \ulb{W^c}, \ulb{I}>,
$$
where $\tW \sqcup \tWc = [1,n] - \tJ -\tI$.

The $\frak{so}(2n)$-isomorphism 
$
r_{\ell} : \textstyle{\bigwedge\limits^{\ell} V \longrightarrow 
\bigwedge\limits^{2n-\ell} V}
$ 
is given by follows.

$$
r_{\ell}(<\ur{J}, \ur{W}, \ulb{W}, \ulb{I} >) = (-1)^{\vert \tI \vert} 
<\ur{J}, \ur{W^c},  \ulb{W^c}, \ulb{I}>,
$$
 where $\tW \sqcup \tWc = [1,n] - \tJ -\tI$.
\end{lem}

\begin{proof}

For $\frak{so}(2n+1)$, we note that we can deduce the above easily from the Theorem \ref{thm phiodd}, but we also can check it directly as in the following argument. So we give proof for $\frak{so}(2n)$ here.

We show that $r_{\ell}$ is  $\frak{so}(2n)$-equivariant.

Then we have $\iota' \circ r_{\ell} = (-1)^{\ell} r_{\ell} \circ \iota'$, since $i + j + 2 w = \ell$.

From the argument just before this lemma, it is enough to show that all the 
$Y_k$ commute with  $r_{\ell}$, since $r_{\ell} X_k = - r_{\ell} \iota' Y_k \iota' = (-1)^{\ell+1}  \iota' r_{\ell} Y_k \iota' = (-1)^{\ell+1}  \iota' Y_k r_{\ell} \iota' = -  \iota' Y_k \iota' r_{\ell} = X_k r_{\ell}$.

Let us prove $r_{\ell} Y_k (<\ur{J}, \ur{W}, \ulb{W}, \ulb{I} >)= Y_k 
r_{\ell}( <\ur{J}, \ur{W}, \ulb{W}, \ulb{I} >)$ for $k = 1,2,\ldots,n-1$.

Here the action of $\frak{so}(2n)$ is the adjoint action and is given by $\ad(Y_k) = E_{k+1,k}-E_{\overline{k},\overline{k+1}}$.

Case 1.  $k, k+1 \in \tJ$ or $k, k+1 \in \tI$.

In this case, both sides are zero and the equality holds.

Case 2. $k \in \tJ$ and $k+1 \in \tI$.

Then we have $k, k+1 \notin \tW$ and 
\begin{align*}
&Y_k (<\ur{J}, \ur{W}, \ulb{W}, \ulb{I} >) 
= \ep <\ur{J-\{k\}}, \ur{W +\{k+1\}}, \ulb{W +\{k+1\}}, \ulb{I-\{k+1\}} > \\
& \qquad \qquad - \ep <\ur{J-\{k\}}, \ur{W +\{k\}}, \ulb{W +\{k\}}, \ulb{I-\{k+1\}} > 
\end{align*}
, where $\ep = (-1)^{\vert \{ j \in \tJ : j > k+1\} \vert + \vert \{ i \in \tI : i > k +1 \} \vert}$.

So we have 
\begin{equation*}
\begin{split}
& r_{\ell} Y_k (<\ur{J}, \ur{W}, \ulb{W}, \ulb{I} >) = 
\ep (-1)^{\vert I \vert -1} <\ur{J-\{k\}}, \ur{W^c +\{k\}}, \ulb{W^c +\{k\}}, 
\ulb{I-\{k+1\}} > \\
& \qquad \qquad - \ep (-1)^{\vert I \vert -1} <\ur{J-\{k\}}, \ur{W^c +\{k+1\}}, \ulb{W^c +\{k+1\}}, \ulb{I-\{k+1\}} > . \\
\end{split}
\end{equation*}

On the other hand, we have 
\begin{equation*}
\begin{split}
& Y_k r_{\ell} (<\ur{J}, \ur{W}, \ulb{W}, \ulb{I} >) = 
\ep' <\ur{J-\{k\}}, \ur{W^c +\{k + 1\}}, \ulb{W^c +\{k + 1\}}, \ulb{I-\{k+1\}} > \\
& \qquad \qquad - \ep'<\ur{J-\{k\}}, \ur{W^c +\{k+1\}}, \ulb{W^c +\{k\}}, 
\ulb{I-\{k\}} >. \\
\end{split}
\end{equation*}

Here $\ep' = (-1)^{\vert \tI \vert + \vert \{ j \in \tJ : j > k+1\} \vert + \vert 
\{ i \in \tI : i > k +1\} \vert }$ and the equality holds.

Case 3. $k \in \tJ$ and $k+1 \in \tW$.

Then we have 
\[
 Y_k (<\ur{J}, \ur{W}, \ulb{W}, \ulb{I} >)
= <\ur{J-\{k\}+\{k+1\}}, \ur{W -\{k+1\} +\{k\}}, \ulb{W - \{k+1\} +\{k\}}, \ulb{I} >.
\]

So we have 
\begin{equation*}
\begin{split}
& r_{\ell} Y_k (<\ur{J}, \ur{W}, \ulb{W}, \ulb{I} >) = 
 (-1)^{\vert \tI \vert} <\ur{J-\{k\}+\{k+1\}}, \ur{W^c}, \ulb{W^c}, \ulb{I} >.
\\
\end{split}
\end{equation*}

On the other hand, since $k, k+1 \notin \tW^c$, we have 
\begin{equation*}
\begin{split}
 Y_k r_{\ell} (<\ur{J}, \ur{W}, \ulb{W}, \ulb{I} >) = 
  (-1)^{\vert \tI \vert} <\ur{J-\{k\}+\{k+1\}}, \ur{W^c }, \ulb{W^c}, \ulb{I} >. \\
\end{split}
\end{equation*}
So we have the equality.

We omit the proof for the other cases since 
the proofs are almost similar.

Let us prove $r_{\ell} Y_n (<\ur{J}, \ur{W}, \ulb{W}, \ulb{I} >)= Y_n r_{\ell}( <\ur{J}, \ur{W}, \ulb{W}, \ulb{I} >)$, where 
$\ad(Y_n) = E_{\overline{n},n-1} -  E_{\overline{n-1},n}$.

Case 1.  $n-1, n \in \tJ$ or $n-1, n \in \tI$.

If  $n-1, n \in \tJ$,  we have 

\noindent
$
 Y_n (<\ur{J}, \ur{W}, \ulb{W}, \ulb{I} >)
= - <\ur{J-\{n\}-\{n-1\}}, \ur{W+\{n-1\}}, \ulb{W+\{n-1\}}, \ulb{I} >
- <\ur{J-\{n-1\}-\{n\}}, \ur{W+\{n\}}, \ulb{W+\{n\}}, \ulb{I} >
$
 and 

\noindent
$
r_{\ell} Y_n (<\ur{J}, \ur{W}, \ulb{W}, \ulb{I} >)=
(-1)^{\vert \tI \vert+1} (<\ur{J-\{n\}-\{n-1\}}, \ur{W^c+\{n\}}, \ulb{W^c+\{n\}}, \ulb{I} >
+ <\ur{J-\{n-1\}-\{n\}}, \ur{W^c+\{n-1\}}, \ulb{W^c+\{n-1\}}, \ulb{I} >).
$ 

The right hand of the above is equal to 
$Y_n r_{\ell}( <\ur{J}, \ur{W}, \ulb{W}, \ulb{I} >).$

So we have the equality.
The proof is almost same for the case $n-1, n \in \tI$.

Case 2. $n-1 \in \tJ$ and $n \in \tI$.

In this case both sides are zero and the equality holds.

Case 3. $n-1 \in \tJ$ and $n \in \tW$. 

$
 Y_n (<\ur{J}, \ur{W}, \ulb{W}, \ulb{I} >)
= - <\ur{J-\{n-1\}}, \ur{W-\{n\}+\{n-1\}}, \ulb{W-\{n\}+\{n-1\}}, \ulb{I+\{n\}} >
$
 and 
$
r_{\ell} Y_n (<\ur{J}, \ur{W}, \ulb{W}, \ulb{I} >)=
(-1)^{\vert \tI \vert} (<\ur{J-\{n-1\}}, \ur{W^c}, \ulb{W^c}, \ulb{I+\{n\}} >.
$ 

The right hand of the above is equal to 
$Y_n r_{\ell}( <\ur{J}, \ur{W}, \ulb{W}, \ulb{I} >).$

Case 4. $n-1 \in \tJ$ and $n \notin \tI, \tW$. 

$
 Y_n (<\ur{J}, \ur{W}, \ulb{W}, \ulb{I} >)
=  <\ur{J-\{n-1\}}, \ur{W}, \ulb{W}, \ulb{I+\{n\}} >
$
 and 
$
r_{\ell} Y_n (<\ur{J}, \ur{W}, \ulb{W}, \ulb{I} >) = 
(-1)^{\vert \tI \vert+1} (<\ur{J-\{n-1\}}, \ur{W^c-\{n\}+\{n-1\}}, \ulb{W^c-\{n\}+\{n-1\}}, \ulb{I+\{n\}} >.
$

Since $n \in \tW^c$, the right hand of the above is equal to 
$Y_n r_{\ell}( <\ur{J}, \ur{W}, \ulb{W}, \ulb{I} >).$

The rest of the proof follows from similar direct verifications and we omit it.
\end{proof}

From the above,  we can define the explicit base elements of the irreducible constituents  $e_n^+$ and $e_n^-$ of $\bigwedge^n V$ under the action of $Spin(2n)$.

The spaces $e_n^+$ and $e_n^-$ are given respectively by the $\pm 1$ eigenspaces of the $Spin(2n)$-equivariant involution $r_n$ on  $\bigwedge^n V$ defined in Lemma \ref{lem rl}.

\begin{defn}\label{def enpm}
For an exterior product $<\ur{J}, \ur{W}, \ulb{W}, \ulb{I} >$ of degree $n$,
we denote the basis elements of  $e_n^{\pm}$ respectively by 
\[
{<\ur{J}, \ur{W}, \ulb{W}, \ulb{I} >}^{(\pm)}
=\dfrac{1}{\sqrt{2}}\left(<\ur{J}, \ur{W}, \ulb{W}, \ulb{I} > \pm (-1)^{\vert \tI 
\vert} <\ur{J}, \ur{W^c}, \ulb{W^c}, \ulb{I} >\right),
\]
 if the right-hand sides are non-zero.
\end{defn}

In the case of $\vert \tI \vert + \vert \tJ \vert =n$, 
we have 
$
{<\ur{J}, \ulb{I} >}^{(+)} = \sqrt{2} <\ur{J}, \ulb{I} >
$,   ${<\ur{J}, \ulb{I} >}^{(-)} =0$ 
if $\vert \tI \vert \equiv 0 \mod{2}$ and 
$
{<\ur{J}, \ulb{I} >}^{(-)} = \sqrt{2} <\ur{J}, \ulb{I} >
$, 
$
{<\ur{J}, \ulb{I} >}^{(+)} = 0
$
 if $\vert \tI \vert \equiv 1 \mod{2}$.

Hence we have 
\[
<\ur{J}, \ur{W}, \ulb{W}, \ulb{I} >
=\dfrac{1}{\sqrt{2}}\left({<\ur{J}, \ur{W}, \ulb{W}, \ulb{I} >}^{(+)} + {<\ur{J}, \ur{W}, \ulb{W}, \ulb{I} >}^{(-)}\right).
\]

We introduce some other duality maps.

For $Spin(2n+1)$, 
since $\Delta$ is self-dual, the isomorphism 
$
\iota_{\Delta}: \Delta \longrightarrow \Delta^*
$
is given by $\iota_{\Delta}([\ur{I}]) = (-1)^{(n+1) \vert I \vert - \sum \tI}
[\ur{I^c}]^*$.  Here $\tI^{\mathtt{c}}$ is the complement of $\tI$ in 
$\{1,2,\ldots,n\}$ i.e., $\tI \sqcup \tI^{\mathtt{c}} = \{1,2,\ldots,n\}$ and $\sum \tI = \sum_{ i \in \tI} i$ and $\vert \tI \vert$ denotes the number of the elements in $\tI$.

Then we can check easily 
$
 \iota_{\Delta} \dphi_{+} = \dphi_{+}^* \iota_{\Delta}.
$

Also for $Spin(2n)$,  
if $n$ is even positive integer, 
$\Delta^{\pm}$ is self-dual and the isomorphism 
$
\iota_{\Delta}: \Delta^{\pm} \longrightarrow (\Delta^{\pm})^*
$
is given by $\iota_{\Delta}([\ur{I}]) = (-1)^{\sum \tI}[\ur{I^c}]^*$.

If $n$ is odd positive integer, then
$(\Delta^{\pm})^* \cong \Delta^{\mp}$ and the isomorphism 
$
\iota_{\Delta}: \Delta^{\pm} \longrightarrow (\Delta^{\mp})^*
$
is given by $\iota_{\Delta}([\ur{I}]) = \pm (-1)^{\sum \tI} [\ur{I^c}]^*$.

Then we can check easily 
$
 \iota_{\Delta} \dphi = \dphi^* \iota_{\Delta}.
$

Also we define an involutive linear automorphism $\iota'$ on $\Delta \bigotimes \Delta^*$ as follows.

\begin{defn}\label{def iota'del}
For any disjoint subsets $\ur{I}$, $\ur{J}$ in $[1,n]$, $\iota'$ is defined by 
$$\iota'( [\ur{I}] \otimes [\ur{J}]^*) =[\ur{J}] \otimes [\ur{I}]^*.$$
\end{defn}

Then  for $\frak{so}(N)$, it is easily checked that $ \iota' X_k = -Y_k \iota'$, $ \iota' Y_k = -X_k \iota'$ and $ \iota' h_k = -h_k \iota'$ on $\Delta \bigotimes \Delta^*$.

If $N =2n$, $\iota'$ is the linear map satisfying 
$\iota'( \Delta^{\ep_1} \bigotimes (\Delta^{\ep_2})^*) = \Delta^{\ep_2} \bigotimes (\Delta^{\ep_1})^*$, where  $\ep_1, \ep_2 \in \{ \pm 1 \}$.

\section{The equivariant embeddings of $\bigwedge^k V$  to $\Delta \bigotimes 
\Delta^*$}\label{sect phi}

In this section we give the key theorems, which bridge the gap between the spin representations and the tensor calculus.

\begin{thm}\label{thm phiodd}
For any $k$ \,$( 1\leqq k \leqq 2n+1)$, there exists an $\frak{so}(2n+1)$-equivariant embedding 
$\phi_k$ from $\bigwedge^k V$ 
 to $\Delta \bigotimes \Delta^*$
 defined as follows. 

\begin{equation}
 \phi_k(<\ur{J}, \ur{W}, \ulb{W}, \ulb{I} >)
 = \sum_{\substack{ [1,n] - \mathtt{J} - \mathtt{I} \supseteq \mathtt{K}}}
\dfrac{(-1)^{\vert \mathtt{W} - \mathtt{W} \cap \mathtt{K} \vert}}{2^{(
n - \vert \tJ \vert - \vert \tI \vert)/2}} [ \ur{I}, \ur{K}] \otimes [ \ur{J}, 
\ur{K}]^*
\tag{\ref{thm phiodd}.1}\label{eqn phiodd}
 \end{equation}

\begin{equation}
 \phi_k(<\ur{J}, \ur{W}, \ur{0'}, \ulb{W}, \ulb{I} >)
 = \sum_{\substack{ [1,n] - \mathtt{J} - \mathtt{I} \supseteq \mathtt{K}}}
\dfrac{(-1)^{\vert \mathtt{K} - \mathtt{K} \cap \mathtt{W} \vert}}{2^{(
n - \vert \tJ \vert - \vert \tI \vert)/2}} [ \ur{I}, \ur{K}] \otimes [ \ur{J}, 
\ur{K}]^*
\tag{\ref{thm phiodd}.2}\label{eqn phiodd0}
 \end{equation}

 Moreover the above $\phi_k$ is an isometric embedding 
 with respect to the invariant hermitian inner products.

\end{thm}

We give proof of the above theorem in the last of this section.

\begin{rem}\label{rem phiodd}
Here  $\Delta$ stands for $\Delta_{+}$  or $\Delta_{-}$ and in the tensor parts corresponding to $\Delta_{-}$, as for the bracket notation, 
we obey the convention of Remark \ref{rem dphi}.
\end{rem}

If we compare the actions of the center $z$ (see Remark \ref{rem dphi}.)  on both sides, 
as $Pin(2n+1)$ modules,  we have 
\begin{equation}
\Delta_{\pm} \textstyle{\bigotimes} \Delta_{\pm} \cong 
\textstyle{\bigoplus\limits_{i=0}^{n}} (\textstyle{\bigwedge\limits^{2i} V}) 
\end{equation}
and 
\begin{equation}
\Delta_{\pm} \textstyle{\bigotimes} \Delta_{\mp} \cong 
\textstyle{\bigoplus\limits_{i=0}^{n}} (\textstyle{\bigwedge\limits^{2i+1} V})
\end{equation}
and 
the above $\phi_k$'s give $Pin(2n+1)$-equivariant embeddings.

So for the group $Pin(2n+1)$, we have the following theorem.

\begin{thm}\label{thm phipinodd}
For any  $k$ \,$( 1\leqq k \leqq 2n+1)$, there exists a $Pin(2n+1)$-equivariant isometric embedding 
$\phi_k$ from $\bigwedge^k V$ 
 to $\Delta_{\ep_1} \bigotimes \Delta_{\ep_2}^*$
 defined by \eqref{eqn phiodd} and \eqref{eqn phiodd0} in Theorem \ref{thm 
phiodd}. 
Here $\ep_1, \ep_2 \in \{ \pm 1\}$ and $(-1)^k = \ep_1 \ep_2$. 
\end{thm}

We move to the case $N=2n$.

For  $\ep_1, \ep_2 \in \{ \pm 1\}$, by $\phi_{\ell}^{\ep_1 \ep_2}$, we denote an $\frak{so}(2n)$-equivariant embedding from $\bigwedge^{\ell} V$ to $\Delta^{\ep_1} \bigotimes (\Delta^{\ep_2})^*$.  For example, $\phi_{\ell}^{\ep_1 \ep_2}= \phi_{\ell}^{+ -}$ if $\ep_1 = 1$ and $\ep_2 = -1$ and so on.

\begin{thm}\label{thm phieven}
If $n$ is even,  we have $(\Delta^+)^* \cong \Delta^+$ and $(\Delta^-)^* \cong 
\Delta^-$ and only the embeddings $\phi_{\ell}^{\ep_1 \ep_2}$ for any $\ell$  $(1 \leqq \ell \leqq 2n)$ satisfying 
the condition $\ep_1 = (-1)^{n-\ell} \ep_2$ occur.
 
If $n$ is odd,  we have $(\Delta^+)^* \cong \Delta^-$ and $(\Delta^-)^* \cong 
\Delta^+$ and only the embeddings $\phi_{\ell}^{\ep_1 \ep_2}$ for any $\ell$  $(1 \leqq \ell \leqq 2n)$ satisfying 
the condition $\ep_1 = -(-1)^{n-\ell} \ep_2$ occur.

For $\ell \ne n$, $\phi_{\ell}^{\ep_1 \ep_2}$ is given by 
 \begin{equation}
 \phi_{\ell}^{\ep_1 \ep_2}(<\ur{J}, \ur{W}, \ulb{W}, \ulb{I} >)
 = \sum_{\substack{ [1,n] - \mathtt{J} - \mathtt{I} \supseteq \mathtt{K} \\ 
 (-1)^{i + k} = \ep_1  \\ (-1)^{j + k} = \ep_2}}
\dfrac{(-1)^{\vert \mathtt{W} - \mathtt{W} \cap \mathtt{K} \vert}}{2^{(
n - \vert \tJ \vert - \vert \tI \vert-1)/2}} [ \ur{I}, \ur{K}] \otimes [ \ur{J}, \ur{K}]^*.
\tag{\ref{thm phieven}.1}\label{eqn phieven}
 \end{equation}

Here by small letters, we denote the number of the elements in the set indexed by their capital letters.

For $e_n^+$ and $e_n^-$, $\phi_{n}^{\ep_1 \ep_2}$ is given by 
\begin{equation}
\phi_{n}^{\ep_1 \ep_2}({<\ur{J}, \ur{W}, \ulb{W}, \ulb{I} >}^{(\ep_1)}) = 
\sum_{\substack{ [1,n] - \mathtt{J} - \mathtt{I} \supseteq \mathtt{K} \\  (-1)^{i + k} = \ep_1  \\ (-1)^{j + k} = \ep_2}}
\dfrac{(-1)^{\vert \mathtt{W} - \mathtt{W} \cap \mathtt{K} \vert}}{2^{(
n - \vert \tJ \vert - \vert \tI \vert-1)/2}} [ \ur{I}, \ur{K}] \otimes [ \ur{J}, \ur{K}]^* .
\tag{\ref{thm phieven}.2}\label{eqn phievenn}
\end{equation}

 Moreover each of the above $\phi_{\ell}^{\ep_1 \ep_2}$'s is an isometric embedding  with respect to the invariant hermitian inner products.
\end{thm}

We give proofs of the above theorems for $Spin(2n+1)$ and $Spin(2n)$ in the last of this section simultaneously.

\begin{rem}
For $e_n^+$ and $e_n^-$, we note that only the above case occurs from \eqref{eqn evendeldel} and \eqref{eqn evendel+del-}.
\end{rem}

The last statements in Theorem \ref{thm phiodd} and Theorem \ref{thm phieven} 
follow easily from the following fact. 

Since $\phi_k$ and $\phi_k^{\ep_1 \ep_2}$ become the $\frak{so}(N,\C)$-equivariant map, then of course it becomes $\frak{so}(N,\C)_{cpt}$ equivariant map.

$\bigwedge^k V$ and $e_n^+$ and $e_n^-$ are irreducible under the action of 
$\frak{so}(N,\C)_{cpt}$. Hence the invariant hermitian products only differ in a scalar multiple. If we take the highest weight vector of $\bigwedge^k V$ and compare its value of the hermitian inner product with that of its image, the last statements follow easily.

Before we prove the theorem for $Pin(2n)$, we state several properties, which follow from the theorems easily.

For $Spin(2n+1)$, 
from Theorem \ref{thm phiodd}, we have 
the isomorphism $\bigoplus_{i=0}^{n} \bigwedge^{2i} V \cong \Delta \bigotimes 
\Delta^*$.

So we give the inverse of the above isomorphism.

Let us compare the same weight spaces in both sides.
For simplicity, we omit $\phi_k$ since the degrees of the exterior products 
tell us which $\phi_k$ we are considering in the argument.
Let $\tI=\{i_1, i_2, \ldots,i_r\}$ and $\tJ=\{j_1, j_2, \ldots,j_s\}$
 be disjoint subsets of $[1,n]$.

Let $\be_{\tJ}-\be_{\tI}$ denote the weight 
$\be_{\tJ}-\be_{\tI} = \be_{j_1} + \be_{j_2} + \ldots + \be_{j_s} - \be_{i_1} - \be_{i_2} - \ldots - \be_{i_r}$.
Then a base of the weight space of weight $\be_{\tJ}-\be_{\tI}$ in 
$\bigoplus_{i=0}^{n} \bigwedge^{2i} V$
is given by $\{<\ur{J}, \ur{W}, \ulb{W}, \ulb{I} >\}$ if $\vert \tJ \vert + \vert \tI \vert \equiv 0 \pmod{2}$
and $\{<\ur{J}, \ur{W},\ur{0'}, \ulb{W}, \ulb{I} >\}$ if $\vert \tJ \vert + \vert \tI \vert \equiv 1 \pmod{2}$.
Here $\tW$ runs over all the subsets in $[1,n] - \tJ -\tI$.

Also a base of the weight space of weight  $\be_{\tJ}-\be_{\tI}$ in 
$\Delta \bigotimes \Delta^*$ is given by $\{[\ur{I},\ur{K}] \otimes 
[\ur{J},\ur{K}]^*\}$, where $\tK$ runs over all the subsets in $[1,n] - \tJ -\tI$.

Since both bases are hermitian orthonormal bases of $\bigoplus_{i=0}^{n} 
\bigwedge^{2i} V$ and $\Delta \bigotimes \Delta^*$ respectively, the transformation matrices ${\ds \frac{1}{2^{(n - \vert \tJ \vert - \vert \tI \vert)/2}} 
((-1)^{\vert \mathtt{W} - \mathtt{W} \cap \mathtt{K} \vert})_{\tW, \tK}}$ are 
unitary matrices.  Moreover all the components in the above are real, so 
the transformation matrices are orthogonal matrices over $\R$, hence 
$H_{\tJ, \tI, n}= ((-1)^{\vert \mathtt{W} - \mathtt{W} \cap \mathtt{K} \vert})_{\tW, \tK}$ ($\tW$ and $\tK$ run over all the subsets in $[1,n] - \tJ -\tI$.) are Hadamard matrices of size $2^{\vert [1,n] - \tJ -\tI \vert}$.
These Hadamard matrices are equivalent to the tensor products of the standard Hadamard matrix of size $2$.

So the inverse matrices of these matrices are given by transposing the original matrices.

Therefore if $\vert \tJ \vert + \vert \tI \vert \equiv 0 \pmod{2}$, we have 
\[
[ \ur{I}, \ur{K}] \otimes [ \ur{J}, \ur{K}]^*
=\sum_{\substack{ [1,n] - \mathtt{J} - \mathtt{I} \supseteq \mathtt{W}}}
\dfrac{(-1)^{\vert \mathtt{W} - \mathtt{W} \cap \mathtt{K} \vert}}{2^{(
n - \vert \tJ \vert - \vert \tI \vert)/2}} 
<\ur{J}, \ur{W}, \ulb{W}, \ulb{I} >
\]
and if $\vert \tJ \vert + \vert \tI \vert \equiv 1 \pmod{2}$, we have 
\[
[ \ur{I}, \ur{K}] \otimes [ \ur{J}, \ur{K}]^*
=\sum_{\substack{ [1,n] - \mathtt{J} - \mathtt{I} \supseteq \mathtt{W}}}
\dfrac{(-1)^{\vert \mathtt{K} - \mathtt{K} \cap \mathtt{W} \vert}}{2^{(
n - \vert \tJ \vert - \vert \tI \vert)/2}} 
<\ur{J}, \ur{W}, \ur{0'}, \ulb{W}, \ulb{I} >.
\]

Similarly for $Spin(2n)$, 
from Theorem \ref{thm phieven}, 
 we have the isomorphism $e_n^{\ep_1} \bigoplus e_{n-2} \bigoplus \cdots \bigoplus e_{0}  \cong \Delta^{\ep_1} \bigotimes (\Delta^{\ep_2})^*$, here $\ep_1, \ep_2 \in \{\pm 1\}$ and $\ep_2 = (-1)^n \ep_1$.

We give the inverse of the above isomorphism.

As before, if we compare the same weight spaces in both sides, 
we obtain the inverse as follows.
We also omit $\phi_k^{\ep_1 \ep_2}$ since the degrees of the exterior products and the parities of the indices in the brackets tell us which $\phi_k^{\ep_1 \ep_2}$ we are considering in the argument.

For the irreducible decomposition $\Delta^{\ep_1} \bigotimes (\Delta^{\ep_2})^*$ with $\ep_2 = (-1)^n \ep_1$, in the case of $(-1)^{k + i} =\ep_1$ and $(-1)^{k + j} = \ep_2$, we have 
\begin{equation*}
 \begin{split}
 & [ \ur{I}, \ur{K}] \otimes [ \ur{J}, \ur{K}]^* = \sum_{\substack{ [1,n] - \mathtt{J} 
- \mathtt{I} \supseteq \mathtt{W} \\ 
 i + j + 2w \leqq n-2 }}
\dfrac{(-1)^{\vert \mathtt{W} - \mathtt{W} \cap \mathtt{K} \vert}}{2^{(
n - \vert \tJ \vert - \vert \tI \vert-1)/2}} <\ur{J}, \ur{W}, \ulb{W}, \ulb{I} > \\
& \qquad \qquad +
\sum_{\substack{ [1,n] - \mathtt{J} - \mathtt{I} \supseteq \mathtt{W} \\ 
 i + j + 2w = n }}
\dfrac{(-1)^{\vert \mathtt{W} - \mathtt{W} \cap \mathtt{K} \vert}}{2^{(
n - \vert \tJ \vert - \vert \tI \vert + 1)/2}} {<\ur{J}, \ur{W}, \ulb{W}, \ulb{I} >}^{(\ep_1)}. \\
\end{split}
\end{equation*}

We note that in the second sum, if $n-i-j >0$, $(-1)^{\vert \mathtt{W} - \mathtt{W} \cap \mathtt{K} \vert} {<\ur{J}, \ur{W}, \ulb{W}, \ulb{I} >}^{(\ep_1)}$
$=$ $(-1)^{\vert \mathtt{W^c} - \mathtt{W^c} \cap \mathtt{K} \vert} {<\ur{J}, \ur{W^c}, \ulb{W^c}, \ulb{I} >}^{(\ep_1)}$ and we count it twice.

For the irreducible decomposition $\Delta^{\ep_1} \bigotimes (\Delta^{\ep_2})^*$ with $\ep_2 = -(-1)^n \ep_1$, 
in the case of $(-1)^{k + i} = \ep_1$ and $(-1)^{k + j} =\ep_2$, we have 

\begin{equation*}
[ \ur{I}, \ur{K}] \otimes [ \ur{J}, \ur{K}]^* = \sum_{\substack{ [1,n] - \mathtt{J} - \mathtt{I} \supseteq \mathtt{W} \\ 
 i + j + 2w \leqq n-1}}
\dfrac{(-1)^{\vert \mathtt{W} - \mathtt{W} \cap \mathtt{K} \vert}}{2^{(
n - \vert \tJ \vert - \vert \tI \vert-1)/2}} <\ur{J}, \ur{W}, \ulb{W}, \ulb{I} 
>.\end{equation*}

The matrix $\left( (-1)^{\vert \mathtt{W} - \mathtt{W} \cap \mathtt{K} \vert} 
\right)_{\tK, \tW}$ is an Hadamard matrix of size $2^{n - \vert \tJ \vert - \vert \tI \vert-1}$, where $\tW$ and $\tK$ run over the subsets in $[1,n] - \mathtt{J} - \mathtt{I}$ satisfying $\vert \tW \vert \leqq (n-1-i-j)/2$ and  
$(-1)^{\vert \tK \vert + i} = \ep_1$.
These Hadamard matrices are also equivalent to the former ones.
\footnote{We thank Prof. E. Bannai for this comment.}

Using the above theorem, we obtain the following $Pin(2n)$ ($O(2n)$) equivariant  isometric embeddings.

\begin{thm}\label{thm phipin}
As $Pin(2n)$ $(O(2n))$ modules, for any even $\ell$ $(1 \leqq \ell \leqq 2n)$, there exists an equivariant isometric embedding 
 $\phi_{\ell} : \bigwedge^{\ell} V \longrightarrow \Delta \bigotimes \Delta^*$ given by 
\begin{equation*}
 \phi_{\ell}(<\ur{J}, \ur{W}, \ulb{W}, \ulb{I} >)
 = \sum_{\substack{ [1,n] - \mathtt{J} - \mathtt{I} \supseteq \mathtt{K}}}
\dfrac{(-1)^{\vert \mathtt{W} - \mathtt{W} \cap \mathtt{K} \vert}}{2^{(
n - \vert \tJ \vert - \vert \tI \vert)/2}} [ \ur{I}, \ur{K}] \otimes [ \ur{J}, 
\ur{K}]^*
\end{equation*}
and 
for any odd $\ell$ ($ 1 \leqq \ell \leqq 2n$), there exists an equivariant isometric embedding  $\phi_{\ell} : \bigwedge^{\ell} V \longrightarrow \Delta \bigotimes \Delta^*$ given by 
\begin{equation*}
 \phi_{\ell}(<\ur{J}, \ur{W}, \ulb{W}, \ulb{I} >)
 = \sum_{\substack{ [1,n] - \mathtt{J} - \mathtt{I} \supseteq \mathtt{K}}}
\dfrac{(-1)^{\vert \mathtt{W} - \mathtt{W} \cap \mathtt{K} \vert +\vert \tI \vert +\vert \tK \vert}}{2^{(n - \vert \tJ \vert - \vert \tI \vert)/2}} [ \ur{I}, \ur{K}] \otimes [ \ur{J}, \ur{K}]^*.
\end{equation*}

\end{thm}

\begin{rem}\label{rem phipin}
This theorem is crucial for us to consider the invariant theory for the group 
$Pin(2n)$ $(O(2n))$.
The above $\phi_{\ell}$ is defined as the sum $\dfrac{1}{\sqrt{2}} 
\left(\phi_{\ell}^{+ +} + \phi_{\ell}^{- -}\right)$ for even $\ell ( \ne n )$  and $\dfrac{1}{\sqrt{2}} \left(\phi_{\ell}^{+ -} - \phi_{\ell}^{- +}\right)$ for odd $\ell ( \ne n )$ and $\phi_{n} = \phi_{n}^{+ +} \bigoplus \phi_{n}^{- -}$ for even $n$ and $\phi_{n} =  \phi_{n}^{+ -} \bigoplus -\phi_{n}^{- +}$ for odd $n$.
\end{rem}

\begin{proof}

From the above remark, $\phi_{\ell}$ is $SO(2n)$-equivariant. So in order to prove the theorem,  it is enough to check that  $\phi_{\ell}$ is $(\ua{n}-\ub{n})$-equivariant.

We note that the action of $(\ua{n}-\ub{n})$ on $V$ is given by the transposition of $\ua{n}$ and $\ub{n}$.

We first prove the case for even $\ell ( \ne n )$.

Case 1.  \quad $n \notin \tI, \tJ, \tW$.

Then we have $(\ua{n}-\ub{n})<\ur{J}, \ur{W}, \ulb{W}, \ulb{I} >=<\ur{J}, \ur{W}, \ulb{W}, \ulb{I} >$ and calculate  $(\ua{n}-\ub{n}) \phi_{\ell}(<\ur{J}, \ur{W}, \ulb{W}, \ulb{I} >)$. 
We take cases for $n \in \tK$ and 
 $n \notin \tK$ in the the expansion of $\phi_{\ell}(<\ur{J}, \ur{W}, \ulb{W}, \ulb{I} >)$ in Theorem \ref{thm phipin}.

Case 1-1.  $n \in \tK$.

We put $\tK_1 = \tK - \{n\}$.  From the Lemma \ref{lem clactdelbase}, the action of $(\ua{n}-\ub{n})$ is given as follows.
(We note that in that lemma, the action is described only for the base elements with the increasing indices in the brackets.)
 
\begin{equation*}
\begin{split}
 &(\ua{n}-\ub{n}) [ \ur{I}, \ur{K}] \otimes [ \ur{J}, \ur{K}]^*  \\
&=\sign \begin{pmatrix} \ur{I}, & \ur{K} \\ \ur{I + K_1}, & \{\ur{n}\} \\ \end{pmatrix} 
\sign \begin{pmatrix} \ur{J}, & \ur{K} \\ \ur{J + K_1}, & \{\ur{n}\} \\ \end{pmatrix}
(-1)^{\vert \tI \vert + \vert \tK_1 \vert + \vert \tJ \vert + \vert \tK_1 \vert} [\ur{I + K_1}] \otimes [\ur{J + K_1}]^* \\
& =(-1)^{\vert \tI \vert  + \vert \tJ \vert } 
[\ur{I}, \ur{K_1}] \otimes [\ur{J},  \ur{K_1}]^*.
\end{split}
\end{equation*}

Here $\ur{I + K_1}$ denotes the union of the set $\ur{I}$ and $\ur{K_1}$ and we have the convention that the indices in the set $\ur{I + K_1}$ are arranged in the increasing order. 
$\sign \begin{pmatrix} \ur{I}, & \ur{K} \\ \ur{I + K_1}, & \{\ur{n}\} 
\\ \end{pmatrix}$ denotes the sign of the permutation obtained by arranging 
$\ur{I}$, $\ur{K}$ into $\ur{I + K_1}$, $\{\ur{n}\}$ in this order.

Case 1-2. $n \notin \tK$.

We put $\tK_1 = \tK$. The action of $(\ua{n}-\ub{n})$ is given by 
\begin{equation*}
\begin{split}
& (\ua{n}-\ub{n}) [ \ur{I}, \ur{K}] \otimes [ \ur{J}, \ur{K}]^*  \\
& =\sign \begin{pmatrix} \ur{I}, & \ur{K_1} \\ \ur{I + K_1}  \\ \end{pmatrix}
\sign \begin{pmatrix} \ur{J}, & \ur{K_1} \\ \ur{J + K_1}  \\ \end{pmatrix}
(-1)^{\vert \tI \vert + \vert \tK_1 \vert - 1 + \vert \tJ \vert + \vert \tK_1 \vert - 1} [\ur{I + K_1}, \ur{n}] \otimes [\ur{J + K_1}, \ur{n}]^* \\
& =(-1)^{\vert \tI \vert  + \vert \tJ \vert } 
[\ur{I}, \ur{K_1+\{n\}}] \otimes [\ur{J},  \ur{K_1+\{n\}}]^*.
\end{split}
\end{equation*}

Since $\vert \tI \vert + \vert \tJ \vert + 2 \vert \tW \vert = \ell$, we have 
$(-1)^{\vert \tI \vert  + \vert \tJ \vert }=1$.

So in this case, the equality holds.

Case 2.  \quad $n \in \tI$.

Then $n \notin \tW, \tJ, \tK$ and $(\ua{n}-\ub{n})<\ur{J}, \ur{W}, \ulb{W}, \ulb{I} 
>=<\ur{J+\{n\}}, \ur{W}, \ulb{W}, \ulb{I-\{n\}} >$.

The image of $\phi_{\ell}$ is 
\begin{equation*}
 \phi_{\ell}(<\ur{J+\{n\}}, \ur{W}, \ulb{W}, \ulb{I-\{n\}} >)
 = \sum_{\substack{ [1,n] - \mathtt{J} - \mathtt{I} \supseteq \mathtt{K}}}
\dfrac{(-1)^{\vert \mathtt{W} - \mathtt{W} \cap \mathtt{K} \vert}}{2^{(
n - \vert \tJ \vert - \vert \tI \vert)/2}} [ \ur{I-\{n\}}, \ur{K}] \otimes 
[ \ur{J+\{n\}}, \ur{K}]^*.
\end{equation*}

We calculate $(\ua{n}-\ub{n}) \phi_{\ell}(<\ur{J}, \ur{W}, \ulb{W}, \ulb{I} >)$. 
Then we have 
\begin{equation*}
\begin{split}
& (\ua{n}-\ub{n}) [ \ur{I}, \ur{K}] \otimes [ \ur{J}, \ur{K}]^*  \\
& =\sign \begin{pmatrix} \ur{I}, & \ur{K} \\ \ur{I-\{n\} + K}, & \{n\} 
\\ \end{pmatrix}
\sign \begin{pmatrix} \ur{J}, & \ur{K} \\ \ur{J + K}  \\ \end{pmatrix}
(-1)^{\vert \tI \vert -1 + \vert \tK \vert + \vert \tJ \vert + \vert \tK \vert+1} 
[\ur{I-\{n\} + K}] \otimes [\ur{J + K}, \ur{n}]^* \\
& =(-1)^{\vert \tI \vert  + \vert \tJ \vert} 
[\ur{I-\{n\}}, \ur{K}] \otimes [\ur{J +\{n\}},  \ur{K}]^*.
\end{split}
\end{equation*}

Here $\ur{I-\{n\} + K}$ denotes the set obtained by the set theoretical difference and sum and we consider the indices in it are arranged in the increasing order.

So we have the equality.

The proof for the case $n \in \tJ$ is almost the same and we omit it.

Case 3. \quad $n \in \tW$.
 
Then $n \notin \tI, \tJ, \tK$ and $(\ua{n}-\ub{n})<\ur{J}, \ur{W}, \ulb{W}, \ulb{I} >=-<\ur{J}, \ur{W}, \ulb{W}, \ulb{I} >$.

The image of $\phi_{\ell}$ of this element is 
$-\sum_{\substack{ [n] - \mathtt{J} - \mathtt{I} \supseteq \mathtt{K}}}
({(-1)^{\vert \mathtt{W} - \mathtt{W} \cap \mathtt{K} \vert}}/{2^{(
n - \vert \tJ \vert - \vert \tI \vert)/2}}) [ \ur{I}, \ur{K}] \otimes [ \ur{J}, \ur{K}]^*$.

We calculate  $(\ua{n}-\ub{n}) \phi_{\ell}(<\ur{J}, \ur{W}, \ulb{W}, \ulb{I} >)$. We take cases for $n \in \tK$ and 
 $n \notin \tK$.

Case 3.1  $n \in \tK$.

We put $\tK_1 = \tK -\{n\}$ and the image of the sum over such $\tK$'s is 
given by 
$\sum_{\substack{ [n] - \mathtt{J} - \mathtt{I} - \{n\} \supseteq \mathtt{K_1}}}({(-1)^{\vert \mathtt{W} - \mathtt{W} \cap \mathtt{K_1} \vert-1}}/{2^{(
n - \vert \tJ \vert - \vert \tI \vert)/2}}) [ \ur{I}, \ur{K_1}] \otimes [ \ur{J}, \ur{K_1}]^*$.

Case 3.2 $n \notin \tK$.

We put $\tK_1 = \tK$ and the image of the sum over such $\tK$'s is 
given by 
$
\sum_{\substack{ [n] - \mathtt{J} - \mathtt{I} - \{n\} \supseteq \mathtt{K_1}}}
({(-1)^{\vert \mathtt{W} - \mathtt{W} \cap \mathtt{K_1} \vert}}/{2^{(
n - \vert \tJ \vert - \vert \tI \vert)/2}}) [ \ur{I}, \ur{K_1+\{n\}}] \otimes [ \ur{J}, \ur{K_1+\{n\}}]^*$.

If we note $\vert \mathtt{W} - \mathtt{W} \cap \mathtt{K_1} \vert \equiv \vert \mathtt{W} - \mathtt{W} \cap \mathtt{(K_1+\{n\})} \vert +1 \mod{2}$, 
we have the equality.

Next we prove the case for odd $\ell ( \ne n )$.

Case 1. \quad $n \notin \tI, \tJ, \tW$.

Then $(\ua{n}-\ub{n})<\ur{J}, \ur{W}, \ulb{W}, \ulb{I} >=<\ur{J}, \ur{W}, \ulb{W}, \ulb{I} >$.

We calculate  $(\ua{n}-\ub{n}) \phi_{\ell}(<\ur{J}, \ur{W}, \ulb{W}, \ulb{I} >)$. 

Case 1.1  $n \in \tK$.

We put $\tK_1 = \tK - \{n\}$, the action of $(\ua{n}-\ub{n})$ is given by 
\begin{equation*}
 (\ua{n}-\ub{n}) [ \ur{I}, \ur{K}] \otimes [ \ur{J}, \ur{K}]^*  
 =(-1)^{\vert \tI \vert  + \vert \tJ \vert } 
[\ur{I}, \ur{K_1}] \otimes [\ur{J},  \ur{K_1}]^* 
 =-[\ur{I}, \ur{K_1}] \otimes [\ur{J},  \ur{K_1}]^*,
\end{equation*}

Since $\vert \tI \vert + \vert \tJ \vert + 2 \vert \tW \vert = \ell$, we have 
$(-1)^{\vert \tI \vert  + \vert \tJ \vert }=-1$.

Case 1.2 $n \notin \tK$.

We put $\tK_1 = \tK$, the action of $(\ua{n}-\ub{n})$ is given by 
\begin{equation*}
 (\ua{n}-\ub{n}) [ \ur{I}, \ur{K}] \otimes [ \ur{J}, \ur{K}]^*  
 =-[\ur{I}, \ur{K_1+\{n\}}] \otimes [\ur{J},  \ur{K_1}+\{n\}]^*.
\end{equation*}

In both cases the signatures are $-1$ and the parities of the numbers of the elements in $\tK$ change.  Since the signatures in the definition of $\phi_{\ell}$ are given by $(-1)^{\vert \mathtt{W} - \mathtt{W} \cap \mathtt{K} \vert +\vert \tI \vert +\vert \tK \vert}$ in this case, we have the equality.

Case 2. \quad $n \in \tI$.

Then $n \notin \tW, \tJ, \tK$ and $(\ua{n}-\ub{n})<\ur{J}, \ur{W}, \ulb{W}, \ulb{I} >=<\ur{J+\{n\}}, \ur{W}, \ulb{W}, \ulb{I-\{n\}} >$.

The image of $\phi_{\ell}$ is 
\begin{equation*}
 \phi_{\ell}(<\ur{J+\{n\}}, \ur{W}, \ulb{W}, \ulb{I-\{n\}} >)
 = \sum_{\substack{ [n] - \mathtt{J} - \mathtt{I} \supseteq \mathtt{K}}}
\dfrac{(-1)^{\vert \mathtt{W} - \mathtt{W} \cap \mathtt{K} \vert +\vert \tI \vert +\vert \tK \vert}}{2^{(
n - \vert \tJ \vert - \vert \tI \vert)/2}} [ \ur{I-\{n\}}, \ur{K}] \otimes [ \ur{J+\{n\}}, \ur{K}]^*.
\end{equation*}

We calculate $(\ua{n}-\ub{n}) \phi_{\ell}(<\ur{J}, \ur{W}, \ulb{W}, \ulb{I} >)$.We have 
\begin{equation*}
 (\ua{n}-\ub{n}) [ \ur{I}, \ur{K}] \otimes [ \ur{J}, \ur{K}]^*  
 =(-1)^{\vert \tI \vert  + \vert \tJ \vert} 
[\ur{I-\{n\}}, \ur{K}] \otimes [\ur{J +\{n\}},  \ur{K}]^* 
 =-[\ur{I-\{n\}}, \ur{K}] \otimes [\ur{J +\{n\}},  \ur{K}]^*.
\end{equation*}

In this case the parity of $\tI$ changes and 
we have the equality.

The proof for the case $n \in \tJ$ is almost the same and we omit it.

Case 3. \quad $n \in \tW$.

Then $n \notin \tI, \tJ, \tK$ and $(\ua{n}-\ub{n})<\ur{J}, \ur{W}, \ulb{W}, \ulb{I} >=-<\ur{J}, \ur{W}, \ulb{W}, \ulb{I} >$.

The image of $\phi_{\ell}$ of this element is 
$
-\sum_{\substack{ [n] - \mathtt{J} - \mathtt{I} \supseteq \mathtt{K}}}
({(-1)^{\vert \mathtt{W} - \mathtt{W} \cap \mathtt{K} \vert +\vert \tI \vert +\vert \tK \vert}}/{2^{(n - \vert \tJ \vert - \vert \tI \vert)/2}}) [ \ur{I}, \ur{K}] \otimes [ \ur{J}, \ur{K}]^*.
$

We calculate  $(\ua{n}-\ub{n}) \phi_{\ell}(<\ur{J}, \ur{W}, \ulb{W}, \ulb{I} >)$. We take cases for $n \in \tK$ and 
 $n \notin \tK$.

Case 3.1 $n \in \tK$.

We put $\tK_1 = \tK -\{n\}$ and the image of the sum over such $\tK$'s is 
given by 
$
\sum_{\substack{ [n] - \mathtt{J} - \mathtt{I} - \{n\} \supseteq \mathtt{K_1}}}
({(-1)^{\vert \mathtt{W} - \mathtt{W} \cap \mathtt{K_1} \vert+\vert \tI \vert +\vert \tK_1 \vert+1}}/{2^{(
n - \vert \tJ \vert - \vert \tI \vert)/2}})  [ \ur{I}, \ur{K_1}] \otimes [ \ur{J}, \ur{K_1}]^*.
$

Case 3.2 $n \notin \tK$.

We put $\tK_1 = \tK$ and the image of the sum over such $\tK$'s is 
given by 
$
\sum_{\substack{ [n] - \mathtt{J} - \mathtt{I} - \{n\} \supseteq \mathtt{K_1}}}
({(-1)^{\vert \mathtt{W} - \mathtt{W} \cap \mathtt{K_1} \vert+\vert \tI \vert +\vert \tK_1 \vert+1}}/{2^{(
n - \vert \tJ \vert - \vert \tI \vert)/2}}) [ \ur{I}, \ur{K_1+\{n\}}] \otimes [ \ur{J}, \ur{K_1+\{n\}}]^*.
$

If we note $\vert \mathtt{W} - \mathtt{W} \cap \mathtt{K_1} \vert+\vert \tI \vert +\vert \tK_1 \vert \equiv \vert \mathtt{W} - \mathtt{W} \cap \mathtt{(K_1+\{n\})} \vert + \vert \tI \vert +\vert (\tK_1+\{n\}) \vert \mod{2}$, 
we have the equality.

In the case of $\ell = n$, only the difference from the other is a scalar multiple $\dfrac{1}{\sqrt{2}}$ and the above proof goes well in this case.

The verifications for the other cases are almost similar and 
 we omit it.

\end{proof}

\begin{rem}

As before, if we compare the same weight spaces in both sides of the above, 
we obtain the inverse as follows.

If $\vert \tJ \vert + \vert \tI \vert \equiv 0 \pmod{2}$, we have 
\[
[ \ur{I}, \ur{K}] \otimes [ \ur{J}, \ur{K}]^*
=\sum_{\substack{ [1,n] - \mathtt{J} - \mathtt{I} \supseteq \mathtt{W}}}
\dfrac{(-1)^{\vert \mathtt{W} - \mathtt{W} \cap \mathtt{K} \vert}}{2^{(
n - \vert \tJ \vert - \vert \tI \vert)/2}} 
<\ur{J}, \ur{W}, \ulb{W}, \ulb{I} >
\]
and if $\vert \tJ \vert + \vert \tI \vert \equiv 1 \pmod{2}$, we have 
\[
[ \ur{I}, \ur{K}] \otimes [ \ur{J}, \ur{K}]^*
=\sum_{\substack{ [1,n] - \mathtt{J} - \mathtt{I} \supseteq \mathtt{W}}}
\dfrac{(-1)^{\vert \mathtt{K} - \mathtt{K} \cap \mathtt{W} \vert +\vert \tI \vert 
+\vert \tK \vert}}{2^{(n - \vert \tJ \vert - \vert \tI \vert)/2}} 
<\ur{J}, \ur{W}, \ulb{W}, \ulb{I} >.
\]
\end{rem}

As for the composition of  $r_{\ell}$ in Lemma \ref{lem rl} and $\phi_{\ell}$, 
$\phi_{\ell}^{\ep_1 \ep_2}$, we have the following.

\begin{lem}\label{lem phirl}
For $\frak{so}(2n+1)$, as a homomorphism from $\bigwedge^{\ell} V$ to $\Delta \bigotimes \Delta^*$, it holds that  $ \phi_{2n+1-\ell} \circ r_{\ell} = \phi_{\ell}$.

For $\frak{so}(2n)$, as a homomorphism from $\bigwedge^{\ell} V$ to $\Delta^{\ep_1} \bigotimes (\Delta^{\ep_2})^*$, 
it holds that 
\[
\phi_{2n-\ell}^{\ep_1 \ep_2} \circ r_{\ell} = \ep_2 (-1)^n \phi_{\ell}^{\ep_1 
\ep_2}.
\]
\end{lem}

\begin{proof}
The verification is straightforward.  
So we omit the proof.

\end{proof}

\begin{rem}

If we admit Theorem \ref{thm phiodd} and Theorem \ref{thm phieven}, from this lemma, we obtain the result in Lemma \ref{lem rl}, which says that  that $r_{\ell}$ is an $\frak{so}(N)$-equivariant isomorphism.
  For,   $\phi_{\ell}$ and $\phi_{N-\ell}$ are equivariant embedding to the same irreducible component in $\Delta \bigotimes (\Delta)^*$ respectively. 
This $r_{\ell}$ comes from the determinant of $SO(N)$.
\end{rem}

Finally we give the proofs of Theorem \ref{thm phiodd} and \ref{thm phieven} simultaneously.

\begin{proof}

Since $\bigwedge^{\ell} V$ ($ \ell \ne n$ for $\frak{so}(2n)$), $e_n^+$ and 
$e_n^-$ are irreducible, in order to prove the theorem,  it is enough to show that $\phi_{\ell}$ and $\phi_{\ell}^{\ep_1 \ep_2}$ are non-zero equivariant homomorphisms.

First we check the images of the highest weight vectors
 under the map $\phi_{\ell}$ and $\phi_{\ell}^{\ep_1 \ep_2}$ are non-zero.

For abbreviation (and to avoid some confusion), we use the following notation 
hereafter.

For $1 \leqq k \leqq r \leqq n$, by $[k,r]$, we denote the set 
$\{k, k+1, \ldots, r\}$ and by $\overline{[r,k]}$, we denote the set 
$\{\overline{r},\overline{r-1},\ldots,\overline{k} \}$.

For $\frak{so}(2n+1)$,
the image of the highest weight for $\ell \leqq n$ is 
$\phi_{\ell}(<\ur{[1,\ell]}>) = \sum_{[\ell +1,n] \supseteq \mathtt{K}}
2^{(n - \ell)/2} [\ur{K}] \otimes [ \ur{J}, \ur{K}]^* \ne 0$ and for $\ell \geqq n$, the image of the highest weight is 
$\phi_{\ell}(<\ur{[1,n]}, \ur{0'}, \ulb{[n, \ell +1]}>) = \sum_{[\ell +1,n] \supseteq \mathtt{K}}
2^{(n - \ell)/2} [\ur{K}] \otimes [ \ur{J}, \ur{K}]^* \ne 0$ , since $\tW=[\ell+1,n]$.

For $\frak{so}(2n)$,
 the image of the highest weight for $\ell < n$ is 
\[
 \phi_{\ell}^{\ep_1 \ep_2}(<\ur{[1,\ell]}>) = \sum_{[\ell+1,n] \supseteq \mathtt{K}}
2^{(n - \ell-1)/2} [\ur{K}] \otimes [ \ur{[1,\ell]}, \ur{K}]^* \ne 0
\]
 and if $\ell > n$, the image of the highest weight is 
\[
 \phi_{\ell}^{\ep_1 \ep_2}(<\ur{[1,n]}, \ulb{[n,\ell+1]}>) = \sum_{[\ell+1,n] 
\supseteq \mathtt{K}}
2^{(n - \ell-1)/2} (-1)^{n-\ell- \vert \tK \vert}[\ur{K}] \otimes [ \ur{[1,\ell]}, \ur{K}]^* \ne 0,
\]
 since $\tW=[\ell+1,n]$.
For $\ell = n$, if we put $\ep =(-1)^n$, we have 
\[
 \phi_{n}^{+ \ep}({<\ur{[1,n]}>}) =  [\ur{\emptyset}] \otimes [ \ur{[1,n]}]^* \ne 0
\]
\[
 \phi_{n}^{- \ep}({<\ur{[1,n-1]},\ulb{n} >}) =  [\ur{n}] \otimes [ \ur{[1,n-1]}]^* \ne 
0.
\]

For $\frak{so}(2n)$ and $\ell = n$, 
if we define the linear homomorphism $\phi_{n}^{even} = \phi_{n}^{+ +} \bigoplus 
\phi_{n}^{- -}$ for even $\ell = n$ and  $\phi_{n}^{even} =  \phi_{n}^{+ -} 
\bigoplus  \phi_{n}^{- +}$ for odd $\ell = n$, then 
$\phi_n^{even}$ is of the form 
\begin{equation*}
\phi_{n}^{even}(<\ur{J}, \ur{W}, \ulb{W}, \ulb{I} >)
 = \sum_{\substack{ [1,n] - \mathtt{J} - \mathtt{I} \supseteq \mathtt{K}}}
\dfrac{(-1)^{\vert \mathtt{W} - \mathtt{W} \cap \mathtt{K} \vert}}{2^{(
n - \vert \tJ \vert - \vert \tI \vert)/2}} [ \ur{I}, \ur{K}] \otimes [ \ur{J}, \ur{K}]^*.
 \end{equation*}

If the above $\phi_n^{even}$ is $\frak{so}(2n)$-equivariant 
isomorphism, then 
$\phi_n^{\ep_1 \ep_2}$ are $\frak{so}(2n)$-equivariant isomorphism.
We note that the above $\phi_n^{even}$ is of the same form to \eqref{eqn phiodd}.

For $\frak{so}(N)$, from the definitions (\ref{def iota'ext}) and (\ref{def 
iota'del}), 
we can check easily  
$\phi_{\ell} \iota' = \iota' \phi_{\ell}$ and 
$\phi_{\ell}^{\ep_1 \ep_2} \iota' = \iota' \phi_{\ell}^{\ep_2 \ep_1}$.

 Therefore to prove the theorem, it is enough to show $X_i \phi_{\ell} = \phi_{\ell} X_i$ and $X_i \phi_{\ell}^{\ep_1 \ep_2} = \phi_{\ell}^{\ep_1 \ep_2} X_i$ for any $i \in \{1,2,\ldots,n\}$, since
 $X_i = - \iota' Y_i \iota'$ in both spaces  $\bigwedge^{\ell} V$ and $\Delta 
\bigotimes \Delta^*$.  
Since $X_i$'s and $Y_i$'s generate the Lie algebra $\frak{so}(N)$, we have the 
theorems.

We note that the definitions of $\phi_{\ell}$, $\phi_{\ell}^{\ep_1 \ep_2}$ and $\phi_n^{even}$ for the elements $<\ur{J}, \ur{W}, \ulb{W}, \ulb{I} >$ are the same up to the scalar and the action of $X_i$ and $Y_i$ for $i=1,2,\ldots, n-1$ are also same for the exterior products of $V$ and $\Delta$ and $\Delta^*$.
 So if we show $X_p \phi_{\ell}(<\ur{J}, \ur{W}, \ulb{W}, \ulb{I} >) = \phi_{\ell} X_p(<\ur{J}, \ur{W}, \ulb{W}, \ulb{I} >)$ $(1 \leqq p \leqq n-1)$ for $\frak{so}(2n+1)$, the proof holds for the rest.

The proof goes by case-by-case verifications.

Proof of $X_p \phi_{\ell}(<\ur{J}, \ur{W}, \ulb{W}, \ulb{I} >) = \phi_{\ell} X_p(<\ur{J}, \ur{W}, \ulb{W}, \ulb{I} >)$ for $1 \leqq p < n$

Case 1. \quad $p+1 \in \tJ$ and $p \in \tW$.

Then $p, p+1 \notin \tI$ and $p \notin \tJ$ and $p+1 \notin \tW$.
We calculate $\phi_{\ell} X_p$.

Since $X_p <\ur{J}, \ur{W}, \ulb{W}, \ulb{I} >$ $=$ $-<\ur{J}, \ur{W}, \ulb{W-\{p\} +\{p+1\}}, \ulb{I} >$ $=$ $<\ur{J-\{p+1\} + \{p\}}, \ur{W-\{p\} +\{p+1\}}, 
\ulb{W-\{p\} +\{p+1\}}, \ulb{I} >$,
 we have 
\begin{align*}
&\phi_{\ell} (X_p <\ur{J}, \ur{W}, \ulb{W}, \ulb{I} > )  \\
&=
 \sum_{\substack{ [1,n] - (\mathtt{J}-\{p+1\} + \{p\}) - \mathtt{I} \supseteq 
\mathtt{K}}}\dfrac{(-1)^{\vert (\mathtt{W}-\{p\} +\{p+1\}) - \mathtt{K} \cap 
(\mathtt{W}-\{p\} +\{p+1\}) \vert}}{2^{(n - \vert \tJ \vert - \vert \tI \vert)/2}} [ \ur{I}, \ur{K}] \otimes [ \ur{J-\{p+1\} + \{p\}}, \ur{K}]^*. \\
\end{align*}

We take two cases for $ p+1 \in \tK$ and $ p+1 \notin \tK$.
(We note that $p \notin \tK$.)

Case 1-1. \quad $ p+1 \in \tK$.

We put $\tK_1 = \tK - \{p+1\}$. Then if we note that the indices in the bracket have alternating property, we have 
$[ \ur{I}, \ur{K}] \otimes [ \ur{J-\{p+1\} + \{p\}}, \ur{K}]^*
= - [ \ur{I}, \ur{K_1+\{p+1\}}] \otimes [ \ur{J}, \ur{K_1+\{p\}}]^*$.
So the coefficient of $[ \ur{I}, \ur{K_1+\{p+1\}}] \otimes [ \ur{J}, 
\ur{K_1+\{p\}}]^*$ is ${(-1)^{\vert \mathtt{W} - \mathtt{K_1} \cap \mathtt{W} \vert}}/{2^{(n - \vert \tJ \vert - \vert \tI \vert)/2}}$.

Case 1-2. \quad $ p+1 \notin \tK$.

We put $\tK_1 = \tK$. 
The coefficient of $[ \ur{I}, \ur{K_1}] \otimes [ \ur{J-\{p+1\} + \{p\}}, 
\ur{K_1}]^*$ is ${(-1)^{\vert \mathtt{W} - \mathtt{K_1} \cap \mathtt{W} \vert}}/{2^{(n - \vert \tJ \vert - \vert \tI \vert)/2}}$.

Next we calculate $X_p \phi_{\ell}$.
\begin{equation*}
X_p \phi_{\ell}(<\ur{J}, \ur{W}, \ulb{W}, \ulb{I} >)
 = \sum_{\substack{ [1,n] - \mathtt{J} - \mathtt{I} \supseteq \mathtt{K}}}
\dfrac{(-1)^{\vert \mathtt{W} - \mathtt{W} \cap \mathtt{K} \vert}}{2^{(
n - \vert \tJ \vert - \vert \tI \vert)/2}} 
X_p ([ \ur{I}, \ur{K}] \otimes [ \ur{J}, \ur{K}]^*).
\end{equation*}

We take two cases for $p \in \tK$ and $p \notin \tK$.
(We note that $p+1 \notin \tK$ since $p+1 \in \tJ$.)

Case 1-$1'$. \quad  $p \in \tK$.

We put $\tK_1 = \tK -\{p\}$.
Then $X_p ([ \ur{I}, \ur{K}] \otimes [ \ur{J}, \ur{K}]^*) =
X_p ([ \ur{I}, \ur{K_1+\{p\}}] \otimes [ \ur{J}, \ur{K_1+\{p\}}]^*)
 = - [ \ur{I}, \ur{K_1+\{p+1\}}] \otimes [ \ur{J}, \ur{K_1+\{p\}}]^*$.
 The coefficient of the base element $[ \ur{I}, \ur{K_1+\{p+1\}}] \otimes [ \ur{J}, 
\ur{K_1+\{p\}}]^*$ is ${(-1)^{\vert \mathtt{W} - \mathtt{W} \cap \mathtt{K_1} \vert}}/{2^{(n - \vert \tJ \vert - \vert \tI \vert)/2}}$.

Case 1-$2'$. \quad  $p \notin \tK$.

We put $\tK_1 = \tK$.
Then $X_p ([ \ur{I}, \ur{K}] \otimes [ \ur{J}, \ur{K}]^*) =
X_p ([ \ur{I}, \ur{K_1}] \otimes [ \ur{J_1+\{p+1\}}, \ur{K_1}]^*)
 = [ \ur{I}, \ur{K_1}] \otimes [ \ur{J_1+\{p\}}, \ur{K_1}]^*$.
 Here we put $\tJ_1 = \tJ -\{p+1\}$.

The coefficient of $[ \ur{I}, \ur{K_1}] \otimes [ \ur{J_1+\{p\}}, \ur{K_1}]^*$ is  
${(-1)^{\vert \mathtt{W} - \mathtt{W} \cap \mathtt{K_1} \vert}}/{2^{(
n - \vert \tJ \vert - \vert \tI \vert)/2}}$.

So comparing the corresponding cases, we obtain  $X_p \phi_{\ell} = \phi_{\ell} X_p$ for this base element.

Case 2. \quad $p+1 \in \tJ$ and $p \in \tI$.

In this case, $p, p+1 \notin \tW$ and we put $\tJ_1 = \tJ - \{p+1\}$ and $\tI_1 = \tI - \{p\}$. 

Then 
$X_p <\ur{J_1+\{p+1\}}, \ur{W}, \ulb{W}, \ulb{I_1+\{p\}} > =
<\ur{J_1+\{p\}}, \ur{W}, \ulb{W}, \ulb{I_1+\{p\}} > - 
<\ur{J_1+\{p+1\}}, \ur{W}, \ulb{W}, \ulb{I_1+\{p+1\}} >$ 
and 
$<\ur{J_1+\{p\}}, \ur{W}, \ulb{W}, \ulb{I_1+\{p\}} >
={\ep} <\ur{J_1}, \ur{W+\{p\}}, \ulb{W+\{p\}}, \ulb{I_1} >$
and 
$<\ur{J_1+\{p+1\}}, \ur{W}, \ulb{W}, \ulb{I_1+\{p+1\}} >
={\ep} <\ur{J_1}, \ur{W+\{p+1\}}, \ulb{W+\{p+1\}}, \ulb{I_1} >$.

Here we put $\ep=(-1)^{ \# \{j \in {\tJ}_1 : j > p\} + \# \{i \in {\tI}_1 : i > p\}}=(-1)^{ \# \{j \in {\tJ}_1 : j > p+1\} + \# \{i \in {\tI}_1 : i > p+1\}}$.

Then we have 
\begin{equation*}
\begin{split}
\phi_{\ell} (X_p <\ur{J}, \ur{W}, \ulb{W}, \ulb{I} > )  
 =& \ep \sum_{\substack{ [1,n] - \mathtt{J_1} - \mathtt{I_1} \supseteq 
\mathtt{K}}}\dfrac{(-1)^{\vert (\mathtt{W} +\{p\}) - \mathtt{K} \cap (\mathtt{W} +\{p\})\vert}}{2^{(n - \vert \tJ_1 \vert - \vert \tI_1 \vert)/2}} [ \ur{I_1}, \ur{K}] \otimes [ \ur{J_1}, \ur{K}]^*  \\
 & - \ep \sum_{\substack{ [1,n] - \mathtt{J_1} - \mathtt{I_1} \supseteq 
\mathtt{K}}}\dfrac{(-1)^{\vert (\mathtt{W} +\{p+1\}) - \mathtt{K} \cap (\mathtt{W} +\{p+1\})\vert}}{2^{(
n - \vert \tJ_1 \vert - \vert \tI_1 \vert)/2}} [ \ur{I_1}, \ur{K}] \otimes [ \ur{J_1}, 
\ur{K}]^*). 
\end{split}
\end{equation*}

We have four cases here, but if  $p, p+1 \in \tK$, or if $p, p+1 \notin \tK$,
the two terms in the right-hand side are canceled out.
So we take two cases for $p \in \tK, p+1 \notin \tK$ and $p \notin \tK, p+1 \in 
\tK$.

Case 2-1. \quad $p \in \tK, p+1 \notin \tK$. 

We put $\tK_1 = \tK - \{p\}$.
Then the base element of the first term in the right-hand side is 
$ [ \ur{I_1}, \ur{K_1+\{p\}}] \otimes [ \ur{J_1}, \ur{K_1+\{p\}}]^* =
\ep [ \ur{I_1+\{p\}}, \ur{K_1}] \otimes [ \ur{J_1+\{p\}}, \ur{K_1}]^*$
 and the coefficient of $[ \ur{I_1+\{p\}}, \ur{K_1}] \otimes [ \ur{J_1+\{p\}}, 
\ur{K_1}]^*$ is 
${(-1)^{\vert \mathtt{W} - \mathtt{K_1} \cap \mathtt{W}\vert}}/{2^{(
n - \vert \tJ_1 \vert - \vert \tI_1 \vert)/2}}$. 

Also the base element of the second term in the right-hand side is 
$[ \ur{I_1+\{p\}}, \ur{K_1}] \otimes [ \ur{J_1+\{p\}}, \ur{K_1}]^*$ and its 
coefficient is 
${(-1)^{ \vert \mathtt{W} - \mathtt{K_1} \cap \mathtt{W}\vert}}/{2^{(
n - \vert \tJ_1 \vert - \vert \tI_1 \vert)/2}}
$. 
So the coefficient are the same and the total coefficient of 
$[ \ur{I_1+\{p\}}, \ur{K_1}] \otimes [ \ur{J_1+\{p\}}, \ur{K_1}]^*$ in this case is given by 
${(-1)^{\vert \mathtt{W} - \mathtt{K_1} \cap \mathtt{W}\vert}}/{2^{(
n -2 - \vert \tJ_1 \vert - \vert \tI_1 \vert)/2}}
={(-1)^{\vert \mathtt{W} - \mathtt{K_1} \cap \mathtt{W}\vert}}/{2^{(
n - \vert \tJ \vert - \vert \tI \vert)/2}}$.

Case 2-2. \quad $p \notin \tK, p+1 \in \tK$. 

We put $\tK_1 = \tK - \{p+1\}$.
By a similar calculation as above, 
 the total coefficient of 
$[ \ur{I_1+\{p+1\}}, \ur{K_1}] \otimes [ \ur{J_1+\{p+1\}}, \ur{K_1}]^*$ in this case is given by 
${(-1)^{1 + \vert \mathtt{W} - \mathtt{K_1} \cap \mathtt{W}\vert}}/{2^{(
n - \vert \tJ \vert - \vert \tI \vert)/2}}$.

Next we calculate $X_p \phi_{\ell}$. We have
\begin{align*}
& X_p \phi_{\ell}(<\ur{J}, \ur{W}, \ulb{W}, \ulb{I} >) 
 = \sum_{\substack{ [1,n] - \mathtt{J} - \mathtt{I} \supseteq \mathtt{K}}}
\dfrac{(-1)^{\vert \mathtt{W} - \mathtt{W} \cap \mathtt{K} \vert}}{2^{(
n - \vert \tJ \vert - \vert \tI \vert)/2}} \times \\
&(-[ \ur{I_1+\{p+1\}}, \ur{K}] \otimes 
[ \ur{J_1+\{p+1\}}, \ur{K}]^* + [ \ur{I_1+\{p\}}, \ur{K}] \otimes [ \ur{J_1+\{p\}}, 
\ur{K}]^*).
\\
\end{align*}

Since $p,p+1 \notin \tK$, we put $\tK =\tK_1$ and compare the coefficients in Case 2-1 and Case 2-2, we obtain  $X_p \phi_{\ell} = \phi_{\ell} X_p$ for this base element $<\ur{J}, \ur{W}, \ulb{W}, \ulb{I} >$.

Case 3. \quad $p+1 \in \tW$ and $p \notin \tI, \tJ, \tW$.

Then $p+1 \notin \tI$ and $p+1 \notin \tJ$ and we put $\tW_1 = \tW -\{p+1\}$.

Since 
$X_p <\ur{J}, \ur{W}, \ulb{W}, \ulb{I} >
={\ep} 
<\ur{J+\{p\}}, \ur{W_1}, \ulb{W_1}, \ulb{I+\{p+1\}} >$, 
we have 
\begin{align*}
&\phi_{\ell}(X_p <\ur{J}, \ur{W}, \ulb{W}, \ulb{I} >) \\
& = \ep
\sum_{\substack{ [1,n] - (\mathtt{J}+\{p\}) - (\mathtt{I}+\{p+1\}) \supseteq 
\mathtt{K}}}
\dfrac{(-1)^{\vert \mathtt{W_1} - \mathtt{W_1} \cap \mathtt{K} \vert}}{2^{(
n - \vert \tJ \vert - \vert \tI \vert-2)/2}} [ \ur{I+\{p+1\}}, \ur{K}] \otimes 
[ \ur{J+\{p\}}, \ur{K}]^* \\
& = 2
\sum_{\substack{ [1,n] - (\mathtt{J}+\{p\}) - (\mathtt{I}+\{p+1\}) \supseteq 
\mathtt{K}}}
\dfrac{(-1)^{\vert \mathtt{W_1} - \mathtt{W_1} \cap \mathtt{K} \vert}}{2^{(
n - \vert \tJ \vert - \vert \tI \vert)/2}} [ \ur{I}, \ur{K+\{p+1\}}] \otimes [ \ur{J}, \ur{K+\{p\}}]^*. \\
\end{align*}
Here we put $\ep =(-1)^{ \# \{j \in \tJ : j > p\} + \# \{i \in \tI : i > p+1\}}$.

We calculate $X_p \phi_{\ell}$. Since $\ur{W}=\ur{W_1}+\{p+1\}$, we have 
\begin{align*}
&\phi_{\ell}( <\ur{J}, \ur{W}, \ulb{W}, \ulb{I} >)  
= \sum_{\substack{ [1,n] - \mathtt{J} - \mathtt{I} \supseteq \mathtt{K}}}
\dfrac{(-1)^{\vert (\mathtt{W_1}+ \{p+1\}) - (\mathtt{W_1}+ \{p+1\}) \cap 
\mathtt{K} \vert}}{2^{(
n - \vert \tJ \vert - \vert \tI \vert)/2}} [ \ur{I}, \ur{K}] \otimes [ \ur{J}, 
\ur{K}]^*. \\
\end{align*}

If $p, p+1 \in \tK$ or if  $p, p+1 \notin \tK$, then the action of $X_p$ on the above element reduces to $0$.
So we have two cases for $p \in \tK, p+1 \notin \tK$
 and $p \notin \tK, p+1 \in \tK$.

Case 3-1. \quad $p \in \tK, p+1 \notin \tK$.

We put $\tK_1 = \tK -\{p\}$. Then we have 
\begin{align*}
&X_p \phi_{\ell}( <\ur{J}, \ur{W}, \ulb{W}, \ulb{I} >) \\
& = \sum_{\substack{ [1,n] - \mathtt{J} - \mathtt{I} \supseteq \mathtt{K}}}
\dfrac{(-1)^{\vert (\mathtt{W_1}+ \{p+1\}) - (\mathtt{W_1}+ \{p+1\}) \cap 
(\mathtt{K_1} + \{p\}) \vert}}{2^{(
n - \vert \tJ \vert - \vert \tI \vert)/2}} X_p([ \ur{I}, \ur{K_1+\{p\}}] \otimes 
[ \ur{J}, \ur{K_1+\{p\}}]^*) \\
& = \sum_{\substack{ [1,n] - \mathtt{J} - \mathtt{I} \supseteq \mathtt{K}}}
\dfrac{(-1)^{\vert \mathtt{W_1} - \mathtt{W_1} \cap \mathtt{K_1} \vert}}{2^{(
n - \vert \tJ \vert - \vert \tI \vert)/2}} [ \ur{I}, \ur{K_1+\{p+1\}}] \otimes [ \ur{J}, \ur{K_1+\{p\}}]^*. \\
\end{align*}

Case 3-2. \quad $p \notin \tK, p+1 \in \tK$.

We put $\tK_1 = \tK -\{p+1\}$.
\begin{equation*}
X_p \phi_{\ell}( <\ur{J}, \ur{W}, \ulb{W}, \ulb{I} >) 
 = \sum_{\substack{ [1,n] - \mathtt{J} - \mathtt{I} \supseteq \mathtt{K}}}
\dfrac{(-1)^{\vert \mathtt{W_1} - \mathtt{W_1} \cap \mathtt{K_1} \vert}}{2^{(
n - \vert \tJ \vert - \vert \tI \vert)/2}} [ \ur{I}, \ur{K_1+\{p+1\}}] \otimes [ \ur{J}, 
\ur{K_1+\{p\}}]^*. 
\end{equation*}

So by adding the Case 3-1 and Case 3-2, we have 
 $X_p \phi_{\ell} = \phi_{\ell} X_p$ for the base element $<\ur{J}, \ur{W}, \ulb{W}, \ulb{I} 
>$.

For $\frak{so}(2n+1)$, we must show that $X_p \phi_{\ell} = \phi_{\ell} X_p$ ($p = 1,2, \ldots, n-1$) for the base element $<\ur{J}, \ur{W}, \ur{0'}, \ulb{W}, \ulb{I} >$, but the proof is almost similar as above, so we omit it.

Next, we show that $X_n \phi_{\ell} = \phi_{\ell} X_n$ for $\frak{so}(2n+1)$.

Proof of $X_n \phi_{\ell} = \phi_{\ell} X_n$.

Case 1. \quad $n \in \tI$.

We put $\tI_1 = \tI - \{n\}$.
Then
$X_n <\ur{J}, \ur{W}, \ulb{W}, \ulb{I} > = -\sqrt{2} <\ur{J}, \ur{W}, \ulb{W}, \ur{0'}, \ulb{I_1} >$ and we have

\begin{equation*}
\phi_{\ell} (X_n <\ur{J}, \ur{W}, \ulb{W}, \ulb{I} > )  
=\sum_{\substack{ [1,n] - \mathtt{J} - \mathtt{I_1} \supseteq \mathtt{K}}}
\dfrac{(-1)^{\vert \mathtt{K} - \mathtt{K} \cap \mathtt{W} \vert + \vert \tW \vert +1}}{2^{(n - \vert \tJ \vert - \vert \tI \vert)/2}} [ \ur{I_1}, \ur{K}] \otimes [ \ur{J}, \ur{K}]^* 
\end{equation*}

We take the cases for $ n \in \tK$ and $n \notin \tK$.

Case 1-1. \quad $n \in \tK$.

We put $\tK_1 = \tK - \{n\}$.
Then in the right-hand side,  the base element
 is $[ \ur{I_1}, \ur{K_1},\ur{n}] \otimes [ \ur{J}, \ur{K_1},\ur{n}]^*$
and its coefficient is ${(-1)^{\vert \mathtt{K_1} - \mathtt{K_1} \cap \mathtt{W} \vert + \vert \tW \vert}}/{2^{(n - \vert \tJ \vert - \vert \tI \vert)/2}}$.

Case 1-2. \quad $n \notin \tK$.

We put $\tK_1 = \tK$.
Then in the right-hand side,  the base element
 is $[ \ur{I_1}, \ur{K_1}] \otimes [ \ur{J}, \ur{K_1}]^*$
and its coefficient is ${(-1)^{\vert \mathtt{K_1} - \mathtt{K_1} \cap \mathtt{W} \vert + \vert \tW \vert+1}}/{2^{(n - \vert \tJ \vert - \vert \tI \vert)/2}}$.

Next, we calculate $X_n \phi_{\ell}$.

\begin{equation*}
 \phi_{\ell}(<\ur{J}, \ur{W}, \ulb{W}, \ulb{I} >)
 = \sum_{\substack{ [1,n] - \mathtt{J} - \mathtt{I} \supseteq \mathtt{K}}}
\dfrac{(-1)^{\vert \mathtt{W} - \mathtt{W} \cap \mathtt{K} \vert}}{2^{(
n - \vert \tJ \vert - \vert \tI \vert)/2}} [ \ur{I}, \ur{K}] \otimes [ \ur{J}, 
\ur{K}]^*
 \end{equation*}
 and since $n \in \tI$, $n \notin \tJ, \tK$. If we write $\tK =\tK_1$, we have 
\begin{align*}
 & X_n \phi_{\ell}(<\ur{J}, \ur{W}, \ulb{W}, \ulb{I} >) \\
& = \sum_{\substack{ [1,n] - \mathtt{J} - \mathtt{I} \supseteq \mathtt{K_1}}}
\dfrac{(-1)^{\vert \mathtt{W} - \mathtt{W} \cap \mathtt{K_1} \vert}}{2^{(
n - \vert \tJ \vert - \vert \tI \vert)/2}} 
((-1)^{\vert \tK_1 \vert +1}[ \ur{I_1}, \ur{K_1}] \otimes [ \ur{J}, \ur{K_1}]^*
 +[ \ur{I}, \ur{K_1}] \otimes [ \ur{J}, \ur{K_1}, \ur{n}]^*).\\
 \end{align*}

Here we note that $X_n [ \ur{I}, \ur{K_1}] = X_n [ \ur{I_1}, \ur{n}, \ur{K_1}] = (-1)^{\vert \tK_1 \vert} X_n [ \ur{I_1},  \ur{K_1}, \ur{n}] =(-1)^{\vert \tK_1 \vert +1}[ \ur{I_1}, \ur{K_1}]$.

So the coefficient of $[ \ur{I_1}, \ur{K_1}] \otimes [ \ur{J}, \ur{K_1}]^*$ is equal to $\dfrac{(-1)^{\vert \mathtt{W} - \mathtt{W} \cap \mathtt{K_1} \vert+\vert \tK_1 \vert +1}}{2^{(n - \vert \tJ \vert - \vert \tI \vert)/2}}$ ,which is equal to that of Case 1-2.
Since  $[ \ur{I}, \ur{K_1}] \otimes [ \ur{J}, \ur{K_1}, \ur{n}]^*
= [ \ur{I_1}, \ur{n}, \ur{K_1}] \otimes [ \ur{J}, \ur{K_1}, \ur{n}]^*
= (-1)^{\vert \tK_1 \vert} [ \ur{I_1}, \ur{K_1}, \ur{n}] \otimes [ \ur{J}, \ur{K_1}, \ur{n}]^*$, 
the coefficient of $[ \ur{I_1}, \ur{K_1}, \ur{n}] \otimes [ \ur{J}, \ur{K_1}, \ur{n}]^*$ is given by 
$
\dfrac{(-1)^{\vert \mathtt{W} - \mathtt{W} \cap \mathtt{K_1} \vert + \vert \tK_1 \vert}}{2^{(n - \vert \tJ \vert - \vert \tI \vert)/2}}$, which is equal to that of Case 1-1.So we have done the verification in this case.

Next we check the equation for the base element $<\ur{J}, \ur{W}, \ur{0'}, \ulb{W}, \ulb{I} >$. (We are still in the Case 1 and assume $\tI \ni n$.)

Since $X_n <\ur{J}, \ur{W}, \ur{0'}, \ulb{W}, \ulb{I} > = \sqrt{2} 
<\ur{J}, \ur{W}, \ur{n}, \ulb{W}, \ulb{I} >
=\sqrt{2} (-1)^{\vert \tW \vert} 
<\ur{J}, \ur{W+\{n\}},  \ulb{W+\{n\}}, \ulb{I_1} >
$,
we have 
\begin{equation*}
\phi_{\ell}(X_n <\ur{J}, \ur{W}, \ur{0'}, \ulb{W}, \ulb{I} >)
=\sqrt{2} (-1)^{\vert \tW \vert}
\sum_{\substack{ [1,n] - \mathtt{J} - \mathtt{I_1} \supseteq \mathtt{K}}}
\dfrac{(-1)^{\vert \mathtt{W}+ \{n\} - (\mathtt{W}+ \{n\}) \cap \mathtt{K} 
\vert}}{2^{(n - \vert \tJ \vert - \vert \tI \vert+1)/2}} [ \ur{I_1}, \ur{K}] \otimes [ \ur{J}, \ur{K}]^*.
 \end{equation*}

We take the cases for $n \in \tK$ and $n \notin \tK$.

Case 1-3. \quad $n \in \tK$. 

We put $\tK_1 = \tK - \{n\}$.
Then in the right-hand side,  the base element
 is $[ \ur{I_1}, \ur{K_1},\ur{n}] \otimes [ \ur{J}, \ur{K_1},\ur{n}]^*$
and its coefficient is ${(-1)^{\vert \mathtt{W} - \mathtt{K_1} \cap \mathtt{W} \vert + \vert \tW \vert}}/{2^{(
n - \vert \tJ \vert - \vert \tI \vert)/2}}$.

Case 1-4. \quad $n \notin \tK$.

We put $\tK_1 = \tK$.
Then in the right-hand side,  the base element
 is $[ \ur{I_1}, \ur{K_1}] \otimes [ \ur{J}, \ur{K_1}]^*$
and its coefficient is ${(-1)^{\vert \mathtt{W} - \mathtt{K_1} \cap \mathtt{W} \vert + \vert \tW \vert+1}}/{2^{(n - \vert \tJ \vert - \vert \tI \vert)/2}}$.

Next we calculate $X_n \phi_{\ell}$.
Since $n \in \tI$, $n \notin \tJ, \tK$, we write $\tK =\tK_1$ and we have 
\begin{align*}
& X_n \phi_{\ell}(<\ur{J}, \ur{W}, \ur{0'}, \ulb{W}, \ulb{I} >) \\
& = X_n \sum_{\substack{ [1,n] - \mathtt{J} - \mathtt{I} \supseteq \mathtt{K}}}
\dfrac{(-1)^{\vert \mathtt{K_1} - \mathtt{K_1} \cap \mathtt{W} \vert}}{2^{(
n - \vert \tJ \vert - \vert \tI \vert)/2}} [ \ur{I}, \ur{K_1}] \otimes [ \ur{J}, \ur{K_1}]^* \\
& = \sum_{\substack{ [1,n] - \mathtt{J} - \mathtt{I} \supseteq \mathtt{K_1}}}
\dfrac{(-1)^{\vert \mathtt{K_1} - \mathtt{K_1} \cap \mathtt{W} \vert}}{2^{(
n - \vert \tJ \vert - \vert \tI \vert)/2}}
((-1)^{\vert \tK_1 \vert +1}[ \ur{I_1}, \ur{K_1}] \otimes [ \ur{J}, \ur{K_1}]^*
 +[ \ur{I}, \ur{K_1}] \otimes [ \ur{J}, \ur{K_1}, \ur{n}]^*).
\end{align*}

So the coefficient of $[ \ur{I_1}, \ur{K_1}] \otimes [ \ur{J}, \ur{K_1}]^*$ is equal to ${(-1)^{\vert \mathtt{K_1} - \mathtt{W} \cap \mathtt{K_1} \vert+\vert \tK_1 \vert +1}}/{2^{(n - \vert \tJ \vert - \vert \tI \vert)/2}}$, 
 which is equal to that of Case 1-4.

As before,  $[ \ur{I}, \ur{K_1}] \otimes [ \ur{J}, \ur{K_1}, \ur{n}]^*
= [ \ur{I_1}, \ur{n}, \ur{K_1}] \otimes [ \ur{J}, \ur{K_1}, \ur{n}]^*
= (-1)^{\vert \tK_1 \vert} [ \ur{I_1}, \ur{K_1}, \ur{n}] \otimes [ \ur{J}, \ur{K_1}, \ur{n}]^*$, 
the coefficient of $[ \ur{I_1}, \ur{K_1}, \ur{n}] \otimes [ \ur{J}, \ur{K_1}, \ur{n}]^*$ is given by 
${(-1)^{\vert \mathtt{K_1} - \mathtt{W} \cap \mathtt{K_1} \vert + \vert \tK_1 
\vert}}/{2^{(n - \vert \tJ \vert - \vert \tI \vert)/2}}$, which is equal to that of Case 1-3.

Case 2. \quad $n \in \tW$.

Then $n \notin \tI$ and $n \notin \tJ$  and we put $\tW_1 = \tW - \{n\}$.

We calculate $\phi_{\ell} X_n$. Since 
$X_n <\ur{J}, \ur{W_1+\{n\}}, \ulb{\{n\}+ W}, \ulb{I} > = -\sqrt{2} <\ur{J}, 
\ur{W_1+\{n\}},  \ur{0'}, \ulb{W_1}, \ulb{I} > =
(-1)^{\vert \tW \vert} \sqrt{2} <\ur{J+\{n\}}, \ur{W_1},  \ur{0'}, \ulb{W_1}, \ulb{I} >$, 
we have 

\begin{align*}
&\phi_{\ell} (X_n <\ur{J}, \ur{W}, \ulb{W}, \ulb{I} > )  
= 2 \sum_{\substack{ [1,n] - \mathtt{J_1} - \mathtt{I} \supseteq 
\mathtt{K}}}\dfrac{(-1)^{\vert \mathtt{K} - \mathtt{K} \cap \mathtt{W_1} \vert + \vert \tW \vert}}{2^{(n - \vert \tJ \vert - \vert \tI \vert)/2}} [ \ur{I}, \ur{K}] \otimes [ \ur{J_1}, \ur{K}]^*. \\
\end{align*}

Here we put $\tJ_1 = \tJ + \{n\}$.
Since $n \in \tJ_1$, $n \notin \tK$ and we put $\tK_1 =\tK$.

Then in the right-hand side,  the base element
 is $[ \ur{I}, \ur{K_1}] \otimes [ \ur{J_1}, \ur{K_1}]^*
 = [ \ur{I}, \ur{K_1}] \otimes [ \ur{J+\{n\}}, \ur{K_1}]^*
 = (-1)^{\vert \tK_1 \vert}  [ \ur{I}, \ur{K_1}] \otimes [ \ur{J}, 
\ur{K_1+\{\ur{n}\}}]^*$
 and the coefficient of $[ \ur{I}, \ur{K_1}] \otimes [ \ur{J}, \ur{K_1+\{\ur{n}\}}]^*$ is 
${(-1)^{\vert \mathtt{K_1} - \mathtt{K_1} \cap \mathtt{W_1} \vert + \vert \tW \vert + \vert \tK_1 \vert}}/{2^{(n-2 - \vert \tJ \vert - \vert \tI \vert)/2}}
= {(-1)^{\vert \mathtt{W} - \mathtt{K_1} \cap \mathtt{W} \vert}}/{2^{(
n -2 - \vert \tJ \vert - \vert \tI \vert)/2}}$.

We calculate $X_n \phi_{\ell}$.
 Since
\begin{equation*}
 \phi_{\ell}(<\ur{J}, \ur{W}, \ulb{W}, \ulb{I} >)
 = \sum_{\substack{ [1,n] - \mathtt{J} - \mathtt{I} \supseteq \mathtt{K}}}
\dfrac{(-1)^{\vert \mathtt{W} - \mathtt{W} \cap \mathtt{K} \vert}}{2^{(
n - \vert \tJ \vert - \vert \tI \vert)/2}} [ \ur{I}, \ur{K}] \otimes [ \ur{J}, 
\ur{K}]^*
 \end{equation*}
and $n \in \tI$, $n \notin \tJ$, 
we take two cases for $ n \in \tK$ and $ n \notin \tK$.

Case 2-1. \quad $n \in \tK$.

We put $\tK_1 =\tK -\{n\}$ and we have 
\begin{equation*}
 X_n \phi_{\ell}(<\ur{J}, \ur{W}, \ulb{W}, \ulb{I} >)
 = \sum_{\substack{ [1,n] - \mathtt{J} - \mathtt{I} \supseteq \mathtt{K}}}
\dfrac{(-1)^{\vert \mathtt{W} - \mathtt{W} \cap \mathtt{K} \vert}}{2^{(
n - \vert \tJ \vert - \vert \tI \vert)/2}} 
(-[ \ur{I_1}, \ur{K_1}] \otimes [ \ur{J}, \ur{K_1+\{n\}}]^*).
 \end{equation*}

So the coefficient of the base element $[ \ur{I_1}, \ur{K_1}] \otimes [ \ur{J}, \ur{K_1+\{n\}}]^*$ is $-{(-1)^{\vert \mathtt{W} - \mathtt{W} \cap \mathtt{K} 
\vert}}/{2^{(n - \vert \tJ \vert - \vert \tI \vert)/2}} 
 = {(-1)^{\vert \mathtt{W} - \mathtt{W} \cap \mathtt{K_1} \vert}}/{2^{(n - \vert \tJ \vert - \vert \tI \vert)/2}}$.

Case 2-2. \quad $n \notin \tK$.

We put $\tK_1 =\tK$ and we have 
\begin{equation*}
 X_n \phi_{\ell}(<\ur{J}, \ur{W}, \ulb{W}, \ulb{I} >)
 = \sum_{\substack{ [1,n] - \mathtt{J} - \mathtt{I} \supseteq \mathtt{K_1}}}
\dfrac{(-1)^{\vert \mathtt{W} - \mathtt{W} \cap \mathtt{K_1} \vert}}{2^{(
n - \vert \tJ \vert - \vert \tI \vert)/2}} 
([ \ur{I_1}, \ur{K_1}] \otimes [ \ur{J}, \ur{K_1+\{\ur{n}\}}]^*).
 \end{equation*}

So the coefficient of the base element $[ \ur{I_1}, \ur{K_1}] \otimes [ \ur{J}, \ur{K_1+\{n\}}]^*$ is ${(-1)^{\vert \mathtt{W} - \mathtt{W} \cap \mathtt{K_1} 
\vert}}/{2^{(n - \vert \tJ \vert - \vert \tI \vert)/2}}.$

By adding the results of Case 2-1 and Case 2-2,  we have $ \phi_{\ell} X_n = X_n \phi_{\ell}$ in this case.

There are several remaining cases, in which we must show $ \phi_{\ell} X_n = X_n \phi_{\ell}$, but the proof is almost similar, so we omit it.

We must show that $ \phi_{\ell}^{\ep_1 \ep_2} X_n = X_n \phi_{\ell}^{\ep_1 \ep_2}$ and $ \phi_{n}^{even} X_n = X_n \phi_{n}^{even}$, but the proofs go in a similar way by (lengthy) case-by-case verifications  and so we omit it.

\end{proof}

\section{Invariant theory and the centralizer algebras }\label{sect inv}

In this section we define a set of elements based on the invariant theory in the centralizer algebra $\mathbf{CP_k}= \Hom_{Pin(N)}(\Delta \bigotimes \bigotimes^k V, \Delta \bigotimes \bigotimes^k V )$.
 They are parameterized by the generalized Brauer diagrams.
In the classical situation for the orthogonal group $O(N)$, R. Brauer (\cite{br}) considered the centralizer algebra $\Hom_{O(N)}( \bigotimes^k V,  \bigotimes^k V )$ and he determined the special basis (parameterized by the Brauer diagrams) and the multiplication rules of the basis elements by using the first and the second fundamental theorem of the invariant theory of $O(N)$. 

Since we have already the explicit isomorphisms between the exterior algebras $\bigwedge^i V$ and the irreducible components of $\Delta \bigotimes \Delta^*$
as $Pin(N)$ ($O(N)$) modules in the previous section, we can follow the above way.

More generally, we define the spaces of the equivariant homomorphisms.

\begin{defn}\label{def cpkl}
$$
\mathbf{CP^k_l}= \Hom_{Pin(N)}(\Delta \textstyle{\bigotimes \bigotimes\limits^k V, \Delta \bigotimes \bigotimes\limits^l V} ).
$$ 
and 
$$
\mathbf{CS^k_l}= \Hom_{Spin(N)}(\Delta \textstyle{\bigotimes \bigotimes\limits^k V, \Delta \bigotimes \bigotimes\limits^l V} ).
$$ 
\end{defn}

We write $\mathbf{CP_k}$ and $\mathbf{CS_k}$ for $\mathbf{CP^k_k}$ and 
$\mathbf{CP^k_k}$ respectively.

\begin{rem}
We note that for $Pin(2n+1)$, $\Delta$ stands for $\Delta_{\pm}$, so there are four cases for $\mathbf{CP^k_l}$. $($But two of them are trivial. See Lemma \ref{lem odd0}.$)$ 
For $Spin(2n)$,  we can change $\Delta$ into $\Delta^{\pm}$ in the above and so there are also four cases for $\mathbf{CS^k_l}$.

We note that from the formulas \eqref{eqn evendelnat}, \eqref{eqn evendelnnat}, for $\ep_1, \ep_2 \in \{\pm 1\}$, the irreducible representation $(1/2 + \delta)^{\ep_1}_{Spin(2n)}$ of $Spin(2n)$ occurs in $\Delta^{\ep_2} \bigotimes \bigotimes^k V$ if and only if  $(-1)^{k - \vert \delta \vert} = \ep_1 \ep_2$. 

So we have $\Hom_{Spin(2n)}(\Delta^{\ep_1} \bigotimes \bigotimes^k V, 
\Delta^{\ep_2} \bigotimes \bigotimes^l V )=0$ if $\ep_1 \ep_2 \ne (-1)^{(k-l)}$.\end{rem}

We interpret them as the elements of the invariant subspace of the tensor spaces under the action of $Pin(N)$.

From Theorem \ref{thm phipinodd} and Theorem \ref{thm phipin} and $V \cong V^*$, we have the following lemma.
\begin{lem}\label{lem equiv}
For $Pin(2n)$, we have 
\begin{equation}
 \Hom_{Pin(2n)}(\textstyle{\Delta \bigotimes  \bigotimes\limits^k V, \Delta 
\bigotimes \bigotimes\limits^l V) 
 \cong 
\bigoplus_{i=0}^{2n} ((\bigwedge\limits^{i} V) \bigotimes  
\bigotimes\limits^{k+l} V})^{O(2n)}.
\tag{\ref{lem equiv}.1}\label{inv equiv1}
\end{equation}

For $Pin(2n+1)$, we have 
\begin{equation}
 \Hom_{Pin(2n+1)}(\Delta_{\pm} \textstyle{\bigotimes  \bigotimes\limits^k V, 
\Delta_{\pm} \bigotimes  \bigotimes\limits^l V ) 
 \cong 
\bigoplus\limits_{i=0}^{n} ((\bigwedge\limits^{2i} V) \bigotimes  
\bigotimes\limits^{k+l} V})^{O(2n+1)}. 
\tag{\ref{lem equiv}.2}\label{inv equiv2}
\end{equation}
and 
\begin{equation}
\Hom_{Pin(2n+1)}(\Delta_{\pm} \textstyle{\bigotimes  \bigotimes\limits^k V, 
\Delta_{\mp} \bigotimes  \bigotimes\limits^l V ) 
 \cong 
\bigoplus_{i=0}^{n} ((\bigwedge\limits^{2i+1} V) \bigotimes   
\bigotimes\limits^{k+l} V})^{O(2n+1)}. 
\tag{\ref{lem equiv}.3}\label{inv equiv3}
\end{equation}

Here the superscript $O(N)$ denotes the 
subspace obtained by taking the $O(N)$-invariants.
\end{lem}

First let us recall the first main theorem and the second main theorem for the  polynomial invariants under the action of the orthogonal groups.

Let $V={K}^N$ be a vector space over a field $K$ of characteristic $0$ and let 
$P(\bigoplus^k V)$ be the polynomial ring over the vector space  $\bigoplus^k V$. Let $(\quad, \quad)$ be a non-degenerate symmetric bilinear form which defines the orthogonal group $O(N)$.
Then \lq The First Main Theorem' (see \cite{weyl} p53 Theorem 2.9A, p64 Theorem 2.11A.) tells us that the polynomial invariants $P(\bigoplus^k V)^{O(N)}$  are 
generated by the scalar products $(\bfv{i}, \bfv{j} )$ ($1 \leqq i, j \leqq k$), where $\bfv{i}$ denotes the vector in the $i$-th direct summand in $\bigoplus^k V$. 
\lq The Second Main Theorem' (see \cite{weyl} p75 Theorem 2.17A .) 
tells us that the relations of the polynomial invariants $(\bfv{i},\bfv{j})$'s in $P(\bigoplus^k V)^{O(N)}$  
are generated by the determinants $\det((\bfv{i_s}, \bfv{j_t}))$ ($0 \leqq s, t \leqq N $) of size $N+1$.
We note that the minimum (total) degree of the relations is $2N+2$.

We consider the invariant tensors $(\bigwedge^{r} V^* \bigotimes^{s} V^*)^{O(N)}$.  
This space is included in the tensor space $(\bigotimes^{r+s} V^*)^{O(N)}$ and we regard this space as the subspace of $P(\bigoplus^{r+s} V )^{O(N)}$.
The elements in this space are multi-linear (i.e., degree $1$) in each of the tensor components and alternating concerning the first $r$ components.
Since the invariants are generated by the scalar products of degree $2$, the multi-linear property tells us that $(\bigotimes^{r+s} V^*)^{O(N)} =0$ if $r + s \equiv 1 \mod{2}$.

So we have the following.
\begin{lem}\label{lem odd0}
If $(-1)^{k + l} \ep_1 \ep_2= -1$, then 
$$\mathbf{CP^k_l} =\Hom_{Pin(2n+1)}(\Delta_{\ep_1} \textstyle{\bigotimes  
\bigotimes\limits^k V, \Delta_{\ep_2} \bigotimes  \bigotimes\limits^l V} ) =0.$$\end{lem}

\begin{rem}\label{rem odd0}
From the above lemma,  we only consider the case $(-1)^{k + l} \ep_1 \ep_2= 1$ for the centralizer algebras of $Pin(2n+1)$. If we write simply $\Delta$ in the formulas for $Pin(2n+1)$, we always assume this condition without mentioning it.
\end{rem}

We can write down the elements in $(\bigwedge^{r} V^* \bigotimes^{s} V^*)^{O(N)}$ as follows.
If $ s < r$, $\bigotimes^{s} V$ does not contain the representations $\bigwedge^{r} V$, so this space must be $0$. 

Let $s \geqq r$. 
For the subsets $\mathbf{t}= \{ {t_1}, {t_2}, \ldots, {t_{r}}\}$ (${t_1} <{t_2} < \ldots < {t_{r}}$ and $r \leqq N$) of $[1,s]$ and the sets of $u=\dfrac{s-r}{2}$ pairs of indices $\{\mathbf{m},\mathbf{l}\}= \{ \{m_1,l_1\}, \{m_2,l_2\},\ldots, \{m_{u},l_{u}\} \}$  such that $[1,s] = \mathbf{t} \sqcup \mathbf{m} \sqcup \mathbf{l}$ (disjoint unions), where we put $\mathbf{m} = \{m_1, \ldots, m_{u}\}$ and $\mathbf{l}=\{l_1, \ldots, l_{u}\}$.), 
we define the invariant elements $\mathbf{T_{t, \{{m},{l}\}}}$  by 
$$
\mathbf{T_{t, \{{m},{l}\}}} = \dfrac{1}{r!} \det((\bfx{i}, \bfy{t_j}))  \times \prod_{j=1}^{u} (\bfy{m_j}, \bfy{l_j}).
$$

Here $\bfx{j}$ ($ j= 1,2 \ldots,r$) are the components of the first $r$ tensors in  $\bigotimes^{r+s} V$ and $\bfy{j}$ ($ j= 1,2 \ldots,s$) are the components of the latter $s$ tensors in  $\bigotimes^{r+s} V$.

Moreover 
$\dfrac{1}{r!} \det((\bfx{i}, \bfy{t_j}))$  
is the invariant element in 
$((\bigwedge^{r} V)^* \bigotimes (\bigwedge^{r} V)^*)^{O(N)}$.

Since $\bigwedge^{r} V$ is an irreducible $O(N)$ module,
the trivial representation occurs in $(\bigwedge^{r} V)^* \bigotimes \bigwedge^{r} V$ with multiplicity one. The invariant element is given by

\begin{equation}
\begin{split}
&\frak{I}_r = \sum_{\substack{[1,n] \supseteq \tJ \sqcup \tI \sqcup \tW \: \text {(disjoint sets)} \\
 \vert \tJ \vert + \vert \tI \vert + 2 \vert \tW \vert =r}} 
 {<\ur{J}, \ur{W}, \ulb{W}, \ulb{I} >}^* \otimes <\ur{J}, \ur{W}, \ulb{W}, \ulb{I} > \\
&(+ \sum_{\substack{[1,n] \supseteq \tJ \sqcup \tI \sqcup \tW \: \text {(disjoint sets)} \\
 \vert \tJ \vert + \vert \tI \vert + 2 \vert \tW \vert = r -1}} 
 {<\ur{J}, \ur{W}, \ur{0'}, \ulb{W}, \ulb{I} >}^* \otimes <\ur{J}, \ur{W}, \ur{0'}, \ulb{W}, \ulb{I} >).\\
\label{eqn ir}
\end{split}
\end{equation}

Here the second term in the parentheses are considered only for $O(2n+1)$.

So $\dfrac{1}{r!} \det((\bfx{i}, \bfy{t_j}))$  
 must be a scalar multiple of $(\id \otimes \iota) (\frak{I}_{r})$ (see Definition \ref{def iota}.) 
and the scalar is given by $(-1)^{\binom{r}{2}}$.
(To obtain the scalar, take cases of the parity of $r$ and compare the values of  both sides at the special elements $\bfx{1} = \ua{1}$, $\bfx{2} = \ua{2}$, etc..)

\begin{lem} \label{lem indep}
For $Pin(2n)$, if we take only $\mathbf{t}$'s such that $\vert \mathbf{t} \vert \equiv k+l \mod{2}$ and $\vert \mathbf{t} \vert \leqq 2n$,
 then the above $\mathbf{T_{t,  \{{m},{l}\}}}$'s  span linearly the whole space of the $\mathbf{CP^k_l}= \Hom_{Pin(2n)}(\Delta \bigotimes \bigotimes^k V, \Delta \bigotimes \bigotimes^l V )$.

Moreover if $k+l < 2n+2 $, 
$\mathbf{T_{t,  \{{m},{l}\}}}$ $($ $\mathbf{t} \sqcup \mathbf{m} \sqcup \mathbf{l} =[1,s]$ $)$ are linearly independent. 

For $Pin(2n+1)$, if
 we take only $\mathbf{t}$'s such that $\vert \mathbf{t} \vert \leqq 2n+1$ and $\vert \mathbf{t} \vert \equiv k+l \equiv 0 \mod{2}$, 
 then the above $\mathbf{T_{t,  \{{m},{l}\}}}$'s  span linearly the whole space of the $\mathbf{CP^k_l}= \Hom_{Pin(2n+1)}(\Delta_{\pm} \bigotimes \bigotimes^k V, \Delta_{\pm} \bigotimes \bigotimes^l V )$.

For $Pin(2n+1)$, if
 we take only $\mathbf{t}$'s such that $\vert \mathbf{t} \vert \leqq 2n+1$ and $\vert \mathbf{t} \vert \equiv k+l \equiv 1 \mod{2} $, 
 then the above $\mathbf{T_{t,  \{{m},{l}\}}}$'s  span linearly the whole space of the $\mathbf{CP^k_l}= \Hom_{Pin(2n+1)}(\Delta_{\pm} \bigotimes \bigotimes^k V, \Delta_{\mp} \bigotimes \bigotimes^l V )$.

In both cases, if $k+l < 2n + 3 $, 
$\mathbf{T_{t,  \{{m},{l}\}}}$\, $(\mathbf{t} \sqcup \mathbf{m} \sqcup \mathbf{l} =[1,s])$ are linearly independent. 
\end{lem}

\begin{proof}
The statements follow from the argument developed above and the Second Main 
Theorem.
\end{proof}

From Lemma \ref{lem equiv}, 
we consider this invariant element $\mathbf{T_{t,  \{{m},{l}\}}}$
as a homomorphism from $\Delta \bigotimes \bigotimes^{k} V$ to $\Delta \bigotimes \bigotimes^{l} V$.

By the original definition, 
$\mathbf{T_{t,  \{{m},{l}\}}}$ are the elements in  $(\bigwedge^{r} V^*) \bigotimes \bigotimes^k V^*  \bigotimes \bigotimes^l V^*$ for $r \equiv k +l \mod{2}$.
 We denote the degree 2 tensor corresponding to the invariant symmetric bilinear form defined by $S$ on $V \times V$ by $\id_{V}$.

Then as the homomorphism in $\Hom_{Pin(N)}(\Delta \bigotimes \bigotimes^k V, \Delta \bigotimes \bigotimes^l V )$,  ${\mathbf{T_{t,  \{{m},{l}\}}}}$
 is given as follows.

The pair $(\bfy{m_j}, \bfy{l_j})$ with $m_j, l_j \leqq k$ corresponds to the 
contraction of the specified tensor components (i.e., $m_j$-th and $l_j$-th 
components) in $\Delta \bigotimes \bigotimes^k V$.  We denote this contraction by $\cont{\{\ur{m_j,l_j}\}}$.

The pair $(\bfy{m_j}, \bfy{l_j})$ with $m_j, l_j > k$ corresponds to the insertion of the invariant tensor $(\id_{V})$ into the specified tensor components (i.e., $m_j-k$-th and $l_j-k$-th components) in $\Delta \bigotimes \bigotimes^l V$.  We denote this insertion by ${\id_{V}}_{\ur{m_j-k,l_j-k}}$.
(We introduce the numbering such that the first $k$ tensors correspond to the domain and the last $l$ tensors correspond to the range in $\Hom_{Pin(N)}(\Delta \bigotimes \bigotimes^k V, \Delta \bigotimes \bigotimes^l V )$.)

The pair $(\bfy{m_j}, \bfy{l_j})$ with $m_j \leqq k$ and $l_j > k$ (resp. with 
$l_j \leqq k$ and $m_j > k$) corresponds to the map which transfers the $m_j$-th (resp. $l_j$-th ) component of the domain tensor  $\Delta \bigotimes \bigotimes^k V$ into $l_j-k$-th (resp. $m_j-k$-th ) component of the image tensor $\Delta \bigotimes \bigotimes^l V$.  Namely this corresponds to the partial permutation of the tensor components.

The above three operators are the same in the usual Brauer centralizer algebras.
For the remaining part 
 $\det( (\bfx{i}, \bfy{t_j}))$ ($1 \leqq i,j \leqq r$)  in  $\mathbf{T_{t,  
\{{m},{l}\}}}$, 
we will give the explicit description as the homomorphism.(See Theorem \ref{thm psi}.)

Let us use a diagrammatic parametrization  of the elements $\{\mathbf{T_{t,  
\{{m},{l}\}}}\}$ of $\Hom_{Pin(N)}(\Delta \bigotimes \bigotimes^k V, \Delta 
\bigotimes \bigotimes^l V )$ as follows.
We consider the diagrams of two lines vertices, with $k$ vertices in  the upper row and $l$ vertices in the lower row, in which vertices are connected with each other as in the usual Brauer diagrams except for admitting isolated vertices. Namely they are graphs with no loops, in which the number of edges connected to each vertex is either $0$ or $1$.

We call such diagrams the generalized Brauer diagrams and denote the set of  the generalized Brauer diagrams by $\mathbf{GB^k_l}$.
(Later we use another parametrization using the same diagrams $\mathbf{GB^k_l}$ from the representation theoretical viewpoint.)

The explicit bijection is given as follows.
Let us number the vertices from $1$ to $k$ from left to right in the upper row and from $1$ to $l$ in the lower row.

\begin{figure}[htbp]
\begin{minipage}{\linewidth}
\begin{center}
\input{spincentfig2.tex}
\end{center}
\caption{An example of the generalized Brauer diagrams with $k=6$ and $l=5$}
\end{minipage}
\label{spincentfig2}
\end{figure}
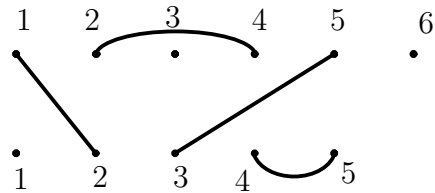

 We correspond each of ${\mathbf{T_{t,  \{{m},{l}\}}}}$ to a generalized Brauer 
diagram as follows. We consider the tensor components in the range $\Delta \bigotimes \bigotimes^l V $ are indexed by $k+1, k+2, \ldots, k+l$ in this argument.  
For the vertices with indices in $\mathbf{t}$, we leave them the isolated vertices and for the other vertices ( namely the vertices with indices in $\mathbf{l}$ and $\mathbf{m}$ ), we connect $l_j$-th vertex and $m_j$-th vertex. 
Except for the isolated vertices, the edges in the diagrams convey the same meaning as in the usual Brauer diagrams.

So we obtain the bijective correspondence between the elements $\{\mathbf{T_{t, 
\{{m},{l}\}}}\}$ of $\mathbf{CP^k_l}$ and the diagrams with the isolated vertices of at most $N$ in $\mathbf{GB^k_l}$.

From the formula \eqref{eqn odddelnat} in Theorem \ref{thm repodd} and \eqref{eqn evendelnat} and \eqref {eqn evendelnnat} in Theorem \ref{thm repeven}, the 
multiplicity of $[\Delta,\lam]_{Pin(N)}$ in $\Delta \bigotimes \bigotimes^k V$ is equal to the number of the \lq up-down-stable ($n$-) tableaux' of degree $k$  starting with the empty set $\emptyset$ and ending with the Young diagram $\lam$ and consisting of the diagrams  with their lengths at most $n$.
We show an example of up-down-stable tableaux.

\begin{figure}[htbp]
\begin{center}
\begin{minipage}{\linewidth}
\input{spincentfig3.tex}
\caption{An example of \lq up-down-stable tableaux' of $\lam =(2,2)$ and $k=8$}
\label{spincentfig3}
\end{minipage}
\end{center}
\end{figure}
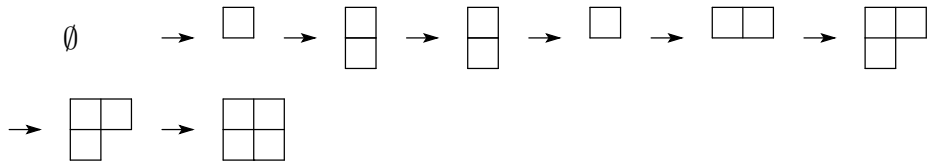

Namely by definition, the \lq up-down-stable ($n$-)tableaux' of degree $k$ are the $k+1$ sequences of Young diagrams, in which the adjacent diagrams are obtained by attaching a cell to or removing a cell from the previous one or are the same (if its length is less than $n$ for $N=2n$).

So figure \ref{spincentfig3} is not an \lq up-down-stable $2$-tableau' for $Pin(4)$ because there exist adjacent diagrams  which are the same and its length is $2$.

Since  $\mathbf{CP_k}$ and $Pin(N)$ are the dual pair on the space $\Delta \bigotimes \bigotimes^k V$, the irreducible representations of $\mathbf{CP_k}$ are in bijective correspondence with the irreducible representations of $Pin(N)$ on this space.

So we can expect the base of the irreducible representation $\lam$ of 
$\mathbf{CP_k}$ (corresponding to $[\Delta, \lam]_{Pin(N)}$)  are parameterized by the \lq up-down-stable ($n$-) tableaux' of degree $k$, starting with $\emptyset$ and ending with $\lam$ and consisting of the diagrams with their lengths at most $n$.

From the above, the number of the \lq up-down-stable ($n$-) tableaux' 
of degree $2k+1$, starting and ending with $\emptyset$ and consisting of the diagrams with their lengths at most $n$, is 
equal to the dimension $\mathbf{CP_k}$. (We note that $\dim \mathbf{CP_k}$ also depends on $n$ if $n \leqq k$.)

It seems to be interesting problem to determine the explicit action of 
$\mathbf{CP_k}$ on the \lq up-down-stable $n$-tableaux' base.

We can identify the elements in $\mathbf{GB^k_k}$ with the involutions in the symmetric group of degree $2k$. 
So we can expect a bijective correspondence between the \lq up-down-stable tableaux' of degree $2k+1$ defined above and the involutions in $\frak{S}_{2k}$. 
The author heard from I. Terada that Tom Roby privately communicated to him that the method used in his Thesis \cite{ro} also gives such a bijection.
 Namely the pictorial rules of the Robinson-Schensted correspondence initiated by Sergey Fomin \cite{fo} give a bijective correspondence between the \lq up-down-stable tableaux' of degree $k$ and the involutions in $\frak{S}_{k}$.

Finally we give the dimension of  $\mathbf{CP_k}$.

\begin{lem}\label{lem dim}
If $k \leqq n$, we have
$$
\dim \mathbf{CP_k} = (2k-1)!! + \binom{2k}{2}(2k-3)!! + \binom{2k}{4}(2k-5)!!
+ \ldots + \binom{2k}{2k-2}(1)!!+\binom{2k}{2k}, 
$$
where $(2k-1)!! = \prod_{i=0}^{k-1}(2k-1-2i)$.
\end{lem}

\section{Representation Theoretic parametrization}\label{sect reppar}

In this section we first give explicit formulas for the projection from $\Delta \bigotimes \bigwedge^k V$ to $\Delta$ and the injection from $\Delta$ to $\Delta \bigotimes \bigwedge^k V$ to $\Delta$.

Then we give another kind of parametrization of the elements in $\mathbf{CP^k_l}$ by the generalized Brauer diagrams, using the above operators.

From Theorem \ref{thm repodd} and Theorem \ref{thm repeven}, 
we have ${ \Delta \bigotimes \bigwedge^k V = \sum\limits_{i=0}^{k} [\Delta, 
(1^i)]_{Spin(N)}}$ for any $k \leqq n$.

 For $Pin(2n)$,  the representations $[\Delta, (1^i)]_{Pin(2n)}$ are determined by their characters and the above decomposition literally holds for $Pin(2n)$.
 But for $Pin(2n+1)$, if we consider the action of the center $z$, we have
the following. For $\ep_1 \in \{ \pm 1\}$, 
${ \Delta_{\ep_1} \bigotimes \bigwedge^k V = \sum\limits_{i=0}^{k} [\Delta, 
(1^i)]_{Pin(2n+1), \ep_2}}$, where $\ep_2 =(-1)^k \ep_1$.

We note that the above decomposition is multiplicity free.
So the $Pin(N)$-equivariant projection from $\Delta \bigotimes \bigwedge^k V$ to $\Delta$ and the $Pin(N)$-equivariant embedding from  $\Delta$ to $\Delta \bigotimes \bigwedge^k V$ are determined uniquely up to scalars.
 We write down those maps explicitly.

The $Pin(N)$-equivariant projection in
$
 \Hom_{Pin(N)} (\Delta \bigotimes \bigwedge\limits^k V, 
\Delta)
$ 
 $(\cong (\Delta^* \bigotimes \Delta \bigotimes (\bigwedge\limits^k V)^*))^{Pin(N)})$ 
 is given by the image of the following element under the map $\phi_k \otimes \id$.  
(See Theorem \ref{thm phipinodd} and \ref{thm phipin}.)
 
\begin{align*}
& \sum_{\substack{[1,n] \supseteq \tJ \sqcup \tI \sqcup \tW \: \text {(disjoint 
sets)} \\
 \vert \tJ \vert + \vert \tI \vert + 2 \vert \tW \vert =k}} 
 <\ur{J}, \ur{W}, \ulb{W}, \ulb{I} > \otimes {<\ur{J}, \ur{W}, \ulb{W}, \ulb{I} >}^* 
\\
&(+ \sum_{\substack{[1,n] \supseteq \tJ \sqcup \tI \sqcup \tW \: \text {(disjoint 
sets)} \\
 \vert \tJ \vert + \vert \tI \vert + 2 \vert \tW \vert = k -1}} 
 <\ur{J}, \ur{W}, \ur{0'}, \ulb{W}, \ulb{I} > \otimes {<\ur{J}, \ur{W}, \ur{0'}, \ulb{W}, \ulb{I} >}^*).\\
\end{align*}

For convenience sake, we multiply $2^{n/2}$ to the above and  denote the 
corresponding homomorphism by $\pr{k}$.
We note that for $Pin(2n+1)$ and odd $k$, $\pr{k}$ is the projection from 
$\Delta_{\pm} \textstyle{\bigotimes \bigwedge\limits^k V}$ to $\Delta_{\mp}$.

\begin{defn}\label{def pr}
For $Pin(2n+1)$, the projection $\pr{k}$ is given by 
$$
\pr{k}( [\ur{T}] \otimes <\ur{I}, \ur{W}, \ulb{W}, \ulb{J} >)
= \begin{cases}
0 \quad & \text{if $\tI \nsubseteq \tT$}, \\ 
\sign \begin{pmatrix} \ur{T} \\ \ur{I}, & \ur{K} \\ \end{pmatrix}
(-1)^{\vert \mathtt{W} - \mathtt{W} \cap \mathtt{K} \vert} 2^{(\vert \tI \vert + \vert \tJ \vert)/2} [ \ur{J}, \ur{K}] \quad & \text{if $ \tI \subseteq \tT$},
\end{cases}
$$
and 
$$
\pr{k}( [\ur{T}] \otimes <\ur{I}, \ur{W}, \ur{0'}, \ulb{W}, \ulb{J} >)
= \begin{cases}
0 \quad & \text{if $ \tI \nsubseteq \tT$ }, \\ 
\sign \begin{pmatrix} \ur{T} \\ \ur{I}, & \ur{K} \\ \end{pmatrix}
(-1)^{\vert \mathtt{K} - \mathtt{W} \cap \mathtt{K} \vert} 2^{(\vert \tI \vert + \vert \tJ \vert)/2} [ \ur{J}, \ur{K}] \quad & \text{if $ \tI \subseteq \tT$}.
\end{cases}
$$
Here we put $\tK = \tT - \tI$ and $\sign \begin{pmatrix} \ur{T} \\ \ur{I}, & \ur{K} \\ \end{pmatrix}$ denotes the sign of the permutation obtained by arranging $\ur{T}$  into $\ur{I}$, $\ur{K}$ in this order.

For $Pin(2n)$, the projection $\pr{k}$ is given as follows. 
If $k$ is even, we define 
$$
\pr{k}( [\ur{T}] \otimes <\ur{I}, \ur{W}, \ulb{W}, \ulb{J} >)
= \begin{cases}
0 \quad & \text{if $\tI \nsubseteq \tT$}, \\ 
\sign \begin{pmatrix} \ur{T} \\ \ur{I}, & \ur{K} \\ \end{pmatrix}
(-1)^{\vert \mathtt{W} - \mathtt{W} \cap \mathtt{K} \vert} 2^{(\vert \tI \vert + \vert \tJ \vert)/2} [ \ur{J}, \ur{K}] \quad & \text{if $ \tI \subseteq \tT$}.
\end{cases}
$$

If $k$ is odd, we define 
$$
\pr{k}( [\ur{T}] \otimes <\ur{I}, \ur{W}, \ulb{W}, \ulb{J} >)
= \begin{cases}
0 \quad & \text{if $\tI \nsubseteq \tT$}, \\ 
\sign \begin{pmatrix} \ur{T} \\ \ur{I}, & \ur{K} \\ \end{pmatrix}
(-1)^{\vert \mathtt{W} - \mathtt{W} \cap \mathtt{K} \vert + \vert \tK \vert +\vert \tJ \vert} 2^{(\vert \tI \vert + \vert \tJ \vert)/2} [ \ur{J}, \ur{K}] \quad & \text{if $ \tI \subseteq \tT$}.
\end{cases}
$$
\end{defn}

Similarly we have 
$\dim (\Hom_{Pin(N)} (\Delta_{(\ep_1)}, \Delta_{(\ep_2)} \bigotimes \bigwedge^k V))=1$. (If $N =2n$, ignore the subscripts $(\ep_1)$ and $(\ep_2)$ of $\Delta$. The subscripts $(\ep_1)$ and $(\ep_2)$ of $\Delta$ are only considered for the case $N = 2n+1$.   From now on, we assume this convention. Here we also assume that $\ep_2 =(-1)^k \ep_1$ is satisfied for $N=2n+1$. )

\begin{equation*}
 \Hom_{Pin(N)} (\textstyle{\Delta_{(\ep_1)}, \Delta_{(\ep_2)} \bigotimes 
\bigwedge\limits^k V) 
 \cong (\Delta_{(\ep_1)}^*  \bigotimes \Delta_{(\ep_2)} \bigotimes 
\bigwedge\limits^k V} )^{Pin(N)} \end{equation*}
and the corresponding invariant element is given by $\frak{I}_k$ (ref. \eqref{eqn ir}).
 So as in the case of $\pr{k}$,  we obtain the corresponding homomorphism 
$\inj{k}$ (after a scalar multiplication by $2^{n/2}$) as follows.
We note that for $Pin(2n+1)$ and odd $k$, $\inj{k}$ is the immersion of $\Delta_{\pm}$ into $\Delta_{\mp} \textstyle{\bigotimes \bigwedge\limits^k V}$.

\begin{defn}\label{def inj}
For $Pin(2n+1)$, the immersion $\inj{k}$ is given by 

\begin{align*}
& \inj{k}( [\ur{T}]) = \\
& \sum_{\substack{ \tT \supseteq \tI \\ \tK = \tT-\tI}}
 \sign \begin{pmatrix} \ur{T} \\ \ur{I}, & \ur{K} \\ \end{pmatrix}
(\sum_{\substack{ \tJ \subseteq ([1,n]-\tT) \\ \tW \subseteq ([1,n]-\tI -\tJ) \\
 \vert \tJ \vert + \vert \tI \vert + 2 \vert \tW \vert =k}}
(-1)^{\vert \mathtt{W} - \mathtt{W} \cap \mathtt{K} \vert} 2^{(\vert \tJ \vert + \vert \tI \vert)/2} [ \ur{J}, \ur{K}] 
  \otimes k! <\ur{J}, \ur{W}, \ulb{W}, \ulb{I} > \\
&+ \sum_{\substack{ \tJ \subseteq ([1,n]-\tT) \\ \tW \subseteq ([1,n]-\tI -\tJ) \\
 \vert \tJ \vert + \vert \tI \vert + 2 \vert \tW \vert +1 =k}}
(-1)^{\vert \mathtt{K} - \mathtt{W} \cap \mathtt{K} \vert} 2^{(\vert \tJ \vert + \vert \tI \vert)/2} [ \ur{J}, \ur{K}] 
  \otimes k! <\ur{J}, \ur{W}, \ur{0'}, \ulb{W}, \ulb{I} >).
\end{align*}

For $Pin(2n)$, the immersion $\inj{k}$ is given as follows. 
For even $k$, we define 
\begin{align*}
& \inj{k}( [\ur{T}]) = \\
& \sum_{\substack{  \tT \supseteq \tI \\ \tK = \tT-\tI}}
 \sign \begin{pmatrix} \ur{T} \\ \ur{I}, & \ur{K} \\ \end{pmatrix}
 (\sum_{\substack{ \tJ \subseteq ([1,n]-\tT) \\ \tW \subseteq ([1,n]-\tI -\tJ) \\
 \vert \tJ \vert + \vert \tI \vert + 2 \vert \tW \vert =k}}
(-1)^{\vert \mathtt{W} - \mathtt{W} \cap \mathtt{K} \vert} 2^{(\vert \tJ \vert + \vert \tI \vert)/2} [ \ur{J}, \ur{K}] 
  \otimes k! <\ur{J}, \ur{W}, \ulb{W}, \ulb{I} >) \\
\end{align*}

For odd $k$, we define 
\begin{align*}
& \inj{k}( [\ur{T}]) = \\
& \sum_{\substack{  \tT \supseteq \tI \\ \tK = \tT-\tI}}
 \sign \begin{pmatrix} \ur{T} \\ \ur{I}, & \ur{K} \\ \end{pmatrix}
 (\sum_{\substack{ \tJ \subseteq ([1,n]-\tT) \\ \tW \subseteq ([1,n]-\tI -\tJ) \\
 \vert \tJ \vert + \vert \tI \vert + 2 \vert \tW \vert =k}}
(-1)^{\vert \mathtt{W} - \mathtt{W} \cap \mathtt{K} \vert +\vert \tJ \vert +\vert \tK \vert} 2^{(\vert \tJ \vert + \vert \tI \vert)/2} [ \ur{J}, \ur{K}] 
  \otimes k! <\ur{J}, \ur{W}, \ulb{W}, \ulb{I} >). \\
\end{align*}
\end{defn}

Using the above projection $\pr{p}$ , we introduce a 
$Pin(N)$-equivariant homomorphism 
$\pr{\ur{T}}$ from $\Delta_{(\ep_1)} \bigotimes \bigotimes^k V$ to 
$\Delta_{(\ep_2)} \bigotimes \bigotimes^{k-p} V$ as follows.
(If $N=2n+1$, we assume $\ep_2 =(-1)^p \ep_1$ is satisfied.)

Let $\tT=\{t_1,t_2\ldots,t_p\}$ \, $(t_1< t_2 < \ldots < t_p )$ be a subset of 
$[1,k]$ and denote the positions in the tensor product  $\bigotimes^k V$.

First we prepare a general notation.

\begin{defn}\label{def alt}
Let $\Alt{\ur{T}}$ be the alternating operator of the tensor components of 
$\bigotimes^k V$ whose positions are in $\tT$. Namely 
\begin{align*}
&\Alt{\ur{T}} (v_1 \otimes v_2 \otimes \ldots \otimes v_k) = \\
& \dfrac{1}{p!} \sum_{ \sigma \in \frak{S}_p } \sign(\sigma)
v_1 \otimes \ldots v_{t_{\sigma^{-1}(1)}} \otimes \ldots \otimes v_{t_{\sigma^{-1}(2)}}  \otimes \ldots \otimes v_{t_{\sigma^{-1}(p)}} \otimes \ldots \otimes v_k.
\end{align*}
\end{defn}

Then $\Alt{\ur{T}}$ acts  naturally on $\Delta \bigotimes \bigotimes^k V$ and 
we consider the operator $\Alt{\ur{T}}$ itself has the alternating property with the indices in $\tT$. 
\begin{defn}\label{def prind}
$\pr{\ur{T}} : \Delta_{(\ep_1)} \bigotimes \bigotimes^k V \longrightarrow 
\Delta_{(\ep_2)} \bigotimes \bigotimes^{k-p} V$  is defined by the composition of $\Alt{\ur{T}}$ and $\pr{p}$, 
namely $\pr{\ur{T}} = \pr{p} \circ \Alt{\ur{T}}$. 
Here $\pr{p}$ acts on the alternating tensors of the position $\ur{T}$ in 
$\Delta \bigotimes \bigotimes^k V$.
\end{defn}

From the definition, $\pr{\ur{T}}$ is the element in $\Hom_{Pin(N)} 
(\Delta_{(\ep_1)} \bigotimes \bigotimes^{k} V, \Delta_{(\ep_2)} \bigotimes 
\bigotimes^{k-p} V)$ and has the alternating property with the indices in $\ur{T}$.

Next we define the element $\inj{\ur{T}}$ in $\Hom_{Pin(N)} (\Delta_{(\ep_1)} 
\bigotimes \bigotimes^{k-p} V, \Delta_{(\ep_2)} \bigotimes \bigotimes^k V)$ as 
follows.

\begin{defn}\label{def injind}
$\inj{\ur{T}} : \Delta_{(\ep_1)} \bigotimes \bigotimes^{k-p} V \longrightarrow 
\Delta_{(\ep_2)} \bigotimes\bigotimes^{k} V$  is defined by the composition of 
$\inj{p}:\Delta \longrightarrow \Delta \bigotimes \bigwedge^{p} V $ and the linear immersion of the alternating tensor to the positions indexed by $\ur{T}$.
Namely the image of this linear immersion of the exterior product $<\ur{J}, \ur{W}, \ulb{W}, \ulb{I} >$ with the $($ordered$)$ index set $\ur{T}$ 
is the tensor obtained by inserting the first component of this alternating tensor in the $t_1$-th position of $\Delta \bigotimes \bigotimes^{k-p} V$ and inserting the second component in the $t_2$-th position of the resulting tensor in $\Delta \bigotimes \bigotimes^{k-p+1} V$ and so on. 
We write this linear immersion by the exterior bracket with the subscript $\ur{T}$
e.g.,  $<\ur{J}, \ur{W}, \ulb{W}, \ulb{I} >_{\ur{T}}$.
\end{defn}

From the definition, $\inj{\ur{T}}$ has also the alternating property with the 
indices in $\ur{T}$.

So we define another parametrization of the elements in $\mathbf{CP^k_l}$ by the generalized Brauer diagrams $\mathbf{GB^k_l}$.

To a given generalized Brauer diagram, 
let us denote the indices of the isolated vertices in the upper row
by $\tT_u$ and in the lower row by $\tT_{\ell}$. 
 We correspond the following $Pin(N)$- 
equivariant homomorphism to it.

To the isolated vertices in the upper row, 
we correspond $\pr{\ur{T_u}}$ and the rest is determined by the same rule as in the usual Brauer diagrams. 
To the isolated vertices in the lower row, 
we correspond $\inj{\ur{T_{\ell}}}$ and the rest is determined by the same rule as in the usual Brauer diagrams.

We call this parametrization the representation theoretical parametrization by the generalized Brauer diagrams and to distinguish the invariant theoretical ones, we add the subscript $rt$ for the representation theoretical ones and the subscript $inv$ for the invariant theoretical ones to the generalized Brauer diagrams. 
At this moment, we have not yet known that these representation theoretical operators span the whole space of the centralizer spaces linearly. 
If $k \leqq n$ and $l \leqq n$, we show this in the section \ref{sect desinv}, namely they become linearly independent and a base of $\mathbf{CP^k_l}$.

\section{Description of the invariant theoretical equivariant 
homomorphisms}\label{sect desinv}

In this section, we give the description of the invariant theoretical equivariant homomorphisms.   First in the invariant theoretical parametrization,  we interpret the invariant element $\frak{I}_{r}$ (ref. \eqref{eqn ir}) corresponding to the isolated vertices in a generalized Brauer diagram as the equivariant 
homomorphism. 
The corresponding element in $\Delta_{(\ep_1)} \bigotimes \Delta_{(\ep_2)}^* 
\bigotimes \bigwedge^{r} V$ is given by $((\phi_r \circ \iota) \otimes \id)(\frak{I}_{r})$.
 (If $N=2n+1$, we assume $\ep_2 =(-1)^r \ep_1$ is satisfied.)

Let us recall that the index set $\mathbf{t}$ which parameterizes $\mathbf{T_{t, \{{m},{l}\}}}$ corresponds to the isolated vertices in the generalized Brauer 
diagrams.  (We note that we numbered the vertices from $1$ to $k+l$ when we defined $\mathbf{T_{t,  \{{m},{l}\}}}$ .)  We put $\tT_{u} = \mathbf{t} \cap [1, k]$ and define $\tT_{\ell}=\{s_1,s_2,\ldots,s_q\} \subseteq [1,l]$ by 
$\{k+s_1,k+s_2,\ldots,k+s_p\}= \mathbf{t} \cap [k+1, k+l]$.
We put $p = \vert \tT_{u} \vert$ and $q = \vert \tT_{\ell} \vert$
 and then $p+q = r$. 
We consider the homomorphism caused by the isolated vertices in $\mathbf{T_{t,  \{{m},{l}\}}}$ and  for convenience sake, we multiply $2^{n/2}$ to it and 
we denote the resulting homomorphism by $\psi^{\ur{T_{u}}}_{\ur{T_{\ell}}}$.
When we need not to show the index set explicitly, we just write 
$\psi^{p}_{q}$. 
As a homomorphism from $\Delta \bigotimes \bigotimes^p V$ to 
$\Delta \bigotimes \bigotimes^q V$ (here we pick up only the $\tT_{u}$-th components in the upper row $\Delta \bigotimes \bigotimes^k V$ and $\tT_{\ell}$-th components in the lower row $\Delta \bigotimes \bigotimes^l V$),  
$\psi^{\ur{T_{u}}}_{\ur{T_{\ell}}}$ is the composition of three homomorphisms, 
namely $\Alt{\ur{T_{u}}} : \Delta \bigotimes \bigotimes^p V \longrightarrow \Delta \bigotimes \bigwedge^p V$ and the homomorphism caused by the $\frak{I}_{r}$ from $\Delta \bigotimes \bigwedge^p V$ to $\Delta \bigotimes \bigwedge^q V$ and finally the linear immersion of the exterior products of degree $q$ to the prescribed positions in the tensor products. 
We note that in the space $\bigwedge^p V \bigotimes \bigwedge^q V$, 
$\bigwedge^r V$ ($r = p + q$) occurs exactly once for any $r \leqq N$.

We divide the exterior product $\bigwedge^{r} V$ occurring in $\frak{I}_{r}$ into the first $p$ tensor product and the remaining $q$ tensor product. 
We also consider  $\psi^{\ur{T_{u}}}_{\ur{T_{\ell}}}$ has the alternating property both in the superscripts and the subscripts respectively.

Direct computation tells us that the equivariant homomorphism 
$\psi^{p}_{q}$ is given as follows.

\begin{thm}\label{thm psi}
For $Pin(2n+1)$, the homomorphism $\psi^{p}_{q}$ is given as follows.

\begin{equation*}
\begin{split}
& \psi^{p}_{q}([\ur{T}] \otimes < \ur{I}, \ur{S_2}, \ulb{S_2}, \ulb{J}>) = 
 \sum_{i_2 + j_2 + s_1 + s_3 +2 s_4 = q}
(-1)^{\sum_{i=1}^4 \vert \tS_i - \tS_i \cap \tK \vert} 2^{(i_1 + j_1 + i_2 + j_2)/2} (-1)^{j_2 s_1 + s_3 j_2 + s_3 +i_1 q}  \\
&  \times \sign \begin{pmatrix} \ur{T} \\ \ur{I_1}, & \ur{I_2}, & \ur{K} 
\\ \end{pmatrix} 
\sign \begin{pmatrix} \ur{I} \\ \ur{I_1}, & \ur{S_3} \\ \end{pmatrix} 
\sign \begin{pmatrix} \ulb{J} \\ \ulb{S_1}, & \ulb{J_1} \\ \end{pmatrix}   
[ \ur{J_1}, \ur{J_2}, \ur{K}] \otimes q! <\ur{J_2}, \ur{S_3}, \ur{S_4}, \ulb{S_4}, \ulb{S_1}, \ulb{I_2} > \\
&+ \sum_{i_2 + j_2 + s_1 + s_3 +2 s_4 = q-1}
(-1)^{\vert \tK -\tK \cap(\cup_{i=1}^{4} \tS_i) \vert} 2^{(i_1 + j_1 + i_2 + j_2)/2}
(-1)^{j_2 (s_1 + s_3 )  + s_2 +i_1 q} \\
& \times \sign \begin{pmatrix} \ur{T} \\ \ur{I_1}, & \ur{I_2} & \ur{K} 
\\ \end{pmatrix} 
\sign \begin{pmatrix} \ur{I} \\ \ur{I_1}, & \ur{S_3} \\ \end{pmatrix} 
\sign \begin{pmatrix} \ulb{J} \\ \ulb{S_1}, & \ulb{J_1} \\ \end{pmatrix}  
[ \ur{J_1}, \ur{J_2}, \ur{K}] \otimes q! <\ur{J_2}, \ur{S_3}, \ur{S_4}, \ur{0'}, \ulb{S_4}, \ulb{S_1}, \ulb{I_2} >. \\
\end{split}
\end{equation*}

In the above , the sums run over all the $\tI_1$ \,$($ $\tI_1 \subseteq  \tT \cap 
\tI$ $)$ and  $\tJ_1$ \, $($ $\tJ_1 \subseteq ([1,n] - \tT) \cap \tJ$ $)$ and  
$\tI_2$\, $( \tI_2 \subseteq \tT - ( \tT \cap (\tI \cup \tS_2 \cup \tJ)))$ and 
$\tJ_2$\, $(\tJ_2 \subseteq ([1,n] - \tT) - (([1,n] - \tT) \cap (\tI \cup \tS_2 \cup \tJ)))$ and $\tS_4$\, $(\tS_4 \subseteq [1,n] - \tI -\tJ -\tS_2 - \tJ_2 - \tI_2)$. 
We also put $\tI-\tI_1 = \tS_3$ and $\tJ-\tJ_1 = \tS_1$ and 
$\tT - \tI_1 -\tI_2 = \tK$ and  we denote the sizes of the sets by their small letters. As before, we assume that the indices in the sets \ur{I}, \ulb{J} and the rest are in the increasing order.

For the element $< \ur{I}, \ur{S_2}, \ur{0'}, \ulb{S_2}, \ulb{J}>$, we have 
\begin{equation*}
\begin{split}
& \psi^{p}_{q}([\ur{T}] \otimes < \ur{I}, \ur{S_2}, \ur{0'}, \ulb{S_2}, \ulb{J}>) = \\& 
\sum_{ i_2 + j_2 + s_1 + s_3 +2 s_4 = q}
(-1)^{\vert \tK -\tK \cap(\cup_{i=1}^{4} \tS_i) \vert} 
2^{(i_1 + j_1 + i_2 + j_2)/2}
(-1)^{j_2 (s_1 + s_3 + 1)  + s_4 +i_1 q} \times \\
& \qquad \sign \begin{pmatrix} \ur{T} \\ \ur{I_1}, & \ur{I_2}, & \ur{K} 
\\ \end{pmatrix} 
\sign \begin{pmatrix} \ur{I} \\ \ur{I_1}, & \ur{S_3} \\ \end{pmatrix} 
\sign \begin{pmatrix} \ulb{J} \\ \ulb{S_1}, & \ulb{J_1} \\ \end{pmatrix} 
 [ \ur{J_1}, \ur{J_2}, \ur{K}] \otimes q! <\ur{J_2}, \ur{S_3}, \ur{S_4}, 
\ulb{S_4}, \ulb{S_1}, \ulb{I_2} > \\
\end{split}
\end{equation*}
 and the notations are the same as above.

For $Pin(2n)$, the homomorphism $\psi^{p}_{q}$ is given as follows.

For even $r = p + q$, we have 

\begin{equation*}
\begin{split}
& \psi^{p}_{q}([\ur{T}] \otimes < \ur{I}, \ur{S_2}, \ulb{S_2}, \ulb{J}>) = \\ 
& \sum_{i_2 + j_2 + s_1 + s_3 +2 s_4 = q}
(-1)^{\sum_{i=1}^4 \vert \tS_i - \tS_i \cap \tK \vert} 2^{(i_1 + j_1 + i_2 + j_2)/2} (-1)^{j_2 s_1 + s_3 j_2 + s_3 +i_1 q} \times \\
& \qquad \sign \begin{pmatrix} \ur{T} \\ \ur{I_1}, & \ur{I_2}, & \ur{K} 
\\ \end{pmatrix} 
\sign \begin{pmatrix} \ur{I} \\ \ur{I_1}, & \ur{S_3} \\ \end{pmatrix} 
\sign \begin{pmatrix} \ulb{J} \\ \ulb{S_1}, & \ulb{J_1} \\ \end{pmatrix}   
[ \ur{J_1}, \ur{J_2}, \ur{K}] \otimes q! <\ur{J_2}, \ur{S_3}, \ur{S_4}, \ulb{S_4}, \ulb{S_1}, \ulb{I_2} >, \\
\end{split}
\end{equation*}
 and the notations are the same as above.

For odd $r = p + q$, we have 

\begin{equation*}
\begin{split}
& \psi^{p}_{q}([\ur{T}] \otimes < \ur{I}, \ur{S_2}, \ulb{S_2}, \ulb{J}>) = \\ 
& \sum_{i_2 + j_2 + s_1 + s_3 +2 s_4 = q}
(-1)^{\sum_{i=1}^4 \vert \tS_i - \tS_i \cap \tK \vert +\vert \tJ_1 \vert + \vert \tJ_2 
\vert + \vert \tK \vert} 2^{(i_1 + j_1 + i_2 + j_2)/2}
(-1)^{j_2 s_1 + s_3 j_2 + s_3 +i_1 q} \times \\
& \qquad \sign \begin{pmatrix} \ur{T} \\ \ur{I_1}, & \ur{I_2}, & \ur{K} 
\\ \end{pmatrix} 
\sign \begin{pmatrix} \ur{I} \\ \ur{I_1}, & \ur{S_3} \\ \end{pmatrix} 
\sign \begin{pmatrix} \ulb{J} \\ \ulb{S_1}, & \ulb{J_1} \\ \end{pmatrix} 
[ \ur{J_1}, \ur{J_2}, \ur{K}] \otimes q! <\ur{J_2}, \ur{S_3}, \ur{S_4}, \ulb{S_4}, \ulb{S_1}, \ulb{I_2} >, \\
\end{split}
\end{equation*}
 and the notations are the same as above.

\end{thm}

For any $\sigma \in \frak{S}_k$ and $\tau \in \frak{S}_l$, we have  
$\tau \circ \psi^{\ur{T_{u}}}_{\ur{T_{\ell}}} \circ \sigma = \psi^{\sigma^{-
1}(\ur{T_{u}})}_{\tau(\ur{T_{\ell}})}$. 
So it is enough to give the explicit description 
for $\psi^{\ur{[1,p]}}_{\ur{[1,q]}}$ 
in terms of the representation theoretical operators.

\begin{thm}\label{thm psitorep}
Let $p \leqq n$ and $q \leqq n$. 

 Then we have 
\begin{equation}
\begin{split}
&\psi^{\ur{[1,p]}}_{\ur{[1,q]}}
= \sum_{i=0}^{\min(p,q)} \ep(i) 
\sum_{\substack{ \sigma \in \frak{S}_q \\ \tau \in \frak{S}_p}}
\sign(\sigma)\sign(\tau) \dfrac{\inj{\{\sigma(\ur{[i+1,q]})\}}}{(q-i)!}   
\dfrac{\begin{pmatrix} \tau(\ur{[1, i]}) \\
\sigma(\ur{[1,i]}) \\
\end{pmatrix}}{i!}
\dfrac{\pr{\{\tau(\ur{[i+1,p]})\}}}{(p-i)!}.
\end{split}
\tag{\ref{thm psitorep}.1}\label{eqn psitorep}
\end{equation}

Here $\ep(i) = (-1)^{(p-i)(q-i)}$ for $Pin(2n+1)$ and $\ep(i) = 1$ for 
$Pin(2n)$ $($ $0 \leqq i \leqq \min(p,q)$ $)$.

Also we have 
\begin{equation}
\inj{\ur{[1,q]}} \circ \pr{\ur{[1,p]}}
= \sum_{i=0}^{\min(p,q)} \ep(i) 
\sum_{\substack{ \sigma \in \frak{S}_q \\ \tau \in \frak{S}_p}}
\sign(\sigma)\sign(\tau) 
\dfrac{\psi^{{\{\tau(\ur{[i+1,p]})\}}}_{{\{\sigma(\ur{[i+1,q]})\}}}}{(q-i)! (p-i)!} \otimes 
\dfrac{\begin{pmatrix} \tau(\ur{[1,i]}) \\
\sigma(\ur{[1, i]}) \\
\end{pmatrix}}{i!}.
\tag{\ref{thm psitorep}.2}\label{eqn reptopsi}
\end{equation}

Here $\ep(i)=(-1)^{i + p q}$ for $Pin(2n+1)$ and $\ep(i)=(-1)^{i}$ for $Pin(2n)$ $($ $0 \leqq i \leqq \min(p,q)$ $)$.

We denote
$\sigma(\ur{[i+1, q]}) =\{ \sigma(\ur{i+1}), \sigma(\ur{i+2}), \ldots, \sigma(\ur{q}) \}$ and 
$\tau(\ur{[i+1, p]}) =\{ \tau(\ur{i+1}), \tau(\ur{i+2}), \ldots, \tau(\ur{p}) \}$ and 
$\begin{pmatrix} \tau(\ur{[1,i]}) \\
\sigma(\ur{[1, i]}) \\
\end{pmatrix}$ stands for the partial permutation of the tensor components 
obtained by sending the $\tau(1)$-th component in $\Delta \bigotimes \bigotimes^p V$ in the upper row to the $\sigma(1)$-th component in $\Delta \bigotimes \bigotimes^q V$ in the lower row and 
so on.

So if $k \leqq n$ and $l \leqq n$,  the generalized Brauer diagrams 
under the representation theoretical parametrization also span the whole space 
$\mathbf{CP^k_l}$ and become a linear base of $\mathbf{CP^k_l}$. 
\end{thm}

\begin{rem}
The right-hand side of  \eqref{eqn psitorep} can be considered as the composition 
of the homomorphisms in this order, but the formula \eqref{eqn reptopsi} cannot. So 
we put $\otimes$ in the right-hand side to show that each summand becomes a 
homomorphism as a whole.

\end{rem}

\begin{proof}

Both sides of the formula \eqref{eqn psitorep} are the homomorphisms from $\Delta \bigotimes \bigotimes^p V$ to $\Delta \bigotimes \bigotimes^q V$ and moreover  both of them are the compositions of the alternating operators $\Alt{\ur[1,p]}$ and the homomorphisms from $\Delta \bigotimes \bigwedge^p V$ to $\Delta \bigotimes \bigwedge^q V$ 
and the linear immersion from $\bigwedge^q V$ to $\bigotimes^q V$.

As we noted before,  we have 
$\dim \Hom_{Pin(N)}(\Delta_{(\ep_1)} \bigotimes \bigwedge^p V, \Delta_{(\ep_2)} \bigotimes \bigwedge^q V) = \min(p,q)+1$. 
Here we assume that $\ep_1 \ep_2 =(-1)^{p+q}$ for $Pin(2n+1)$. (ref. Lemma \ref{lem odd0}.)

On the way to proving the formula \eqref{eqn psitorep},
 we show that the terms occurring in the right-hand side of the formula are linearly independent.

The proof goes as follows.

For $0 \leqq i \leqq \min(p,q)$, let us put 
$$
 {}^i\! \phi^{\ur{[1,p]}}_{\ur{[1,q]}} = \sum_{\substack{ \sigma \in \frak{S}_q \\ \tau \in \frak{S}_p}}
\sign(\sigma)\sign(\tau) \dfrac{\inj{\{\sigma(\ur{[i+1, q]})\}}}{(q-i)!} 
 \dfrac{\begin{pmatrix} \tau(\ur{[1, i]}) \\
\sigma(\ur{[1, i]}) \\
\end{pmatrix}}{i!}
\dfrac{\pr{\{\tau(\ur{[i+1, p]})\}}}{(p-i)!}.
$$

We find a sequence of pairs of elements $\{x_i, y_i\}$ (here $x_i \in \Delta \bigotimes \bigotimes^p V$, $y_i \in \Delta \bigotimes \bigotimes^q V$ and $0 \leqq i \leqq \min(p,q)$) satisfying the following condition.

For $Pin(2n+1)$, the coefficients of $y_i$ in the image ${}^j\! 
\phi^{\ur{[1,p]}}_{\ur{[1,q]}} (x_i)$ is $0$  if $j \ne i$ and non-zero if $j = i$.

For $Pin(2n)$, the coefficients of $y_{2i}$ in the image ${}^j\! 
\phi^{\ur{[1,p]}}_{\ur{[1,q]}} (x_{2i})$ is $0$  if $j \ne 2i$ and non-zero if $j = 2i$ and
 the coefficients of $y_{2i+1}$ in the image ${}^j\! \phi^{\ur{[1,p]}}_{\ur{[1,q]}} (x_{2i+1})$ is $0$  if $j \ne 2i, 2i+1$ and non-zero if $j = 2i, 2i+1$.

If we find such a sequence, ${}^i\! \phi^{\ur{[1,p]}}_{\ur{[1,q]}}$'s ($0 \leqq i \leqq \min(p,q)$) are linearly independent and span the intertwining space 
$\Hom_{Pin(N)}(\Delta_{(\ep_1)} \bigotimes \bigwedge^p V, \Delta_{(\ep_2)} 
\bigotimes \bigwedge^q V)$.

We put 
\begin{equation*}
\psi^{\ur{[1,p]}}_{\ur{[1,q]}}
= \sum_{i=0}^{\min(p,q)} d_i \;{}^i\! \phi^{\ur{[1,p]}}_{\ur{[1,q]}}
\end{equation*}
and give the above sequence $\{x_i, y_i\}$ and 
apply both sides to the elements $x_i$ and compare the coefficients of $y_i$. 
Then we determine the coefficients $d_i$ inductively.

We first prove the formula \eqref{eqn psitorep} for $Pin(2n+1)$.

For abbreviation, we write  $\otimes \{\ur{i_1, i_2, \ldots, i_r} \}$
for the tensor product of the base elements $u_{i_1} \otimes u_{i_2} \otimes \ldots \otimes u_{i_r}$ and
 denote $\otimes \{\ur{k, k+1, \ldots, k+ r} \}$ by $\otimes \{\ur{[k,k+r]} \}$ and  $\otimes \{\ur{\overline{k+r}, \overline{k+r-1}, \ldots, \overline{k}}\}$ by $\otimes \{\ulb{[k+r,k]} \}$ and we use the juxtaposition of the above 
notations.

We take cases for the parities of $p$ and $q$.

Case 1 \quad $p=2r$ and $q=2s$

We note $r+s \leqq n$.

Case 1-0. \quad The determination of $d_0$.

We take $y_0 =[\ur{\emptyset}] \otimes \{\ur{[r+1, r+s]}, \ulb{[r+s,r+1]}\}$, $x_0 =[\ur{\emptyset}] \otimes \{\ur{[1, r]}, \ulb{[r, 1]}\}$.  

To calculate the image of
$
\psi^{\ur{[1,p]}}_{\ur{[1,q]}}(x_0)$, we put $\tS_2=\{[1, r]\}$ and $\tT=\tJ = \tI=\emptyset$. So $\tI_1 = \tJ_1=\tI_2=\tS_3=\tS_1=\tK=\emptyset$
and $\tJ_2 \subseteq [1,n]- \tS_2$ and $\tS_4 \subseteq [1,n]-\tJ_2- \tS_2$ and apply Theorem \ref{thm psi}.
Since we consider the coefficients of $y_0$, $\tJ_2 =\emptyset$ and 
$\tS_4= \{[r+1, r+s]\}$ and the coefficient is given by 
$(-1)^{r+s}$.
On the other hand, since there are no common numbers in the sets $\{[1, r]\}$ and $\{[r+1, r+s]\}$, only the term coming from $i=0$ contributes to the coefficient.

We calculate $\inj{\{\ur{[1,q]} \}} \circ \pr{\{\ur{[1,p]}\}}(x_0)$.
Since $\tI = \tJ = \tT=\emptyset$ and $\tW = \{\ur{[1, r]}\}$ in the definition \ref{def pr}, 
$\pr{\{\ur{[1,p]}\}}(x_0) = (-1)^r [\ur{\emptyset}]$.
We calculate the coefficient of $y_0$ 
in the image $\inj{\{\ur{[1,q]} \}}((-1)^r [\ur{\emptyset}])$.
Since $\tT = \emptyset$, $\tI = \tK = \emptyset$ and we only consider the case 
$\tJ = \emptyset$ and $\tW=\{[r+1, r+s]\}$ in the definition \ref{def inj}. 
So the coefficient is given by $(-1)^{r+s}$. 
Therefore we have $d_0 = 1 = (-1)^{p q}$ in this case.

Case 1-1. \quad The determination of $d_1$.

We take  $y_1=[\ur{r+1}] \otimes \{\ur{1}, \ur{[r+1, r+s]}, \ulb{[r+s, r+2]}\}$, $x_1 =[\ur{\emptyset}] \otimes \{\ur{[1, r]}, \ur{0'}, \ulb{[r, 2]}\}$.

In the left-hand side, since $\tI=\{1\}$, $\tS_2 =\{[2, r]\}$, $\tJ=\emptyset$, $\tT=\emptyset$, we have $\tJ_1=\emptyset$, $\tS_1=\emptyset$, $\tK=\emptyset$, $\tI_1=\emptyset$, $\tI_2=\emptyset$.
Only the case of $\tS_3=\{1\}$, $\tJ_2=\{r+1\}$ and $\tS_4 =\{[r+2, r+s]\}$ contributes to the coefficient of $[\ur{r+1}] \otimes \{\ur{1}, \ur{[r+1, r+s]}, \ulb{[r+s, r+2]}\}$. 
Namely we have 
$\psi^{\ur{[1,p]}}_{\ur{[1,q]}}([\ur{\emptyset}] \otimes \{\ur{[1, r]}, \ur{0'}, \ulb{[r, 2]}\})= \sqrt{2} (-1)^{s-1+1+1}
[\ur{r+1}] \otimes q! <\ur{r+1}, \ur{1}, \ur{[r+2, r+s]}, \ulb{[r+s, r+2]}> + \, \text{the other terms}$. 
So the coefficient is $\sqrt{2} (-1)^{s}$ in the left-hand side.

On the other hand since $\{1, [r+1, r+s], \overline{[r+s, r+2]}\} \cap \{[1, r], 0', \overline{[r, 2]}\} = \{1\}$, the possible terms to contribute to the coefficient are only $i=0$ or $i=1$. Moreover if we consider the weight of $[\ur{\emptyset}] \otimes \{\ur{[1, r]}, \ur{0'}, \ulb{[r, 2]}\}$,  this is not in the weight set of $\Delta$, so the image of $\pr{}$ must be $0$ and $i=0$ term also does not contribute to the coefficient.

In the case of $i=1$, to obtain the coefficient, $\sigma(\ur{1})$ and $\tau(\ur{1})$ must be $\ur{1}$ and in the $i=1$ term, we have only to consider the term
$
d_1 \inj{\{\ur{[2, q]} \}} \begin{pmatrix} \ur{1} \\ \ur{1} \\
\end{pmatrix} \pr{\{\ur{[2, p]}\}}.
$

Since $\pr{\{\ur{[2,p]}\}} ([\ur{\emptyset}] \otimes <\ur{[2, r]}, \ur{0'}, \ulb{[r, 2]}>) = [\ur{\emptyset}]$ and $\inj{\{\ur{[2,q]} \}}([\ur{\emptyset}]) = \sqrt{2} (-1)^{s-1}  [r+1] \otimes (q-1)! <\ur{[r+1, r+s]}, \ulb{[r+s, r+2]}>$, the coefficient of the right-hand side is $d_1 \sqrt{2} (-1)^{s-1}$.
So we have $d_1 = -1$. 

The similar argument goes on for the general case and we only give the elements $x_i$ and $y_i$ and $d_i$ for each of the cases.

Case 1-$2i$ \quad The determination of $d_{2i}$.

We take $y_{2i}=[\ur{\emptyset}] \otimes \{\ur{[1, 2i]}, \ur{[r+i+1, r+s]}, \ulb{[r+s, r+i+1]}\}$, $x_{2i}=[\ur{\emptyset}] \otimes \{\ur{[1, r+i]}, \ulb{[r+i, 2i+1]}\}$. 
In this case $d_{2i} = 1$.

Case 1-$2i+1$ \quad The determination of $d_{2i+1}$.

We take $y_{2i+1}=[\ur{r+i+1}] \otimes \{\ur{[1, 2i+1]}, \ur{[r+i+1, r+s]}, \ulb{[r+s, r+i+2]}\}$, $x_{2i+1}=[\ur{\emptyset}] \otimes \{\ur{[1, r+i]}, \ur{0'}, \ulb{[r+i, 2i+2]}\}$.
 In this case $d_{2i+1} = -1$.

Case 2. \quad $p=2r+1$ and $q=2s$.

In this case, we have $r + s < n$.

Case 2-$2i$. \quad The determination of $d_{2i}$.

We take $y_{2i}=[\ur{\emptyset}] \otimes \{\ur{[1, 2i]}, \ur{[r+i+1, r+s]}, \ulb{[r+s, r+i+1]}\}$, $x_{2i}=[\ur{\emptyset}] \otimes \{\ur{[1, r+i]}, \ur{0'}, \ulb{[r+i, 2i+1]}\}$.  
In this case $d_{2i} = 1$.

Case 2-$2i+1$. \quad The determination of $d_{2i+1}$.

We take $y_{2i+1}=[\ur{\emptyset}] \otimes \{\ur{[1, 2i+1]}, \ur{[r+i+2,r+s]}, \ur{0'}, \ulb{[r+s, r+i+2]}\}$, $x_{2i+1}=[\ur{\emptyset}] \otimes \{\ur{[1, r+i+1]} \ulb{[r+i+1, 2i+2]}\}$.   

In this case $d_{2i+1} = 1$.

Case 3. \quad $p=2r$ and $q=2s+1$.

In this case, we have $r + s < n$.

Case 3-$2i$. \quad The determination of $d_{2i}$.

We take $y_{2i} = [\ur{\emptyset}] \otimes \{\ur{[1, 2i]}, \ur{[r+i+1, r+s]}, \ulb{[r+s, r+i+1]}\}$,  $x_{2i}= [\ur{2i+1}] \otimes \{\ur{[1, r+i]}, \ur{0'}, \ulb{[r+i, 2i+2]}\}$. 

In this case $d_{2i} = 1$.

Case 3-$2i+1$. \quad The determination of $d_{2i+1}$.

We take $y_{2i+1} = [\ur{\emptyset}] \otimes \{\ur{[1, 2i+1]}, \ur{[r+i+1,r+s]}, \ulb{[r+s,r+i+1]}\}$,  $x_{2i+1} = [\ur{\emptyset}] \otimes \{\ur{[1, r+i]}, \ur{0'}, \ulb{[r+i, 2i+2]}\}$. 

In this case $d_{2i+1} = 1$.

Case 4. \quad $p=2r+1$ and $q=2s+1$.

In this case, we note that $r + s \leqq n-1$.

Case 4-$2i$. \quad The determination of $d_{2i}$.

We take $y_{2i} = [\ur{r+i+1}] \otimes \{\ur{[1, 2i]}, \ur{[r+i+1, r+s+1]}, \ulb{[r+s+1, r+i+2]}\}$, $x_{2i} =[\ur{\emptyset}] \otimes \{\ur{[1, r+i]}, \ur{0'}, \ulb{[r+i, 2i+1]}\}$.

In this case $d_{2i} = -1$.

Case 4-$2i+1$. \quad The determination of $d_{2i+1}$.

We take $y_{2i+1} = [\ur{\emptyset}] \otimes \{\ur{[1, 2i+1]}, \ur{[r+i+2, r+s+1]}, \ulb{[r+s+1, r+i+2]}\}$,  $x_{2i+1} = [\ur{\emptyset}] \otimes \{\ur{[1, r+i+1]},  \ulb{[r+i+1, 2i+2]}\}$.  

In this case $d_{2i+1} = 1$.

We verify all the cases and we have proved the formula \eqref{eqn psitorep} for $Pin(2n+1)$.

We will prove the formula \eqref{eqn psitorep} for $Pin(2n)$.

Case 1. \quad $p=2r$ and $q=2s$.

We note $r+s \leqq n$.

Case 1-0. \quad The determination of $d_0$.

We take $y_0 =[\ur{\emptyset}] \otimes \{\ur{[r+1, r+s]}, \ulb{[r+s, r+1]}\}$,  $x_0 = [\ur{\emptyset}] \otimes \{\ur{[1,r]}, \ulb{[r,1]}\}$.
In this case $d_0 = 1$.

Case 1-1. \quad The determination of $d_1$.

We take $y_1 =[\ur{1},\ur{r+1}] \otimes \{\ur{1}, \ur{[r+1, r+s]}, \ulb{[r+s, r+2]}\}$, $x_1 =[\ur{\emptyset}] \otimes \{\ur{[1,r]}, \ulb{[r,1]}\}$. 

In this case,  the coefficient of the left-hand side is $0$.

On the right-hand side, since $\{\ur{1}, \ur{[r+1, r+s]}, \ulb{[r+s, r+2]}\} \cap \{\ur{[1,r]}, \ulb{[r,1]}\} = \{\ur{1}\}$, the possible terms to contribute to the coefficient are only $i=0$ or $i=1$.

We first calculate the contribution coming from the term $i=1$.
The image of $\pr{2r-1}([\ur{\emptyset}] \otimes \{\ur{[2,r]}, \ulb{[r,1]}\}$ is $(-1)^{r-1+1} \sqrt{2}[1]$.
The image $\inj{2s-1}([1])$ is given by  $(-1)^{s-1+1+1}\sqrt{2} [\ur{r+1},\ur{1}] \otimes (2s-1)! <\ur{[r+1, r+s]}, \ulb{[r+s, r+2]}> +\text{the other terms}$. So the contribution from this term to the 
coefficient is $2 (-1)^{r+s}$.

We calculate the contribution coming from the term $i=0$.
The image of $\pr{2r}([\ur{\emptyset}] \otimes \{\ur{[1,r]}, \ulb{[r,1]}\}$ is $(-1)^{r} [\ur{\emptyset}]$.
The image $\inj{2s}([\ur{\emptyset}])$ is given by  $(-1)^{s-1} 2 [\ur{1},\ur{r+1}] \otimes (2s)! 
<\ur{1},\ur{[r+1, r+s]}, \ulb{[r+s, r+2]}> +\text{the other terms}$. So the contribution from this term to the coefficient is 
$2 (-1)^{r+s-1}$.
Hence we have $d_1=1$.

Case 1-$2i$. \quad The determination of $d_{2i}$.

We take $y_{2i} =[\ur{\emptyset}] \otimes \{\ur{[1,2i]}, \ur{[r+i+1, r+s]}, \ulb{[r+s, r+i+1]}\}$,  $x_{2i} = [\ur{\emptyset}] \otimes \{\ur{[1,r+i]}, \ulb{[r+i, 2i+1]}\}$.

In this case $d_{2i} = 1$.

Case 1-$2i+1$. \quad The determination of $d_{2i+1}$.

We take $y_{2i+1} =[\ur{2i+1}, \ur{r+i+1}] \otimes \{\ur{[1,2i+1]}, \ur{[r+i+1, r+s]}, \ulb{[r+s, r+i+2]}\}$, $x_{2i+1} = [\ur{\emptyset}] \otimes \{\ur{[1,r+i]}, \ulb{[r+i, 2i+1]}\}$. In this case $d_{2i+1}=1$.

For the remaining cases of parities of $p$ and $q$, 
the similar argument goes on and so we omit it.

Next we show the formula \eqref{eqn reptopsi}.

For $0 \leqq i \leqq \min(p,q)$, let us put 
$$
 {}^i\! \psi^{\ur{[1,p]}}_{\ur{[1,q]}} = \sum_{\substack{ \sigma \in \frak{S}_q \\ \tau \in \frak{S}_p}}
\sign(\sigma)\sign(\tau) 
\dfrac{\psi^{{\{\tau(\ur{[i+1, p]})\}}}_{{\{\sigma(\ur{[i+1, q]})\}}}}{(q-i)! (p-i)!} \otimes \dfrac{\begin{pmatrix} \tau(\ur{[1,i]}) \\
\sigma(\ur{[1,i]}) \\
\end{pmatrix}}{i!}
$$

We find a sequence of pairs of elements $\{x_i, y_i\}$ (here $x_i \in \Delta \bigotimes \bigotimes^p V$, $y_i \in \Delta \bigotimes \bigotimes^q V$ and $0 \leqq i \leqq \min(p,q)$) such that the coefficients of $y_i$ in the image ${}^j\! 
\psi^{\ur{[1,p]}}_{\ur{[1,q]}} (x_i)$ is $0$  if $j > i$ and non-zero if $j = i$.

In this case ${}^i\! \psi^{\ur{[1,p]}}_{\ur{[1,q]}}$ is equal to  $\Alt{\ur[1,q]} \circ (\phi^{\ur{[i+1,p]}}_{\ur{[i+1,q]}} \bigotimes 
\binom{\ur{[1,i]}}{\ur{[1,i]}} ) \circ \Alt{\ur[1,p]}$ up to scalar and so ${}^i\! \psi^{\ur{[1,p]}}_{\ur{[1,q]}}$'s ($0 \leqq i \leqq \min(p,q)$) are linearly independent. (See Lemma \ref{lem indep}.)

We first put 
\begin{equation*}
 \inj{\ur{[1,q]}} \circ \pr{\ur{[1,p]}}
= \sum_{u=0}^{\min(p,q)} {d_u}^{'}\; {}^u\! \psi^{\ur{[1,p]}}_{\ur{[1,q]}}
\end{equation*}
and give the above sequence $\{x_u, y_u\}$ and 
apply both sides to the elements $x_u$ and compare the coefficients of  $y_u$ as before. 
Then we can determine the coefficients ${d_u}^{'}$ inductively.

We take cases for the parities of $p$ and $q$.

Case 1. \quad $p=2r$ and $q=2s$.

We note $r+s \leqq n$ because of $p \leqq n$ and $q \leqq n$.

Case 1-0 \quad The determination of ${d_{0}}^{'}$.

We take $y_0 =[\ur{\emptyset}] \otimes \{\ur{[r+1, r+s]}, \ulb{[r+s, r+1]}\}$, $x_0 
=[\ur{\emptyset}] \otimes \{\ur{[1, r]}, \ulb{[r, 1]}\}$.  
In this case ${d_0}^{'}=1$.

Case 1-1. \quad The determination of ${d_{1}}^{'}$.

We take $y_1 =[\ur{1},\ur{r+1}] \otimes \{\ur{1},\ur{[r+1, r+s]}, \ulb{[r+s, r+2]}\}$, $x_1 = [\ur{1}] \otimes \{\ur{[1, r]}, \ur{0'}, \ulb{[r, 2]}\}$.  
The coefficient of the left-hand side is $ 2 \sqrt{2} (-1)^{s-1} $.

In the right-hand side, by considering the number of the common letters,
the $u$ term with $ u > 1$ never contributes to the coefficient.
We calculate the contribution coming from the $u=0$ term.
Then we have $\tI=\{1\}$, $\tS_2=\{[2,r]\}$,  $\tJ=\emptyset$ and $\tT=\{1\}$.
Only the case of $\tJ_2=\{r+1\}$, $\tK=\{1\}$, $\tS_4=\{[r+2, r+s]\}$, $\tS_3=\{1\}$, 
$\tI_1=\emptyset$, $\tJ_1=\emptyset$,  $\tI_2=\emptyset$ and $\tS_1 
=\emptyset$ contributes to the coefficient.
The image of the $u=0$ term is given by
$
\sqrt{2} (-1)^{s-1}[\ur{r+1},\ur{1}] \otimes q!<\ur{r+1},\ur{1}, \ur{S_4}, \ulb{S_4}> + \text{the other terms}
$
 and the coefficient is given by $\sqrt{2} (-1)^{s-1} {d_0}^{'}$. 
Next we consider the $u=1$ term.
Then only the case $\tau(\ur{1})=\sigma(\ur{1})=\ur{1}$ contributes to the coefficient and 
we calculate ${d_{1}}^{'} \psi^{\{\ur{[2,p]}\}}_{\{\ur{[2,q]}\}}$.
At this time only the case 
$\tI=\emptyset$, $\tS_2=\{[2,r]\}$, $\tT=\{1\}$, $\tJ_2=\{r+1\}$, $\tJ=\tJ_1=\emptyset$, $\tS_1=\emptyset$, $\tS_4=\{[r+2,r+s]\}$
, $\tK=\{1\}$, $\tI_1=\emptyset$, $\tI_2=\emptyset$ and 
$\tS_3=\emptyset$ contributes to the coefficient.   The image of the $u=1$ term is given by 
$
\psi^{\{\ur{[2,p]}\}}_{\{\ur{[2,q]}\}}(
[\ur{1}] \otimes < \ur{S_2}, \ur{0'}, \ulb{S_2}>) = 
(-1)^{1+1+s-1} \sqrt{2} [\ur{r+1},\ur{1}] \otimes (q-1)!<\ur{r+1},\ur{S_4},\ulb{S_4}> +\text{the other terms}
$
 and the coefficient is $\sqrt{2} (-1)^s {d_1}^{'}$ because of $[\ur{r+1},\ur{1}]=-[\ur{1},\ur{r+1}]$ in $\Delta$.

So we have ${d_1}^{'}=-1$.

We prove ${d_i}^{'} =(-1)^i$ by induction on $i$ in this case.

We assume that the above holds for $i=1,2, \ldots, 2u-1$ and prove the case
of $i=2u$.

We take $y_{2u} = [\ur{[1,2u]}] \otimes \{\ur{[1,2u]}, \ur{[r+u+1, r+s]}, \ulb{[r+s, r+u+1]}\}$, $x_{2u} = [\ur{[1,2u]}] \otimes \{\ur{[1,r+u]}, \ulb{[r+u, 2u+1]}\}$.

The coefficient of the left-hand side is $ 2^{2u} (-1)^{r+s} $.

We calculate the right-hand side.

The contribution to the coefficient under consideration from the $i=0$ term comes from $\tI=\{[1,2u]\}$, $\tS_2=\{[2u+1,r+u]\}$, $\tT=\{[1,2u]\}$, $\tJ=\tJ_1=\emptyset$, $\tS_1=\emptyset$, $\tS_4=\{[r+u+1,r+s]\}$, 
$\tK=\{[1,2u]\}$, $\tI_1=\emptyset$, $\tI_2=\emptyset$, $\tJ_2=\emptyset$ and $\tS_3=\{[1,2u]\}$.   The coefficient is given by $(-1)^{s-u+r-u+2s} {d_0}^{'}$.

We consider the contribution from the $i=1$ term.
Then we must have $\ur{1} \leqq \tau(\ur{1})=\sigma(\ur{1})=\ur{j} \leqq \ur{2u}$ and  for each $j$, the right-hand side of the $i=1$ term is reduced to 
$
\psi^{\{\ur{[1,p]-\{j\}}\}}_{\{\ur{[1,q]-\{j\}}\}}
$
 from the alternating property with the indices.

The contribution to the coefficient from this term comes from $\tI=\{[1,2u]-\{j\}\}$, $\tS_2=\{[2u+1,r+u]\}$, $\tT=\{[1,2u]\}$, $\tJ_2=\emptyset$, 
$\tJ=\tJ_1=\emptyset$, $\tS_1=\emptyset$, $\tS_4=\{[r+u+1,r+s]\}$, $\tK=\tT$, $\tI_1=\emptyset$, $\tI_2=\emptyset$ and $\tS_3=\tI$.   
Then we have 
$
\psi^{\{\ur{[1,p]-\{j\}}\}}_{\{\ur{[1,q]-\{j\}}\}}
([\ur{[1,2u]}] \otimes \{\ur{[1,r+u]-\{j\}}, \ulb{[r+u, 2u+1]}\})
 = (-1)^{s-u+r-u+2u-1}[\ur{T}] \otimes (q-1)! <\ur{S_3}, \ur{S_4}, \ulb{S_4}>_{\ur{\{[1,q]-\{j\}\}}} +\text{the other terms}.
 $

So the total contribution to the coefficient from the $i=1$ term is given by 
$\binom{2u}{1} (-1)^{s+r-1}{d_1}^{'}$.

We consider the contribution from the $i=2$ term.
Then we must have $\tau(\ur{1})=\sigma(\ur{1})$ and  $\tau(\ur{2})=\sigma(\ur{2})$ and the right-hand side of the $i=2$ term is given by 
$
\sum_{1 \leqq j_1 < j_2 \leqq 2u} \psi^{\{\ur{[1,p]-\{j_1\}-\{j_2\}}\}}_{\{\ur{[1,q]-\{j_1\}-\{j_2\}}\}} \otimes \begin{pmatrix} \ur{j_1} & \ur{j_2} \\ \ur{j_1} & \ur{j_2} \\ \end{pmatrix}.
$

The contribution to the coefficient from this term comes from $\tI=\{[1,2u]-\{j_1\}-\{j_2\}\}=\tS_3$, $\tS_2=\{[2u+1,r+u]\}$, $\tT=\{[1,2u]\}=\tK$, 
$\tJ_2=\emptyset$, $\tJ=\tJ_1=\emptyset$, $\tS_1=\emptyset$, 
$\tS_4=\{[r+u+1,r+s]\}$ and $\tI_1=\tI_2=\emptyset$.
The coefficient of the image under consideration by this map is given by 
$(-1)^{s-u+r-u+2u-2}$ and the total contribution coming from the $i=2$ term is given by $\binom{2u}{2} (-1)^{s+r}{d_2}^{'}$.

So from the same argument, we can deduce that the total contribution to the coefficient from the $i$ term is given by 
$\binom{2u}{i} (-1)^{s+r-i}{d_i}^{'}$.

From the induction hypothesis, we have ${d_{2u}}^{'}=1$ since $2^{2u} = 
\sum_{i=0}^{2u} \binom{2u}{i}$.

Next we assume that the ${d_i}^{'}=(-1)^i$ holds for $i=1,2, \ldots, 2u$ and prove the case of $i=2u+1$.

We take $y_{2u+1}=[\ur{[1,2u+1]},\ur{r+u+1}] \otimes \{\ur{[1,2u+1]}, \ur{[r+u+1, r+s]}, \ulb{[r+s,r+u+2]}\}$, $x_{2u+1}=[\ur{[1,2u+1]}] \otimes \{\ur{[1,r+u]}, \ur{0'}, \ulb{[r+u, 2u+2]}\}$.

The coefficient of the left-hand side is $ 2^{2u+1+1/2} (-1)^{s-u-1}$.

We calculate the right-hand side.

The contribution to the coefficient under consideration from the $i=0$ term comes from $\tI=\{[1,2u+1]\}=\tS_3$, $\tS_2=\{[2u+2,r+u]\}$, $\tT=\{[1,2u+1]\}=\tK$, $\tJ_2=\{r+u+1\}$, $\tJ=\tJ_1=\emptyset$, $\tS_1=\emptyset$, $\tS_4=\{[r+u+2, r+s]\}$, $\tI_1=\emptyset$ and $\tI_2=\emptyset$.   
The image of $\psi$ is 
$\sqrt{2} (-1)^{1 (2u+1+1) +s-u-1} [\ur{r+u+1}, \ur{T}] \otimes q! <\ur{r+u+1}, \ur{S_3}, \ur{S_4}, \ulb{S_4} > + \text{the other terms}$
 and the coefficient is given by $\sqrt{2}(-1)^{s-u-1} {d_0}^{'}$.

We consider the contribution from the $i=1$ term.
 As before, we must have $\ur{1} \leqq \tau(\ur{1})=\sigma(\ur{1})=\ur{j} \leqq \ur{2u}$ and  
for each $j$, 
the contribution to the coefficient from this term comes from 
$\tI=\{[1,2u+1] -\{j\}\}=\tS_3$, $\tS_2=\{[2u+2,r+u]\}$, 
$\tT=\{[1,2u+1]\}=\tK$, $\tJ_2=\{r+u+1\}$, $\tJ=\tJ_1=\emptyset$, 
$\tS_1=\emptyset$, $\tS_4=\{[r+u+2,r+s]\}$, $\tI_1=\emptyset$ and 
$\tI_2=\emptyset$.   
Then we have 
$
\psi^{\{\ur{[1,p]-\{j\}}\}}_{\{\ur{[1,q]-\{j\}}\}}
([\ur{[1,2u+1]}] \otimes \{\ur{[1,r+u]-\{\ur{j}\}}, \ur{0'}, \ulb{[r+u, 2u+2]}\})
 = \sqrt{2} (-1)^{1+ 1 (2u+1) +s-u-1}[\ur{r+u+1},\ur{T}] \otimes (q-1)! <\ur{r+u+1}, \ur{S_3}, \ur{S_4}, \ulb{S_4}>_{\ur{\{[1,q]-\{j\}\}}} +\text{the other terms}$. From the alternating property, we have $[\ur{r+u+1},\ur{T}] \otimes (q-1)! <\ur{r+u+1}, \ur{S_3}, \ur{S_4}, \ulb{S_4}>=-[\ur{T},\ur{r+u+1}] \otimes (q-1)! < \ur{S_3}, \ur{r+u+1}, \ur{S_4}, \ulb{S_4}>$ and the total contribution to the coefficient from the $i=1$ term is given 
by $\binom{2u+1}{1} \sqrt{2} (-1)^{s-u}{d_1}^{'}$.

Similarly we can calculate the contribution coming from the $i=2$ term.
As before, we must have $\tau(\ur{1})=\sigma(\ur{1})$ and  $\tau(\ur{2})=\sigma(\ur{2})$ and the 
right-hand side of the $i=2$ term is given by 
$
\sum_{1 \leqq j_1 < j_2 \leqq 2u}  \psi^{\{\ur{[1,p]-\{j_1\}-\{j_2\}}\}}_{\{\ur{[1,q]-\{j_1\}-\{j_2\}}\}} \otimes \begin{pmatrix} \ur{j_1} & \ur{j_2} \\ \ur{j_1} & \ur{j_2} \\ \end{pmatrix}.$

The contribution to the coefficient from this term comes from $\tI=\{[1,2u+1]-\{j_1\}-\{j_2\}\}=\tS_3$, $\tS_2=\{[2u+2,r+u]\}$, $\tT=\{[1,2u+1]\}=\tK$, $\tJ_2=\{r+u+1\}$, $\tJ=\tJ_1=\emptyset$, $\tS_1=\emptyset$, 
$\tS_4=\{[r+u+2,r+s]\}$ and $\tI_1=\tI_2=\emptyset$.
So from the same argument, the total contribution to the coefficient from the 
$i=2$ term is given by 
$\binom{2u+1}{2} \sqrt{2} (-1)^{s-u-1}{d_2}^{'}$.

Similarly we can calculate the contribution to the coefficient from the $i$ term and 
it is given by 
$\binom{2u+1}{i} \sqrt{2} (-1)^{s-u-1+i}{d_i}^{'}$.

So from the induction hypothesis, we obtain ${d_i}^{'}=(-1)^i$ in this case.

For the remaining cases of parities of $p$ and $q$, the proof goes almost in a similar way.  So we omit it and 
 summarizing all the cases, we have  ${d_i}^{'} =(-1)^{pq+i}$.

For $Pin(2n)$, the proof goes almost same and we omit it.

From the above two formulas, if $n \geqq k$, the generalized Brauer diagrams 
under the representation theoretical parametrization also span the full centralizer algebra $\mathbf{CP_k}$ and 
they become a linear base of $\mathbf{CP_k}$ 
since $\dim \mathbf{CP_k}$ is equal to the number of the  generalized Brauer 
diagrams.  

\end{proof}

\begin{cor}\label{cor phii}
Let $p \leqq n$ and $q \leqq n$. 
For $0 \leqq i \leqq \min(p,q)$, if we put 
$$
 {}^i\! \phi^{\ur{[1,p]}}_{\ur{[1,q]}} = \sum_{\substack{ \sigma \in \frak{S}_q \\ \tau \in \frak{S}_p}}
\sign(\sigma)\sign(\tau) \dfrac{\inj{\{\sigma(\ur{[i+1, q]})\}}}{(q-i)!} 
 \dfrac{\begin{pmatrix} \tau(\ur{[1, i]}) \\
\sigma(\ur{[1, i]}) \\
\end{pmatrix}}{i!}
\dfrac{\pr{\{\tau(\ur{[i+1, p]})\}}}{(p-i)!},
$$
 then 
${}^i\! \phi^{\ur{[1,p]}}_{\ur{[1,q]}}$'s $(0 \leqq i \leqq \min(p,q))$ are linearly independent and span the intertwining space 
$\Hom_{Pin(N)}(\Delta_{(\ep_1)} \bigotimes \bigwedge^p V, \Delta_{(\ep_2)} 
\bigotimes \bigwedge^q V)$.
\end{cor}

\begin{exmp}
We show a few example of the relations of two parametrization 
in figure \ref{spincentfig4} for $n \geqq 2$.

\begin{figure}[htbp]
\begin{minipage}{\linewidth}
\begin{center}
\input{spincentfig4.tex}
\end{center}
\end{minipage}
\caption{The relations of two parametrization for $Pin(2n+1)$}
\label{spincentfig4}
\end{figure}
\end{exmp}

In figure \ref{spincentfig4}, if we change simultaneously the signatures of the diagrams which have exactly one vertical edge (the edge from the upper row to the lower row), we have the result for $Pin(2n)$. 

If $p > n$,  we have $\Delta_{(\ep)} \bigotimes \bigwedge^p V \cong \Delta_{(\ep)} \bigotimes \bigwedge^{N-p} V \bigotimes \det \cong 
\Delta_{(-\ep)} \bigotimes \bigwedge^{N-p} V$.

For $Pin(2n+1)$, the isomorphism $A : \Delta_{\ep} \bigotimes \det \cong 
\Delta_{-\ep}$ is given by the map $A [\tI] = [\tI]$ for any $\tI$. Here we obey the convention in Remark \ref{rem dphi}.
Then we have the $Pin(2n+1)$ isomorphism 
$A \otimes r_p : \Delta_{(\ep)} \bigotimes \bigwedge^p V \longrightarrow 
\Delta_{(\ep)} \bigotimes \bigwedge^{2n+1-p} V$. 
Here $r_p$ is defined in Lemma \ref{lem rl}.

 For $Pin(2n)$, let $A$ be the associator in Lemma \ref{lem assoc}, the isomorphism from $\Delta \bigotimes \det$ to $\Delta$.
 Then we have the isomorphism $A \otimes r_p : \Delta \bigotimes \bigwedge^p V \longrightarrow \Delta \bigotimes \bigwedge^{2n-p} V$.

As before, we define a homomorphism $r^{\ur{[1,p]}}_{\ur{[1,2n-p]}}$ by the 
composition of the alternating operator $\Alt{\ur{[1,p]}}$ and $r_p$ and the 
immersion of the alternating tensor to the positions $[1,2n-p]$.

\begin{thm}\label{thm psitorep2}
Let $p > n$ and $p + q \leqq N$. 
As a homomorphism from $\Delta_{(\ep_1)} \bigotimes \bigotimes^p V$ to 
$\Delta_{(\ep_2)} \bigotimes \bigotimes^q V$ $($ we assume $(-1)^{p+q} = \ep_1 \ep_2$ for $Pin(2n+1)$. $)$, we have 

\begin{equation}
\psi^{\ur{[1,p]}}_{\ur{[1,q]}}
= \ep \sum_{\substack{ \sigma \in \frak{S}_{N-p} \\ \tau \in \frak{S}_q}}
\sign(\sigma)\sign(\tau)  \dfrac{\begin{pmatrix} \sigma(\ur{[1, q]}) \\
\tau(\ur{[1,q]}) \\
\end{pmatrix}}{q!}
\dfrac{\pr{\{\sigma(\ur{[q+1,N-p]})\}}}{(N-p-q)!}
(A \otimes r^{\ur{[1,p]}}_{\ur{[1,N-p]}}).
\tag{\ref{thm psitorep2}.1}\label{eqn psitorep2}
\end{equation}

Here $\ep = 1$ for $Pin(2n+1)$ and $\ep = (-1)^{q+n}$ for $Pin(2n)$.

Let $q > n$ and $p + q \leqq N$.  As homomorphisms from $\Delta_{(\ep_1)} \bigotimes \bigotimes^p V$ to $\Delta_{(\ep_2)} \bigotimes \bigotimes^q V$ $($ we assume $(-1)^{p+q} = \ep_1 \ep_2$ for $Pin(2n+1)$. $)$ is given by

\begin{equation}
\psi^{\ur{[1,p]}}_{\ur{[1,q]}}
= \ep \dfrac{q!}{(N-q)!}(A \otimes  r^{\ur{[1,N-q]}}_{\ur{[1,q]}}) \sum_{\substack{ \sigma \in \frak{S}_{p} \\ \tau \in \frak{S}_{N-q}}}
\sign(\sigma)\sign(\tau)  
\dfrac{\inj{\{\tau(\ur{[p+1,N-q]})\}}}{(N-p-q)!}
\dfrac{\begin{pmatrix} \sigma(\ur{[1, p]}) \\
\tau(\ur{[1, p]}) \\
\end{pmatrix}}{p!}.
\tag{\ref{thm psitorep2}.2}\label{eqn psitorep3}
\end{equation}

Here $\ep = 1$ for $Pin(2n+1)$ and $\ep = (-1)^{p+n}$ for $Pin(2n)$.
\end{thm}

\begin{proof}
We prove the formula \eqref{eqn psitorep2} for $Pin(2n+1)$.
 Let $p > n$ and $q \leqq n$. 
Since $q \leqq 2n+1-p$ and Corollary \ref{cor phii}, 
 we can put 
\begin{equation*}
\psi^{\ur{[1,p]}}_{\ur{[1,q]}}
= \sum_{i=0}^{q} d_i \;{}^i\! \phi^{\ur{[1,2n+1-p]}}_{\ur{[1,q]}}
(A \otimes r^{\ur{[1,p]}}_{\ur{[1,2n-p]}}).
\end{equation*}

As before, we find a sequence of pairs $\{x_i, y_i\}$ and apply both sides of 
the above to $x_i$'s and compare the coefficients of $y_i$'s and 
 determine the coefficients $d_i$ from $i = q$ to $i=0$ inductively.

We take cases for the parities of $p$ and $q$.

Case 1 \quad $ p \equiv 1 \mod{2} $ and $ q \equiv 0 \mod{2}$.

We put $2n+1-p = 2r$ and $q = 2s$. We have $ r \geqq s$.

We determine the coefficients $d_q$ first.

Case 1- $2s$. \quad $d_{2s}=1$.

We take $y_{2s} = [\ur{\emptyset}] \otimes \{\ur{[1,2s]} \}$, $x_{2s} = [\ur{\emptyset}] \otimes 
\{\ur{[1,2s]},\ur{[r+s+1, n]}, \ur{0'}, \ulb{[n,r+s+1]}\}$.
If $i < 2s$, the term $i$ never contributes to the coefficient since the image of $\pr{}$ is $0$ because of the weight argument. And the direct calculation tells us that $d_{2s}=1$.

Case 1- $2s-1$. \quad $d_{2s-1}=0$.

We take $y_{2s-1} =[\ur{\emptyset}] \otimes \{\ur{[1,2s-1]}, \ur{0'} \}$, $x_{2s-1} =[\ur{\emptyset}] \otimes \{\ur{[1,2s-1]},\ur{[r+s, n]}, \ulb{[n,r+s]}\}$.
Only the terms $2s$ and $2s-1$  contribute to the coefficient because of the weight argument. And the direct calculation tells us that $d_{2s-1}=0$

We assume that $d_{2s}=1$, $d_{2s-1}= \ldots = d_{2s-2i+1}=0$ and we show that 
$d_{2s-2i}=0$.

Case 1- $2s-2i$. \quad $d_{2s-2i}=0$.

We take $y_{2s-2i} =[\ur{\emptyset}] \otimes \{\ur{[1, 2s-i]}, \ulb{[2s-i, 2s-2i+1]} \}$, 
$x_{2s-2i} =[\ur{\emptyset}] \otimes \{\ur{[1,2s-2i]},\ur{[r+s+1-i, n]}, \ur{0'}, \ulb{[n, r+s+1-i]}\}$. 
From the induction hypothesis and the weight argument, only the terms $2s$ and 
$2s-2i$  contribute to the coefficient and the direct calculation tells us that $d_{2s-2i}=0$.

We assume that $d_{2s}=1$, $d_{2s-1}= \ldots = d_{2s-2i}=0$ and we show that 
$d_{2s-2i-1}=0$.

Case 1- $2s-2i-1$ \quad $d_{2s-2i-1}=0$ ($i \geqq 1$).

We take $y_{2s-2i-1} = [\ur{\emptyset}] \otimes \{\ur{[1, 2s-i-1]},\ur{0'}, \ulb{[2s-i-1, 2s-2i]} \}$, $x_{2s-2i-1} =[\ur{\emptyset}] \otimes \{\ur{[1,2s-2i-1]},\ur{[r+s-i, n]}, \ulb{[n, r+s-i]}\}$.

For the remaining cases of parities of $p$ and $q$, the proof goes almost in a similar way. So we omit it. Summarizing all the cases, we have $d_q = 1$ and $d_i =0$ if $i < q$ in the formula \eqref{thm psitorep2}.

Next we prove the formula \eqref{thm psitorep2} for $Pin(2n+1)$.

Let $p \leqq n$ and $q > n$. 

We can put 
\begin{equation*}
\psi^{\ur{[1,p]}}_{\ur{[1,q]}}
= \sum_{i=0}^{p} d_i
\dfrac{q!}{(2n+1-q)!}(A \otimes r^{\ur{[1,2n+1-q]}}_{\ur{[1,q]}}) \;
{}^i\! \phi^{\ur{[1,p]}}_{\ur{[1,2n+1-q]}}
\end{equation*}

We determine the coefficients $d_i$ just as in the above.

We take cases for the parities of $p$ and $q$.

Case 1 \quad $ p \equiv 0 \mod{2} $ and $ q \equiv 1 \mod{2}$.

We put $p = 2r$ and $2n+1-q = 2s$. We have $ s \geqq r$.

We determine the coefficients $d_p$ first.

Case 1- $2r$. \quad $d_{2r}=1$.

We take $y_{2r} =[\ur{\emptyset}] \otimes \{\ur{[1,2r]},\ur{[r+s+1, n]}, \ur{0'}, \ulb{[n, r+s+1]}\}$, $x_{2r} =[\ur{\emptyset}] \otimes \{\ur{[1,2r]} \}$. We have  $d_{2r}=1$.

Case 1- $2r-1$. \quad $d_{2r-1}=0$.

We take $y_{2r-1} =[\ur{\emptyset}] \otimes \{\ur{[1,2r-1]},\ur{[r+s, n]}, \ulb{[n, r+s]}\}$, 
$x_{2r-1} =[\ur{\emptyset}] \otimes \{\ur{[1,2r-1]}, \ur{0'} \}$.  We have  $d_{2r-1}=0$.

We assume that $d_{2r}=1$, $d_{2r-1}= \ldots = d_{2r-2i+1}=0$ and we show that 
$d_{2r-2i}=0$.

Case 1- $2r-2i$. \quad $d_{2r-2i}=0$. ($i \geqq 1$).

We take $y_{2r-2i} =[\ur{\emptyset}] \otimes \{\ur{[1,2r-2i]},\ur{[r+s-i+1, n]}, \ur{0'}, \ulb{[n, r+s-i+1]}\}$, $x_{2r-2i} =[\ur{\emptyset}] \otimes \{\ur{[1, 2r-i]}, \ulb{[2r-i,2r-2i+1]} \}$.  We have  $d_{2r-2i}=0$.

We assume that $d_{2r}=1$, $d_{2r-1}= \ldots = d_{2r-2i}=0$ and we show that 
$d_{2r-2i-1}=0$.

Case 1- $2r-2i-1$. \quad $d_{2r-2i-1}=0$. ($i \geqq 1$).

We take $y_{2r-2i-1} =[\ur{\emptyset}] \otimes \{\ur{[1,2r-2i-1]},\ur{[r+s-i, n]}, \ulb{[n, r+s-i]}\}$, $x_{2r-2i-1} =[\ur{\emptyset}] \otimes \{\ur{[1, 2r-i-1]}, \ur{0'}, \ulb{[2r-i-1, 2r-2i]} \}$.  We have  $d_{2r-2i-1}=0$.

For the remaining cases, proof goes almost in a similar way. So we omit it and  summarizing all the cases, we have $d_p = 1$ and $d_i =0$ if 
$i < p$ in the formula \eqref{thm psitorep2}.

For $Pin(2n)$, the proof is almost similar to the above and we omit it.

\end{proof}

\section{Relations between $Pin(N)$- equivariant homomorphisms}\label{sect rel}

In this section we give the relations  between the compositions of $\pr{}$ and 
$\inj{}$ and the contractions and the immersion of the invariant forms, from which we can deduce the multiplication rules of the generalized Brauer diagrams.

We recall some notations. Let  $\cont{\{\ur{i,j}\}}$  be the contraction operator of the $i$-th and $j$-th tensor components $($with respect to $S$ $)$ and let ${\id_{V}}_{\{\ur{i,j}\}}$ be the linear immersion of the invariant element 
$\id_{V} = \sum_{i=1}^{n} ({\ua{i}} \otimes {\ub{i}} + {\ub{i}} \otimes {\ua{i}}) (+ {\uo} \otimes {\uo})$ to the $i$-th and the $j$-th position.

As in Corollary \ref{cor phii}, we define the equivariant homomorphism from $\Delta \bigotimes \bigotimes^p$ $($here we consider the tensor $\bigotimes^{p} V$ sits in the positions $\{[r+1,r+p]\}$ $)$ to $\Delta \bigotimes 
\bigotimes^q$ $($here we consider the tensor $\bigotimes^{q} V$ sits in the positions $\{[s+1,s+q]\}$ $)$ by 
$$
 {}^i\! \phi^{\ur{[r+1,p+r]}}_{\ur{[s+1,s+q]}} = \sum_{\substack{ \sigma \in 
\frak{S}_q[s+1,s+q] \\ \tau \in \frak{S}_p[r+1,r+p]}}
\sign(\sigma)\sign(\tau) \dfrac{\inj{\{\sigma(\ur{[s+i+1, s+q]})\}}}{(q-i)!} 
 \dfrac{\begin{pmatrix} \tau(\ur{[r+1, r+i]}) \\
\sigma(\ur{[s+1, s+i]}) \\
\end{pmatrix}}{i!}
\dfrac{\pr{\{\tau(\ur{[r+i+1, r+p]})\}}}{(p-i)!}.
$$
Here $\frak{S}_q[s+1,s+q]$ and $\frak{S}_p[r+1,r+p]$ denote the symmetric groups of degree $q$ and $p$ acting on the sets $\{s+1, s+2, \ldots, s+q\}$ and 
$\{r+1, r+2, \ldots, r+p\}$ respectively.

Then from Corollary \ref{cor phii}, ${}^i\! \phi^{\ur{[r+1,p+r]}}_{\ur{[s+1,s+q]}}$'s ($ 0 \leqq i \leqq \min(p, q)$) become a base of $\Hom_{Pin(N)}(\Delta \bigotimes \bigwedge^p V, \Delta \bigotimes \bigwedge^q V)$.

\begin{thm}\label{thm reprel}

The following formulas hold for $\pr{}$, $\inj{}$,  $\cont{\{i,j\}}$ and 
${\id_{V}}_{\{i,j\}}$.

\begin{enumerate}
\item
If $p \leqq n$,  as the homomorphisms from $\Delta \bigotimes \bigotimes^p V$ to $\Delta$,  $($here we consider the tensor $\bigotimes^{p} V$ sits in the positions $\{[q+1,p+q]\}$ $)$, 
we have 
\begin{equation}
\pr{\{\ur{[1,q+p]}\}} \circ \inj{\{\ur{[1,q]}\}} = (N-p)_{q}\; \ep
\; \pr{\{\ur{[q+1, q+p]}\}}.
\tag{\ref{thm reprel}.1}\label{eqn prinj1}
\end{equation}

Here $\ep = 1$ for $Pin(2n+1)$ and $\ep = (-1)^{(p+q)q}$ for $Pin(2n)$.

$(N-p)_{q}$ denotes the lower factorial, that is, for any $x$ and non-negative integer $i$,  we define 
$(x)_{i} = x (x-1) (x-2) \cdots (x-(i-1))$.

The above also holds for $p=0$ and in this case we regard $\pr{}$ in the right-hand side as the identity map of $\Delta$.

\item
If $p \leqq n$,  as the homomorphisms from $\Delta$ to $\Delta \bigotimes 
\bigotimes^p V$\, $($here we consider the tensor $\bigotimes^{p} V$ sits in the positions $\{[q+1,p+q]\}$ $)$, 
we have 
\begin{equation}
\pr{\{\ur{[1,q]}\}} \circ \inj{\{\ur{[1, q+p]}\}} = (N-p)_{q}\; \ep \;
\inj{\{\ur{[q+1,q+p]}\}}.
\tag{\ref{thm reprel}.2}\label{eqn injpr2}
\end{equation}

Here $\ep = 1$ for $Pin(2n+1)$ and $\ep = (-1)^{(p+q)q}$ for $Pin(2n)$.

\item
If $p \leqq n$ and $q \leqq n$,  as the homomorphisms from $\Delta$ to $\Delta 
\bigotimes \bigotimes^{p+q} V$, 
we have 
\begin{equation}
\begin{split}
& \inj{\{\ur{[1,q]}\}} \circ \inj{\{\ur{[q+1,q+p]}\}} = \\
&\sum_{i=0}^{\min(p,q)} \ep(i) \; 
\sum_{\substack{\sigma \in \frak{S}_q \\ \tau \in \frak{S}_p[q+1,q+p]}}
\sign(\sigma) \sign(\tau) \dfrac{ \prod_{u=1}^{i} 
{\id_{V}}_{\{\sigma(\ur{u}),\tau(\ur{q+u})\}}}{i!} 
\dfrac{\inj{\{\sigma(\ur{[i+1, q]}),\tau(\ur{[q+i+1,q+p]})\}}}{(q-i)! (p-i)!}.
\end{split}
\tag{\ref{thm reprel}.3}\label{eqn injinj}
\end{equation}

Here $\ep(i) = (-1)^{q i+\binom{i+1}{2}}$ for $Pin(2n+1)$ and $\ep(i) = (-1)^{q i+\binom{i+1}{2}+pq}$ for $Pin(2n)$.

\item
If $p \leqq n$ and $q \leqq n$,  as the homomorphisms from $\Delta \bigotimes 
\bigotimes^{p+q} V$ to  $\Delta$, 
we have 
\begin{equation}
\begin{split}
& \pr{\{\ur{[1,q]}\}} \circ \pr{\{\ur{[q+1,q+p]}\}} = \\
&\sum_{i=0}^{\min(p,q)} \ep(i)  
\sum_{\substack{\sigma \in \frak{S}_q \\ \tau \in \frak{S}_p[q+1,q+p]}}
\sign(\sigma) \sign(\tau) \dfrac{\pr{\{\sigma(\ur{[i+1,q]}),\tau(\ur{[q+i+1,q+p]})\}}}{(q-
i)! (p-i)!}
\dfrac{ \prod_{u=1}^{i} {\cont{\{\sigma(\ur{u}),\tau(\ur{q+u})\}}}}{i!} .
\end{split}
\tag{\ref{thm reprel}.4}\label{eqn prpr}
\end{equation}

Here $\ep(i)= (-1)^{q p + q i+\binom{i}{2}}$ for $Pin(2n+1)$ and $\ep(i) = (-1)^{q i+\binom{i}{2}}$ for $Pin(2n)$.

\item
If $p \geqq t$ and $p-t \leqq n$ and $q \leqq n$,  as the homomorphisms from $\Delta \bigotimes \bigotimes^{p-t} V$ to  $\Delta \bigotimes \bigotimes^{q} V$, $($ here we consider the tensor $\bigotimes^{p-t} V$ sits in the positions 
$\{[q+t+1,q+p]\}$.$)$ 
we have 
\begin{equation}
 \pr{\{\ur{[q+1,q+p]}\}} \circ \inj{\{\ur{[1, q+t]}\}}  = 
\sum_{i=0}^{\min(p-t,q)} \ep(i) (\sum_{u=0}^{i} \binom{i}{u} (N-p-q+t+i-u)_{t})  \;
{}^i\! \phi^{\ur{[q+t+1,q + p]}}_{\ur{[1,q]}}
\tag{\ref{thm reprel}.5} \label{eqn prinj4}
\end{equation}

Here $\ep(i)= (-1)^{(q-i)(p-i)+it}$ for $Pin(2n+1)$ and $\ep(i) = (-1)^{p q + t (p+q)}$ for $Pin(2n)$ and  $\binom{i}{u}$ in the parentheses denotes the ordinary binomial coefficient. If $t=0$, $(N-p-q+0+i-u)_{0}=1$ and the sum in the parentheses is equal to $2^i$.

\item
If $p \leqq n$ and $q \leqq n$,  as the homomorphisms from $\Delta \bigotimes 
\bigotimes^{q} V$ to  $\Delta \bigotimes \bigotimes^{p} V$, $($ here we consider the tensor $\bigotimes^{q} V$ sits in the positions $\{[p+q+1,p+2q]\}$.$)$ 
we have 
\begin{equation}
\prod_{i=1}^{q} \cont{\{\ur{i,p+q+i}\}}  \inj{\{\ur{[1, q+p]} \}} = 
\sum_{i=0}^{\min(p,q)} \ep(i) \;
{}^i\! \phi^{\ur{[p+q+1,p + 2q]}}_{\ur{[q+1,q+p]}}.
\tag{\ref{thm reprel}.6}\label{eqn continj}
\end{equation}

Here $\ep(i) = (-1)^{p q +\binom{q}{2} +i (p+q-1)}$ for $Pin(2n+1)$ and $\ep(i) = (-1)^{\binom{q}{2}}$ for $Pin(2n)$.

\item
If $p \leqq n$ and $q \leqq n$,  as the homomorphisms from $\Delta \bigotimes 
\bigotimes^{q} V$ to  $\Delta \bigotimes \bigotimes^{p} V$, $($here we consider the tensor $\bigotimes^{q} V$ sits in the positions $[1,q]$ and $\bigotimes^{p} V$ sits in the positions $\{[p+q+1,2p+q]\}$.$)$ 
we have 
\begin{equation}
\pr{\{[1,q+p] \}} \prod_{i=1}^{p} {\id_V}_{\{\ur{q+i,p+q+i}\}} = 
\sum_{i=0}^{\min(p,q)} \ep(i) \;
{}^i\! \phi^{\ur{[1, q]}}_{\ur{[p+q+1,q+2p]}}.
\tag{\ref{thm reprel}.7}\label{eqn prid}
\end{equation}

Here $\ep(i) = (-1)^{\binom{p}{2} +i (p+q+1)}$ for $Pin(2n+1)$ and $\ep(i) =(-
1)^{\binom{p}{2} + p q}$ for $Pin(2n)$.

\end{enumerate}
\end{thm}

\begin{rem}
If we change $N$ in the above formulas into an indeterminate $X$ simultaneously, we can define the \lq generic' centralizer algebra of $\mathbf{CP_k}$ just as in the ordinary Brauer algebras.
\end{rem}

\begin{proof}
We prove the formula \eqref{eqn prinj1}.  If we note that  
$\Alt{\ur{[1,p+q]}}( \Alt{\ur{[1,q]}} \times \Alt{\ur{[q+1, p+q]}}) = 
\Alt{\ur{[1,p+q]}}(\Alt{\ur{[1,q]}}) = \Alt{\ur{[1,p+q]}}$, 
  both sides factor through the  homomorphisms from $ \Delta \bigotimes 
\bigwedge^p V $ to $\Delta$. Because of $\dim \; \Hom_{Pin(N)} (\Delta \bigotimes \bigwedge^p V, \Delta) =1$, they only differ in some scalar if they are non-zero. 
We compare the coefficients of $[\ur{\emptyset}]$ of the images of 
 $[\ur{[1,p]}] \otimes \{\ur{[1,p]}\}$
 under both sides of the formula \eqref{eqn prinj1}.
Direct calculation tells us that the scalar is given by $(2n+1-p)_{q}$ for $Pin(2n+1)$ and the scalar is given by $(-1)^{(p+q)q} (2n-p)_{q}$ for $Pin(2n)$.

We prove the formula \eqref{eqn injpr2}.

Both sides of the above can be considered as the  homomorphisms from $\Delta$ to $ \Delta \bigotimes \bigwedge^p V $. As before, they only differ in some scalar.
We compare the coefficients of $[\ur{[1,p]}] \otimes \{\ur{[1,p]}\}$ of the images of $[\ur{\emptyset}]$ under both sides of the formula \eqref{eqn injpr2}.

Direct calculation tells us that the scalar is given by $(2n+1-p)_{q}$ for $Pin(2n+1)$ and $(-1)^{(p+q)q}(2n+1-p)_{q}$ for $Pin(2n)$.

Next we prove the formula \eqref{eqn injinj}.

From Theorem \ref{thm repodd} and Theorem \ref{thm repeven}, 
$\Delta$ occurs in the representation $\Delta \bigotimes \mu_{O(N)}$ if and only if $\mu_{O(N)} = (1^i)$ and then the multiplicity is one. Also the character theory tells us that in the representation $\bigwedge^q V \bigotimes \bigwedge^p V$, 
$\bigwedge^{p+q-2i} V$ ($0 \leqq i \leqq \min(p,q)$) occurs exactly once and the other exterior products never appears.

So the left-hand side in \eqref{eqn injinj} 
can be considered as the homomorphism in $\Hom_{Pin(N)}( \Delta, \Delta \bigotimes \bigwedge^q V \bigotimes \bigwedge^p V) \cong \bigoplus_{i=0}^{\min(p,q)} 
\Hom_{Pin(2n)}( \Delta, \Delta \bigotimes \bigwedge^{p+q-2i} V)$.

A $Pin(2n)$ equivariant embedding from $\bigwedge^{p+q-2i} V$ to $ \bigwedge^q V 
\bigotimes \bigwedge^p V$ is given as follows.
For any subset $\mathtt{U} \subseteq \{[1,n],(0'), 
\overline{[n,1]}\}$ of size $p+q-2i$, we define 
\begin{equation*}
 <\ur{U}> \longrightarrow 
 \sum_{ \substack{ 
\sigma \in \frak{S}_q \\
\tau \in \frak{S}_p[q] }} 
\sign(\sigma) \sign(\tau) \prod_{u=1}^{i} {\id_{V}}_{\{\sigma(\ur{u}),\tau(\ur{q+u})\}} 
\times   <\ur{U}>_{\{\sigma(\ur{[i+1, q]}), \tau(\ur{[q+i+1, q+p]})\}}.
\end{equation*}
Here the indices in $\ur{U}$ are in the increasing order.

So in the  formula \eqref{eqn injinj}, the $i$ term  in the right-hand side 
 is the homomorphism in the direct summand $\Hom_{Pin(N)}( \Delta, \Delta \bigotimes \bigwedge^{p+q-2i} V)$.

So we can put

\begin{equation*}
\begin{split}
& \inj{\{\ur{[1,q]}\}} \circ \inj{\{\ur{[q+1,q+p]}\}} = \\
&\sum_{i=0}^{\min(p,q)} d_i
\sum_{\substack{\sigma \in \frak{S}_q \\ \tau \in \frak{S}_p[q+1,q+p]}}
\sign(\sigma) \sign(\tau) \dfrac{ \prod_{u=1}^{i} 
{\id_{V}}_{\{\sigma(\ur{u}),\tau(\ur{q+u})\}}}{i!} 
\dfrac{\inj{\{\sigma(\ur{[i+1, q]}),\tau(\ur{[q+i+1,q+p]})\}}}{(q-i)! (p-i)!}.
\end{split}
\end{equation*}

As before, we find a sequence of pairs $\{x_i, y_i\}$ and apply both sides of the above to $x_i$'s and compare the coefficients of $y_i$'s and 
 determine the coefficients $d_i$ inductively.

We prove the formula \eqref{eqn injinj} for $Pin(2n+1)$.

We take the cases for the parities of $p$ and $q$.

Case 1. \quad $p=2r$ and $q=2s$.

We note $r+s \leqq n$ because of $p \leqq n$ and $q \leqq n$.

Case 1-0. \quad The determination of $d_{0}$.

We take $y_0 = [\ur{\emptyset}] \otimes \{\ur{[1,s]}, \ulb{[s, 1]}, \ur{[s+1, r+s]}, \ulb{[r+s, s+1]}\}$,  $x_0 = [\ur{\emptyset}]$.  
In this case  ${d_0}=1$.

Case 1-1. \quad The determination of ${d_{1}}$.

We take $y_1 =[\ur{s+1}] \otimes \{\ulb{1},\ur{[2,s]}, \ur{0'}, \ulb{[s, 2]}, \ur{1}, \ur{[s+1, r+s]}, \ulb{[r+s, s+2]}\}$,  $x_1 =[\ur{\emptyset}]$.  
In this case ${d_1}=-1$.

We prove ${d_i} =(-1)^{\binom{i+1}{2}}$ by induction on $i$ in this case.

Case 1-$2u$. \quad The determination of ${d_{2u}}$.

We assume that the above holds for $i=1,2, \ldots, 2u-1$ and prove the case
of $i=2u$.
We take
\begin{align*}
&y_{2u} = [\ur{\emptyset}] \otimes \\
&\{\ulb{[2u,1]}, \ur{[2u+1, s+u]}, \ulb{[s+u, 2u+1]}, \ur{[1,2u]}, \ur{[s+u+1,r+s]}, \ulb{[r+s, s+u+1]}\}\\
\end{align*}
 and $x_{2u} = [\ur{\emptyset}]$.
Then in the right-hand side, only the terms from $0$ to $2u$ contribute to the 
coefficient and we have $d_{2u} =(-1)^{\binom{2u+1}{2}}$.

Case 1-$2u+1$. \quad The determination of ${d_{2u+1}}$.

We take
\begin{align*}
&y_{2u+1} = [\ur{s+u+1}] \otimes \\
&\{\ulb{[1, 2u+1]}, \ur{[2u+2, s+u]}, \ur{0'}, \ulb{[s+u, 2u+2]}, \ur{[1, 2u+1]}, \ur{[s+u+1, r+s]}, \ulb{[r+s, s+u+2]}\}\\
\end{align*}
and $x_{2u+1}=[\ur{\emptyset}]$.

For the remaining cases for $Pin(2n+1)$ and for $Pin(2n+1)$, the proof goes in a similar way and  so we omit it.

Next we prove the formula \eqref{eqn prpr}.
As in the proof of the formula \eqref{eqn injinj}, 
the left-hand side can be considered as the homomorphism in $\Hom_{Pin(N)}( \Delta \bigotimes \bigwedge^q V \bigotimes \bigwedge^p V, \Delta) \cong 
\bigoplus_{i=0}^{\min(p,q)} \Hom_{Pin(N)}(\Delta \bigotimes \bigwedge^{p+q-2i} V, \Delta)$. 
A $Pin(N)$-equivariant embedding from $ \bigwedge^q V \bigotimes \bigwedge^p V$ to $\bigwedge^{p+q-2i} V$  is given as follows.
For any subsets $\mathtt{Q}, \mathtt{P} \subseteq \{[1,n],(0'), \overline{[n,1]} \}$ of size $q$, $p$ respectively, we define

\begin{align*}
& <\ur{Q}> \otimes <\ur{P}> \longrightarrow 
 \Alt{\{\ur{[i+1, q]}, \ur{[q+i+1,q+p]}\}}
  \prod_{u=1}^{i} \cont{\{\ur{u,q+u}\}} \Alt{\{\ur{[1,q]}\}} (\otimes{Q}) \otimes \Alt{\{\ur{[q+1,q+p]}\}} (\otimes{P}) \\
&\quad  =\sum_{\substack{ \sigma \in \frak{S}_q \\ \tau \in 
\frak{S}_p[q+1,q+p]}}
\dfrac{1}{q! p!} \sign(\sigma) \sign(\tau) 
\Alt{\{\sigma(\ur{[i+1,q]}), \tau(\ur{[q+i+1,q+p]})\}}
\prod_{u=1}^{i} \cont{\{\sigma(\ur{u}),\tau(\ur{q+u}) \}} (\otimes{Q}) \otimes (\otimes{P})
\end{align*}

Here if $Q=\{\ell_{1} \ell_{2} ,\ldots, \ell_{q}\}$,  $(\otimes{Q})$ denotes the ordinary tensor product (not exterior one) $\otimes \{\ur{\ell_{1} \ell_{2} ,\ldots, \ell_{q}}\}$. 

So in the  formula \eqref{eqn prpr}, the $i$ term  in the right-hand side 
 corresponds to the homomorphism in the direct summand $\Hom_{Pin(N)}(\Delta 
\bigotimes \bigwedge^{p+q-2i} V, \Delta)$.

So we can put 
\begin{equation*}
\begin{split}
&\pr{\{\ur{[1,q]}\}} \circ \pr{\{\ur{[q+1,q+p]}\}} = \\
&\sum_{i=0}^{\min(p,q)} d_i
\sum_{\substack{\sigma \in \frak{S}_q \\ \tau \in \frak{S}_p[q+1,q+p]}}
\sign(\sigma) \sign(\tau) \dfrac{\pr{\{\sigma(\ur{[i+1,q]}),\tau(\ur{[q+i+1, 
q+p]})\}}}{(q-i)! (p-i)!}
\dfrac{ \prod_{u=1}^{i} {\cont{\{\sigma(\ur{u}),\tau(\ur{q+u})\}}}}{i!}. \\
\end{split}
\end{equation*}

The rest of the proof is almost similar to \eqref{eqn injinj} and we omit it.

Next we prove the formula \eqref{eqn prinj4}.

Both sides of the above can be considered  as the homomorphisms from $\Delta 
\bigotimes\bigwedge^{p-t} V$ to  $\Delta \bigotimes \bigwedge^{q} V$.

So we can put (see Corollary \ref{cor phii}.)
\begin{equation*}
 \pr{\{\ur{[q+1,q+p]}\}} \circ \inj{\{\ur{[1, q+t]}\}}  = \sum_{i=0}^{\min(p-t,q)} d_i \;
{}^i\! \phi^{\ur{[q+t+1,q + p]}}_{\ur{[1,q]}}.
\end{equation*}

We show $d_i =(-1)^{(q-i)(p-i)+it} (\sum_{u=0}^{i} \binom{i}{u} (2n+1-p-q+t+i-
u)_{t})$ for $Pin(2n+1)$ and $d_i =(-1)^{p q + t (p+q)} (\sum_{u=0}^{i} \binom{i}{u} (2n+1-p-q+t+i-u)_{t})$ for $Pin(2n)$.

As before, we find a sequence of pairs $\{x_i, y_i\}$ and apply both sides of the above to $x_i$'s and compare the coefficients of $y_i$'s and 
 determine the coefficients $d_i$ inductively.

We prove the formula \eqref{eqn prinj4} for $Pin(2n+1)$.

We take the cases for the parities of $p-t$ and $q$.

Case 1. \quad $p-t=2r$ and $q=2s$.

We note $r+s \leqq n$. 

Case 1-0. \quad The determination of $d_{0}$.

We take $y_0 = [\ur{\emptyset}] \otimes \{\ur{[1,s]}, \ulb{[s, 1]}\}$, $x_0 = [\ur{\emptyset}] \otimes \{ \ur{[n-r+1, n]}, \ulb{[n, n-r+1]}\}$.  

In the left-hand side, from the definition \ref{def inj}, we have 

\begin{equation}
\begin{split}
& \inj{q+t}( [\ur{\emptyset}]) = 
 \sum_{\substack{ \tJ \subseteq [1,n] \\ \tW \subseteq ([1,n]-\tJ) \\
 j + 2 w =q+t}}
(-1)^{w} 2^{j/2} [ \ur{J}] 
  \otimes (q+t)! <\ur{J}, \ur{W}, \ulb{W} > \\
& \qquad \quad + \sum_{\substack{ \tJ \subseteq [1,n] \\ \tW \subseteq ([1,n] -\tJ) 
\\
 j + 2 w +1 =q+t}}
 2^{j/2} [ \ur{J}] 
  \otimes (q+t)! <\ur{J}, \ur{W}, \ur{0'}, \ulb{W} >.
\end{split}
\label{eqn injq+t}
\end{equation}

In the above only the $\tW$'s which include the set $\{[1,s]\}$ contribute to the 
coefficient and for them, we put $\tW_1 = \tW - \{[1,s]\}$.
The exterior products can be rewritten as 
$
(q+t)! <\ur{J}, \ur{W}, \ulb{W} > =
\otimes \{ \ur{[1,s]}, \ulb{[s,1]} \} \otimes t! <\ur{J}, \ur{W_1}, \ulb{W_1} >
 + \text{the other terms} $
 and 
$
(q+t)! <\ur{J}, \ur{W}, \ur{0'}, \ulb{W} > =
(-1)^s \otimes \{ \ur{[1,s]}, \ulb{[s,1]} \} \otimes t! <\ur{J}, \ur{W_1}, \ur{0'}, \ulb{W_1} > + \text{the other terms}$.

We tensor $\otimes \{\ur{[n-r+1, n]}, \ulb{[n, n-r+1]}\}$ after the exterior products in the above and alternate them. 
If $(\tJ \cup \tW_1)\cap \{[n-r+1, n]\} \ne \emptyset$, then the resulting exterior products are 0. So we can assume  $(\tJ \cup \tW_1)\cap \{[n-r+1, n]\} = \emptyset$.
Then we have $t! <\ur{J}, \ur{W_1}, \ur{[n-r+1, n]}, \ulb{[n, n-r+1]}, \ulb{W_1} >$ and $t! (-1)^r <\ur{J}, \ur{W_1}, \ur{[n-r+1, n]}, \ur{0'}, \ulb{[n, n-r+1]}, \ulb{W_1} >$.  

So the coefficient in the left-hand side is equal to 

\begin{align*}
& \sum_{\substack{\{[s+1, n-r]\} \supseteq \tW_1 \sqcup  \tJ \:\text{(disjoint)}  \\
 j + 2 w_1 =t}}
(-1)^{w_1+s} 2^{j/2} t! 
2^{j/2} (-1)^{w_1+r} + \sum_{\substack{\{[s+1, n-r]\} \supseteq \tW_1 \sqcup  \tJ \:\text{(disjoint)}  \\
 j + 2 w_1+1 =t}}
 2^{j/2} t!
 (-1)^{s+r}  2^{j/2} \\
& = (-1)^{r+s}  t! ((1+2x+x^2)^{n-r-s}[t] 
+ (1+2x+x^2)^{n-r-s}[t-1]) \\
&=(-1)^{r+s}  t! (1+x)^{2n+1-p+t-q}[t]=(-1)^{r+s} (2n+1-p+t-q)_{t}.
\end{align*}

Here for the formal power series $f(x)$, by $f(x)[t]$ we denote the coefficient of the power $x^t$ in $f(x)$.

In the right-hand side, only the $i=0$ term contributes to the coefficient and it is given by $(-1)^{r+s} d_0$ and we have $d_0=(2n+1-p+t-q)_{t}$.

Case 1-1. \quad The determination of $d_{1}$.

We take $y_1 = [\ur{s+1}] \otimes \{\ur{[1,s+1]}, \ulb{[s, 2]}\}$, $x_1 = [\ur{\emptyset}] \otimes \{\ur{1}, \ur{[n-r+2, n]}, \ur{0'},\ulb{[n, n-r+2]}\}$.

We calculate the left-hand side.

In the equation \eqref{eqn injq+t}, only the first sum contributes to the coefficient and moreover 
in the first sum the $\tW$'s and $\tJ$'s  such that  $\tW \supseteq 
\{2,\ldots,s\}$ and $\tJ \cup \tW \supseteq \{1,s+1\}$  contribute to the 
coefficient. 
There are four subcases, 1) $\tW \supseteq \{1,s+1\}$ (if $t \geqq 2$) \quad 2) $\tJ 
\supseteq \{1\}$, $\tW \supseteq \{s+1\}$ \quad  3) $\tJ \supseteq \{s+1\}$, $\tW \supseteq \{1\}$ \quad 4)  $\tJ \supseteq \{1,s+1\}$.
We put $\tJ_1=\tJ - \tJ \cap \{1,s+1\}$ and $\tW_1=\tW - \tW \cap [s+1]$. 

In the subcase 1), we have the coefficient 
\begin{equation*}
\sum_{\substack{ \{[s+2, n-r+1]\} \supseteq \tJ_1 \sqcup \tW_1 \:\text{(disjoint)} 
\\
 j_1 +2+2 w_1 =t}}
- 2^{j_1 +1/2} t!   =-\sqrt{2} t! (1+x)^{2n-2r-2s}[t-2].
\end{equation*}

In the subcase 2) we have the coefficient 
$${ 
\sum_{\substack{ \{[s+2, n-r+1]\} \supseteq \tJ_1 \sqcup \tW_1 \:\text{(disjoint)} 
\\
 j_1 +1+2 w_1 =t}}
- 2^{j_1+1 +1/2} t! =-2 \sqrt{2} t! (1+x)^{2n-2r-2s}[t-1].
}$$
 
In the subcase 3) we have the coefficient 
$${ 
\sum_{\substack{ \{[s+2, n-r+1]\} \supseteq \tJ_1 \sqcup \tW_1 \:\text{(disjoint)} 
\\
 j_1 +1+2 w_1 =t}}
- 2^{j_1 +1/2} t! =- \sqrt{2} t! (1+x)^{2n-2r-2s}[t-1].
}$$
 
In the subcase 4) we have the coefficient 
$${
\sum_{\substack{ \{[s+2, n-r+1]\} \supseteq \tJ_1 \sqcup \tW_1 \:\text{(disjoint)} 
\\
 j_1 +2 w_1 =t}}
- 2^{j_1+1 +1/2} t! =-2 \sqrt{2} t! (1+x)^{2n-2r-2s}[t].
}$$

So the total coefficient of the left-hand side is 
$-\sqrt{2} t!( (1+x)^{2n+1-2r-2s}[t] + (1+x)^{2n+1-2r-2s+1}[t]
=-\sqrt{2}((2n+1-p-q+t)_{t} + (2n+1-p-q+t+1)_{t})$.
If $t=1$, the subcase 1) does not occur, but the result also holds in this case  since 
$(1+x)^{2n-2r-2s}[-1]=0$.

In the right-hand side, only the $i=1$ term contributes to the coefficient and it is given by $\sqrt{2} d_1$ and we have $d_1=-((2n+1-p-q+t)_{t} + (2n+1-p-q+t+1)_{t})$.

Case 1-2. \quad The determination of $d_{2}$.

We take $y_2 =[\ur{\emptyset}] \otimes \{[1,s+1], \overline{[s+1, 3]}\}$,  $x_2 =[\ur{\emptyset}] \otimes \{1,2, [n-r+2, n], \overline{[n, n-r+2]}\}$.

We calculate the left-hand side.

In the equation \eqref{eqn injq+t}, only the $\tW$'s and $\tJ$'s such that  $\tW \supseteq \{[3, s+1]\}$ and $\tJ \cup \tW \supseteq \{1,2\}$  contribute to the 
coefficient.
There are four subcases, 1) $\tW \supseteq \{1,2\}$ (if $t \geqq 2$) \quad 2) $\tJ \supseteq \{1\}$, $\tW \supseteq \{2\}$ \quad  3) $\tJ \supseteq \{2\}$, $\tW 
\supseteq \{1\}$ \quad 4)  $\tJ \supseteq \{1,2\}$.
We put $\tJ_1=\tJ - \tJ \cap \{1,2\}$ and $\tW_1=\tW - \tW \cap [s+1]$. 

In the subcase 1) we have the coefficient 
\begin{equation*}
\sum_{\substack{ \{[s+2, n-r+1]\} \supseteq \tJ_1 \sqcup \tW_1 \:\text{(disjoint)} 
\\
 j_1 +2+2 w_1 =t}}
(-1)^{s+r} 2^{j_1} t! 
+\sum_{\substack{ \{[s+2, n-r+1]\} \supseteq \tJ_1 \sqcup \tW_1 \:\text{(disjoint)} 
\\
 j_1 +2+2 w_1+1 =t}}
(-1)^{s+r} 2^{j_1} t!.
\end{equation*}

In the subcase 2) we have the coefficient 
\begin{equation*}
\sum_{\substack{ \{[s+2, n-r+1]\} \supseteq \tJ_1 \sqcup \tW_1 \:\text{(disjoint)} \\
 j_1 +1+2 w_1 =t}}
(-1)^{s+r} 2^{j_1+1} t! 
+\sum_{\substack{ \{[s+2, n-r+1]\} \supseteq \tJ_1 \sqcup \tW_1 \:\text{(disjoint)} \\
 j_1 +2+2 w_1 =t}}
(-1)^{s+r} 2^{j_1+1} t!. 
\end{equation*}

In the subcase 3) we have the coefficient 
\begin{equation*}
\sum_{\substack{ \{[s+2, n-r+1]\} \supseteq \tJ_1 \sqcup \tW_1 \:\text{(disjoint)} \\
 j_1 +1+2 w_1 =t}}
(-1)^{s+r} 2^{j_1+1} t! 
+\sum_{\substack{ \{[s+2, n-r+1]\} \supseteq \tJ_1 \sqcup \tW_1 \:\text{(disjoint)} \\
 j_1 +2+2 w_1 =t}}
(-1)^{s+r} 2^{j_1+1} t!. 
\end{equation*}

In the subcase 4) we have the coefficient 
\begin{equation*}
\sum_{\substack{ \{[s+2, n-r+1]\} \supseteq \tJ_1 \sqcup \tW_1 \:\text{(disjoint)} 
\\
 j_1 +2 w_1 =t}}
(-1)^{s+r} 2^{j_1+2} t! 
+\sum_{\substack{ \{[s+2, n-r+1]\} \supseteq \tJ_1 \sqcup \tW_1 \:\text{(disjoint)} 
\\
 j_1+2 w_1+1 =t}}
(-1)^{s+r} 2^{j_1+2} t!. 
\end{equation*}

So the total coefficient in the left-hand side is given as follows.

$
(-1)^{s+r} t!\{
(1+x)^{2n-p-q+t}[t-2] +2 (1+x)^{2n-p-q+t}[t-1] +2 (1+x)^{2n-p-q+t}[t-1] +4 
(1+x)^{2n-p-q+t}[t] +(1+x)^{2n-p-q+t}[t-3] +2 (1+x)^{2n-p-q+t}[t-2] +2 
(1+x)^{2n-p-q+t}[t-2] +4 (1+x)^{2n-p-q+t}[t-1] \}.
$

Here we prove the following lemma.
\begin{lem}
For any indeterminate $X$ and for any non-negative integers $t$ and $i$, we have$$
 \sum_{u=0}^{i} \binom{i}{u} 2^u \binom{X}{t-i+u} = \sum_{u=0}^{i} \binom{i}{u}
\binom{X+i-u}{t}.
$$
\end{lem}
\begin{proof}
This follows easily from the formula
$$
\binom{X+i}{t}= \sum_{u=0}^{i} \binom{X}{t-u} \binom{i}{u}.
$$

\end{proof}

From this lemma, the coefficient of the left-hand side is 
$
(-1)^{s+r} \{(2n+3-p-q+t)_{t} +2 (2n+2-p-q+t)_{t}  +(2n+1-p-q+t)_{t} \}.
$

In the right-hand side, only the $i=2$ term contributes to the coefficient and it is given by $(-1)^{s-1+r-1}d_2$.
So $d_2 = (2n+3-p-q+t)_{t} +2 (2n+2-p-q+t)_{t}  +(2n+1-p-q+t)_{t}$.

Case 1-$2i+1$. \quad The determination of $d_{2i+1}$.

We take $y_{2i+1} = [\ur{s+i+1}] \otimes \{\ur{[1,s+i+1]}, \ulb{[s+i, 2i+2]}\}$,  $x_{2i+1} = [\ur{\emptyset}] \otimes \{\ur{[1, 2i+1]}, \ur{[n-r+i+2, n]}, \ur{0'}, \ulb{[n, n-r+i+2]}\}$.

In the equation \eqref{eqn injq+t}, only the $\tW$'s and $\tJ$'s such that  $\tW \supseteq \{[2i+2, s+i]\}$ and $\tJ \cup \tW \supseteq \{[1, 2i+1],s+i+1\}$  contribute to the coefficient.

We put $\tJ_2=\tJ \cap [1,2i+1]$ and $\tJ_1=\tJ \cap \{[s+i+2, n-r+i+1]\}$ and 
$\tW_2=\tW \cap [1,2i+1]$ and $\tW_1=\tW \cap \{[s+i+2, n-r+i+1]\}$ . 
We fix a subset $\tJ_2 \subseteq [1,2i+1]$ with $\vert \tJ_2 \vert = v$ and  take cases for $s+i+1 \in \tJ$ and $s+i+1 \notin \tJ$.

Subcase  1-$2i+1$-1. \quad $s+i+1 \notin \tJ$.

In this subcase we have the coefficient
\begin{align*}
& -t!  \sum_{\substack{\{[s+i+2, n-r+i+1]\} \supseteq \tW_1 \sqcup  \tJ_1 
\:\text{(disjoint)}  \\
 j_1 + 2 w_1+2i+2-v =t}}
2^{j_1+v+1/2} \\
&=-t!  2^{v+1/2} \sum_{ j_1 + 2 w_1+2i+2-v =t}\binom{n-r-s}{j_1}\binom{n-r-s-
j_1}{w_1} 2^{j_1}\\
&=-t! 2^{v+1/2}(1+2x+x^2)^{n-r-s}[t-2i-2+v]. \\
\end{align*}

Subcase  1-$2i+1$-2. \quad $s+i+1 \in \tJ$.

In this subcase we have the coefficient
\begin{align*}
& -t!  \sum_{\substack{\{[s+i+2, n-r+i+1]\} \supseteq \tW_1 \sqcup  \tJ_1 
\:\text{(disjoint)}  \\
 j_1 + 2 w_1+2i+1-v =t}}
2^{j_1+v+1/2} \\
&=-t! 2^{v+1/2}(1+2x+x^2)^{n-r-s}[t-2i-1+v].\\
\end{align*}

Therefore for a fixed $\tJ_2 \subseteq [2i+1]$, 
we have the coefficient $-t! 2^{v+1/2}(1+x)^{2n+1-p-q+t}[t-2i-1+v]$ and the total coefficient is 
\begin{align*}
&-t! \sqrt{2} \sum_{v=0}^{2i+1}\binom{2n+1-p-q+t}{t-2i-1+v}2^{v}
\binom{2i+1}{v} \\
&\qquad \qquad = -t! \sqrt{2} \sum_{v=0}^{2i+1} \binom{2i+1}{v}
\binom{2n+1-p-q+t+2i+1-v}{t}.\\
\end{align*}

In the right-hand side, only the $2i+1$ term contributes to the coefficient and it is given by $\sqrt{2}d_{2i+1}$.
So $d_{2i+1}=-\sum_{v=0}^{2i+1} \binom{2i+1}{v}
(2n+1-p-q+t+2i+1-v)_{t}$.

Case 1-$2i$. \quad The determination of $d_{2i}$.

We take $y_{2i} = [\ur{\emptyset}] \otimes \{\ur{[1, s+i]}, \ulb{[s+i, 2i+1]}\}$,\newline $x_{2i} 
= [\ur{\emptyset}] \otimes \{\ur{[1, 2i]}, \ur{[n-r+i+1, n]}, \ulb{[n, n-r+i+1]}\}$.

In the equation \eqref{eqn injq+t}, only the $\tW$'s and $\tJ$'s such that  $\tW \supseteq \{2i+1,\ldots,s+i\}$ and $\tJ \cup \tW \supseteq \{[1,2i]\}$  contribute to the coefficient.

We put $\tJ_2=\tJ \cap [1,2i]$ and $\tJ_1=\tJ \cap \{[s+i+1, n-r+i]\}$ and $\tW_2=\tW \cap [1,2i]$ and $\tW_1=\tW \cap \{[s+i+1, n-r+i]\}$ . 
We fix a subset $\tJ_2 \subseteq [1,2i]$ with $\vert \tJ_2 \vert = v$.
We take the cases.

Subcase  1-$2i$-1. \quad The terms do not contain $0'$.

In this subcase we have the coefficient
\begin{align*}
& t!  \sum_{\substack{\{[s+i+1, n-r+i]\} \supseteq \tW_1 \sqcup  \tJ_1 
\:\text{(disjoint)}  \\
 j_1 + 2 w_1+2i-v =t}}
(-1)^{s-i+r-i} 2^{j_1+v} \\
&=(-1)^{s+r} t! 2^{v} \sum_{ j_1 + 2 w_1+2i-v =t}\binom{n-r-s}{j_1}\binom{n-r-s-j_1}{w_1}2^{j_1} \\
&=(-1)^{s+r} t! 2^{v}(1+x)^{2n-p-q+t}[t-2i+v]. \\
\end{align*}

Subcase  1-$2i$-2. \quad The terms contain $0'$.

In this subcase we have the coefficient
\begin{align*}
& t! 2^{v}(-1)^{s+r} \sum_{\substack{\{[s+i+1, n-r+i]\} \supseteq \tW_1 \sqcup  \tJ_1 \:\text{(disjoint)}  \\
 j_1 + 2 w_1+2i+1-v =t}}
2^{j_1} \\
&= t! 2^{v}(-1)^{s+r}(1+2x+x^2)^{n-r-s}[t-2i-1+v].\\
\end{align*}

Therefore for a fixed $\tJ_2 \subseteq [1,2i]$, 
we have the coefficient $t! 2^{v}(-1)^{s+r}(1+x)^{2n+1-p-q+t}[t-2i+v]$ and the total coefficient is 
\begin{align*}
&t! (-1)^{s+r} \sum_{v=0}^{2i}\binom{2n+1-p-q+t}{t-2i+v}2^{v}
\binom{2i}{v} \\
& \qquad \qquad = t! (-1)^{s+r} \sum_{v=0}^{2i} \binom{2i}{v}
\binom{2n+1-p-q+t+2i-v}{t}. \\
\end{align*}

In the right-hand side, only the $2i$ term contributes to the coefficient and it is given by $(-1)^{s+r}d_{2i}$. 
So $d_{2i}=\sum_{v=0}^{2i} \binom{2i}{v}
(2n+1-p-q+t+2i-v)_{t}$.

Therefore in this case we have $d_{i}=(-1)^i \sum_{v=0}^{i} \binom{i}{v}
(2n+1-p-q+t+i-v)_{t}$.

The proof goes in a similar way for the remaining cases and so we omit it.
 $d_i$ is given by $ (-1)^{(q-i)(p-i)+it} \sum_{u=0}^{i}  \binom{i}{u}
(2n+1-p-q+t+i-u)_{t}$ for $Pin(2n+1)$ and 
 $d_i$ is given by $ (-1)^{ q p + t(p+q)} \sum_{u=0}^{i} \binom{i}{u}
(2n-p-q+t+i-u)_{t}$ for $Pin(2n)$.

Next we prove the formula \eqref{eqn continj}.

Both sides of the formula can be considered as the homomorphisms from $\Delta \bigotimes \bigotimes^{q} V$ to  $\Delta \bigotimes \bigotimes^{p} V$.

So we can put  
\begin{equation*}
\prod_{i=1}^{q} \cont{\{\ur{i,p+q+i}\}}  \inj{\{\ur{[1, q+p]} \}} = 
\sum_{i=0}^{\min(p,q)} d_i\;
{}^i\! \phi^{\ur{[p+q+1,p + 2q]}}_{\ur{[q+1,q+p]}}
\end{equation*}

As before, we find a sequence of pairs $\{x_i, y_i\}$ and apply both sides of the above to $x_i$'s and compare the coefficients of $y_i$'s and 
 determine the coefficients $d_i$ inductively.

We determine $d_i$ for $Pin(2n+1)$.

Case 1 \quad $p=2r$ and $q=2s$.

We note $r+s \leqq n$.

Case 1-$2i$. \quad The determination of $d_{2i}$.

We take $y_{2i} = [\ur{\emptyset}] \otimes \{\ur{[1,2i]},\ur{[s+i+1,s+r]}, \ulb{[s+r, s+i+1]}\}$,  $x_{2i} = [\ur{\emptyset}] \otimes \{\ur{[1, s+i]}, \ulb{[s+i, 2i+1]}\}$. 
In this case $d_{2i}=(-1)^{\binom{q}{2}}$.

Case 1-$2i+1$ \quad The determination of $d_{2i+1}$.

We take $y_{2i+1} =[\ur{s+i+1}] \otimes \{\ur{[1, 2i+1]}, \ur{[s+i+1, s+r]}, \ulb{[s+r, s+i+2]}\}$, $x_{2i+1} =[\ur{\emptyset}] \otimes \{\ur{[1, s+i]}, \ur{0'}, \ulb{[s+i, 2i+2]}\}$. 
In this case $d_{2i+1}=(-1)^{s+1}$.

The proof goes in a similar way for the remaining cases and so we omit it.
Summarizing all the calculation, we have $d_i =(-1)^{p q +\binom{q}{2} +i (p+q-
1)}$ for $Pin(2n+1)$.

For $Pin(2n)$, a similar calculation gives us $d_i =(-1)^{\binom{q}{2}}$.

Finally we prove the formula \eqref{eqn prid}.

Both sides of the formula can be considered as the homomorphisms from $\Delta \bigotimes 
\bigotimes^{q} V$ to  $\Delta \bigotimes \bigotimes^{p} V$.

So we can put
\begin{equation*}
 \pr{\{\ur{[1, q+p]} \}} \prod_{i=1}^{p} {\id_V}_{\{\ur{q+i,p+q+i}\}} = 
\sum_{i=0}^{\min(p,q)} d_i\;
{}^i\! \phi^{\ur{[1, q]}}_{\ur{[p+q+1,q+2p]}}.
\end{equation*}

As before, we find a sequence of pairs $\{x_i, y_i\}$ and apply both sides of the above to $x_i$'s and compare the coefficients of $y_i$'s and 
 determine the coefficients $d_i$ inductively.

We determine $d_i$ for $Pin(2n+1)$.

Case 1. \quad $p=2r$ and $q=2s$.

We note $r+s \leqq n$. 

Case 1-$2i$. \quad The determination of $d_{2i}$.

We take $y_{2i} = [\ur{\emptyset}] \otimes \{\ur{[1, 2i]}, \ur{[s+i+1, s+r]}, \ulb{[s+r, s+i+1]}\}$, $x_{2i} = [\ur{\emptyset}] \otimes \{\ur{[1, s+i]}, \ulb{[s+i, 2i+1]}\}$.

In this case $d_{2i}=(-1)^{r}$.

Case 1-$2i+1$. \quad The determination of $d_{2i+1}$.

We take $y_{2i+1} = [\ur{s+i+1}] \otimes \{\ur{[1, 2i+1]}, \ur{[s+i+1, s+r]}, \ulb{[s+r, s+i+2]}\}$, $x_{2i+1} = [\ur{\emptyset}] \otimes \{\ur{[1, s+i]}, \ur{0'}, \ulb{[s+i, 2i+2]}\}$.

In this case $d_{2i+1}=(-1)^{r-1}$.

The proof goes in a similar way for the remaining cases and so we omit it.
Summarizing all the calculation, we have $d_i =(-1)^{\binom{p}{2} +i (p+q+1)}$ for $Pin(2n+1)$.

For $Pin(2n)$, a similar calculation gives us $d_i =(-1)^{\binom{p}{2} + p q}$.

\end{proof}

\section{Examples of products of the generalized Brauer diagrams}\label{sect exam}

From the result of the previous section, 
we can calculate the product of the generalized Brauer diagrams.

As we remarked after the statement of Theorem \ref{thm reprel}, 
we change $N$ in the  formulas of the theorem into an indeterminate $X$ simultaneously and we write down the relations of  the \lq generic' centralizer algebra of $\mathbf{CP_k}$.

We summarize other relations between the contractions and the immersions and 
$\pr{}$ and $\inj{}$ which follow easily from the definitions.
 We fix the index set of the tensor positions of $\bigotimes^k V$ from $1$ to $k$ and after the contraction $\cont{\{\ur{i,j}\}}$ we consider the  $i$-th and the $j$-th components are occupied by the empty set and ${\id_V}_{\{\ur{i,j}\}}$ is allowed if the positions $i$ and $j$ are occupied by the empty set.

Then we have the followings.
\begin{lem}\label{lem releasy}

\begin{enumerate}
\item
\[
\cont{\{\ur{s,t}\}} {\id_V}_{\{\ur{i,j}\}}  = {\id_V}_{\{\ur{i,j}\}}  \cont{\{\ur{s,t}\}} \qquad 
\text{ if $\{s,t\} \cap \{i,j\} = \emptyset$.}
\]
\item
\[
\cont{\{\ur{s,t}\}} {\id_V}_{\{\ur{t,j}\}} = \begin{pmatrix} \ur{s} \\ \ur{j} \\ \end{pmatrix}
\qquad \text{ if $s \ne j$.}
\]
Here $\begin{pmatrix} \ur{s} \\ \ur{j} \\ \end{pmatrix}$ denotes the (partial) permutation 
which sends $s$-th component to the $j$-th position.

\item 
\[
\cont{\{\ur{s,t}\}} {\id_V}_{\{\ur{s,t}\}} = N\; \id.
\]
\item 
If $s, t \in \tT$, we have  
\[
\cont{\{\ur{s,t}\}} \inj{\tT} = 0, \qquad \pr{\tT} {\id_{V}}_{\{\ur{s,t}\}}= 0.
\]

\item 
If $s, t_1, t_2, \ldots, t_r \in [1,n]$  are different from each other, we have  
\[
\begin{pmatrix} \ur{t_i} \\ \ur{s} \\ \end{pmatrix} \inj{\{\ur{t_1,\ldots, t_i, \ldots, t_r}\}} = 
\inj{\{\ur{t_1},\ldots, \overset{i}{\ur{s}}, \ldots, \ur{t_r}\}},
\quad
\begin{pmatrix} \ur{t_1} \\ \ur{s} \\ \end{pmatrix} {\id_{V}}_{\{\ur{t_1,t_2}\}} = 
{\id_{V}}_{\{\ur{s,t_2}\}}.
\]

\item 
If $s, t_1, t_2, \ldots, t_r \in [1,n]$  are different from each other, we have  
\[
\pr{\{\ur{t_1,\ldots, t_i, \ldots, t_r}\}} \begin{pmatrix} \ur{s} \\ \ur{t_i} \\ \end{pmatrix} = 
\pr{\{\ur{t_1},\ldots, \overset{i}{\ur{s}}, \ldots, \ur{t_r}\}},
\quad
\cont{\{\ur{t_1, t_2}\}} \begin{pmatrix} \ur{s} \\ \ur{t_1} \\ \end{pmatrix} = \cont{\{\ur{s, t_2}\}}.
\]

\item 
If $s, t, t_1, t_2, \ldots, t_r \in [1,n]$  are different from each other, 
${\id_{V}}_{\{\ur{s,t}\}}$  commutes with $\inj{\ur{T}}$ and $\pr{\ur{T}}$ and 
$\cont{\{\ur{s,t}\}}$ commutes with $\inj{\ur{T}}$ and $\pr{\ur{T}}$, where 
$\ur{T} =\{t_1, t_2, \ldots, t_r\}$.

\end{enumerate}
\end{lem}

\begin{proof}
The first three formulas are easy and the fourth formula follows directly from the alternating property with the indices in $\tT$ of $\pr{}$ and $\inj{}$. The rest are obvious from the definitions.
\end{proof}

\begin{rem}
If we put $e_{\ur{i}}= {\id_V}_{\{\ur{i,i+1}\}} \cont{\{\ur{i,i+1}\}}$, we can deduce easily $e_{\ur{i}} e_{\ur{i+1}} e_{\ur{i}} = e_{\ur{i}}$ from the above lemma. These relations were used to define the 
Birman-Wenzl algebras (q-analogs of the Brauer centralizer algebras).(\cite{biwen}.) 
\end{rem}

\begin{rem}
We note that 
\[
\inj{\{\ur{[1,q]}\}} \circ \inj{\{\ur{[q+1,q+p]}\}} \ne 
\inj{\{\ur{[q+1,q+p]}\}} \circ \inj{\{\ur{[1,q]}\}}
\]
and
\[
\pr{\{\ur{[1,q]}\}} \circ \pr{\{\ur{[q+1,q+p]}\}} \ne 
\pr{\{\ur{[q+1,q+p]}\}} \circ \pr{\{\ur{[1,q]}\}}.
\]
\end{rem}

Let us show some examples.

\begin{exmp}\label{exmp product}
In this example we always assume $n \geqq k$ and use the representation 
theoretical parametrization of the generalized Brauer diagrams, so we omit the subscript  $rt$  here.
We calculate the product of $y_5 y_8$ in figure \ref{spincentfig1}.

\begin{figure}[h]
\begin{minipage}{\linewidth}
\begin{center}
\input{spincentfig5.tex}
\end{center}
\end{minipage}
\caption{The result of the product of $y_5 y_8$ for $Pin(2n+1)$}
\label{fig y5y8}
\end{figure}

Here $y_8= \inj{\{\ur{1,2}\}} \cont{\{\ur{1,2}\}}$ and $y_5 = \inj{\{\ur{1}\}} \begin{pmatrix} \ur{2} 
\\ \ur{2} \\ \end{pmatrix} \pr{\{\ur{1}\}}$.
From the formula (\ref{eqn injpr2}), 
we have $\pr{\{\ur{1}\}}\inj{\{\ur{1,2}\}} =(X-1)_1\inj{\{\ur{2}\}}$ (we put $X$ for $2n+1$.) and 
the targeting homomorphism is $\inj{\{\ur{1}\}} \begin{pmatrix} \ur{2} \\ \ur{2} \\ \end{pmatrix} 
(X-1)\inj{\{\ur{2}\}} \cont{\{\ur{1,2}\}} = (X-1) \inj{\{\ur{1}\}} \inj{\{\ur{2}\}} \cont{\{\ur{1,2}\}}$.
From the formula (\ref{eqn injinj}), we have 
$\inj{\{\ur{1}\}} \inj{\{\ur{2}\}} = \inj{\{\ur{1,2}\}} + {\id_{V}}_{\{\ur{1,2}\}}$ and the final 
formula is given in figure \ref{fig y5y8}.

We calculate a more complicated example of figure \ref{fig calc}.

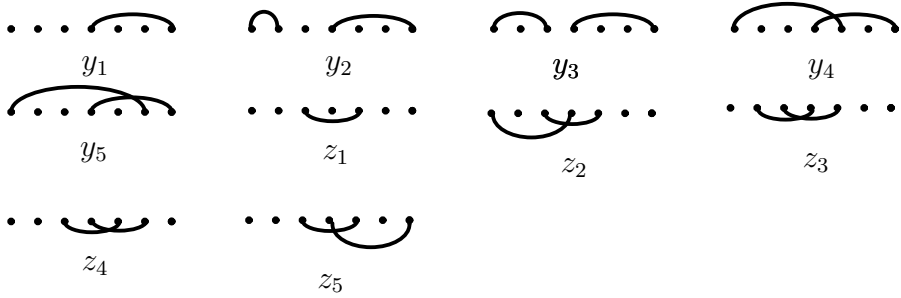
\begin{figure}[h]
\begin{minipage}{\linewidth}
\begin{center}
\input{spincentfig6.tex}
\end{center}
\caption{The result of the product of a more complicated example for $Pin(2n)$}
\label{fig calc}
\end{minipage}
\end{figure}

In this case we must calculate 
$$
{\id_{V}}_{\{\ur{3,5}\}}\inj{\{\ur{1,2,6,7}\}} 
 \begin{pmatrix} \ur{5} \\ \ur{4} \\ \end{pmatrix} \pr{\{\ur{1,2,4,7}\}}\cont{\{\ur{3,6}\}}
{\id_{V}}_{\{\ur{5,7}\}}\inj{\{\ur{1,2,4,6}\}} 
 \begin{pmatrix} \ur{1} \\ \ur{3} \\ \end{pmatrix} \pr{\{\ur{2,3,5,6}\}}\cont{\{\ur{4,7}\}}.
$$

The inside homomorphism of the above, i.e., the homomorphism 
obtained by taking the middle part between the contractions and 
$\pr{}$ of the first element
 and the immersion and $\inj{}$ of the second element, is given by 
$\begin{pmatrix} \ur{5} \\ \ur{4} \\ \end{pmatrix} \pr{\{\ur{1,2,4,7}\}}\cont{\{\ur{3,6}\}} {\id_{V}}_{\{\ur{5,7}\}}\inj{\{\ur{1,2,4,6}\}} 
 \begin{pmatrix} \ur{1} \\ \ur{3} \\ \end{pmatrix}$. 
This element is in $\mathbf{CP^1_1}= \Hom_{Pin(2n)}(\Delta \bigotimes V, \Delta \bigotimes V)$. $($ $V$ sits in the first place in the upper row and sits in the forth place in the lower row.$)$
 So it can be written as a sum of the basis of $\mathbf{CP^1_1}$ i.e., $\mathbf{GB^1_1}$. 
In general, this fact tells us that the linear spaces spanned by all the diagrams in $\mathbf{GB^k_k}$ with at most $s$ $(0 \leqq s \leqq k)$ vertical edges $($the edge from the upper row to the lower row$)$, become  two-sided ideals in $\mathbf{CP_k}$.

From the formula \eqref{eqn prid}, we have 
$$
\pr{\{\ur{1,2,3,4}\}}{\id_{V}}_{\{\ur{4,5}\}} = -\inj{\{\ur{5}\}} \pr{\{\ur{1,2,3}\}}
- \begin{pmatrix} \ur{1} \\ \ur{5} \\ \end{pmatrix} \pr{\{\ur{2,3}\}}
+ \begin{pmatrix} \ur{2} \\ \ur{5} \\ \end{pmatrix} \pr{\{\ur{1,3}\}}
- \begin{pmatrix} \ur{3} \\ \ur{5} \\ \end{pmatrix} \pr{\{\ur{1,2}\}}.
$$
If we apply the conjugation of the permutation $\eta = \begin{pmatrix} \ur{3} & \ur{4} \\ \ur{4} & \ur{7} \\ \end{pmatrix}$, we have 
$$
\pr{\{\ur{1,2,4,7}\}}{\id_{V}}_{\{\ur{5,7}\}} = -\inj{\{\ur{5}\}} \pr{\{\ur{1,2,4}\}}
- \begin{pmatrix} \ur{1} \\ \ur{5} \\ \end{pmatrix} \pr{\{\ur{2,4}\}}
+ \begin{pmatrix} \ur{2} \\ \ur{5} \\ \end{pmatrix} \pr{\{\ur{1,4}\}}
- \begin{pmatrix} \ur{4} \\ \ur{5} \\ \end{pmatrix} \pr{\{\ur{1,2}\}}.
$$

From the formula \eqref{eqn continj}, we have 
$$
\cont{\{\ur{1,5}\}} \inj{\{\ur{1,2,3,4}\}}=  \inj{\{\ur{2,3,4}\}}\pr{\{\ur{5}\}}
+ \inj{\{\ur{3,4}\}} \begin{pmatrix} \ur{5} \\ \ur{2} \\ \end{pmatrix} - 
\inj{\{\ur{2,4}\}} \begin{pmatrix} \ur{5} \\ \ur{3} \\ \end{pmatrix} + 
\inj{\{\ur{2,3}\}} \begin{pmatrix} \ur{5} \\ \ur{4} \\ \end{pmatrix}.
$$
If we apply the conjugation of the permutation $\eta = \begin{pmatrix} \ur{1} & \ur{2} & \ur{3} & \ur{4} & \ur{5} \\ \ur{6} & \ur{2} & \ur{4} & \ur{1} & \ur{3} \\ \end{pmatrix}$, we have 
$$
\cont{\{\ur{3,6}\}} \inj{\{\ur{1,2,4,6}\}}= - \inj{\{\ur{1,2,4}\}}\pr{\{\ur{3}\}}
+ \inj{\{\ur{1,4}\}} \begin{pmatrix} \ur{3} \\ \ur{2} \\ \end{pmatrix} - 
\inj{\{\ur{1,2}\}} \begin{pmatrix} \ur{3} \\ \ur{4} \\ \end{pmatrix} - 
\inj{\{\ur{2,4}\}} \begin{pmatrix} \ur{3} \\ \ur{1} \\ \end{pmatrix}.
$$

Since $\pr{\{\ur{2,3}\}}\inj{\{\ur{1,2}\}} = -(X-2) \inj{\{\ur{1}\}}\pr{\{\ur{3}\}} - \{(X-1-2+1+1-0) + (X-1-2+1+1-1)\}  \begin{pmatrix} \ur{3} \\ \ur{1} \\ \end{pmatrix}$ from the formula \eqref{eqn prinj4}, we have 
$\pr{\{\ur{2,4}\}}\inj{\{\ur{1,4}\}}= (X-2) \inj{\{\ur{1}\}}\pr{\{\ur{2}\}} + (2X-3)  \begin{pmatrix} \ur{2} \\ \ur{1} \\ \end{pmatrix}$. 
The calculation goes on in almost similar way and 
the final result consists of $26$ terms in which only one term contains a vertical line. It is given in figure \ref{fig calc}.
\end{exmp}

\section{Dual Pair and Spin representations}\label{sect spinrep}

In this section we define the subspace in the tensor space $T_{k} =\Delta \bigotimes \bigotimes^k V$, on which the symmetric group of degree $k$ and $Pin(N)$ (or $Spin(N)$) act as the dual pair.

By $\mathbf{GB_k}$, we denote $\mathbf{GB_k^k}$.

From  Lemma \ref{lem indep},  the elements of $\mathbf{GB_k}$ under the invariant theoretical parametrization with at most $N$ isolated vertices (we denote a set of these elements by $\mathbf{GB_{k, inv}}$.) span linearly $\mathbf{CP_k}$.
If $n \geqq k$,  they become a base of $\mathbf{CP_k}$ and in this case, the elements of $\mathbf{GB_k}$ under the representation theoretical parametrization (we denote a set of these elements by $\mathbf{GB_{k, rt}}$.) also become a base of $\mathbf{CP_k}$.

Since $\mathbf{CP_k}$ contains $\R[\frak{S}_k]$ naturally,  For $0 \leqq s \leqq \min(k,n)$, 
let us define the subspace $T_{k,s}^{0}$ of $\Delta \bigotimes \bigotimes^k V$ as follows.

\begin{defn}\label{def T0}
Let $T_{k,s}^{0}$ denote the intersection of the all the kernels of the contractions $\cont{\{\ur{i,j}\}}$ ( $ 1 \leqq i < j \leqq k $) and the projections $\pr{\{\ur{i_1,i_2, \ldots, i_r}\}}$ ($r > 0$ and $ 1 \leqq i_1 < i_2 < \ldots < i_r \leqq k $) 
and the alternating operators  $\Alt{\{\ur{i_1,i_2, \ldots, i_r}\}}$ 
(see Definition \ref{def alt} and $r > s $) of degree greater than $s$. 
\end{defn}

If $n \geqq k$, then we ignore the alternating operators in the above and put $T_{k,k}^{0} =T_{k}^{0}$.

\begin{lem}
The space $T_{k,s}^{0}$ becomes a $\mathbf{CP_k} \times Pin(N)$ module.
\end{lem}

\begin{proof}
It is obvious that $T_{k,s}^{0}$ is a $Pin(N)$ module from the definition.

Let $y_{inv} \in \mathbf{GB_{k, inv}}$ and let $\psi^{\ur{T_{u}}}_{\ur{T_{\ell}}}$ be the homomorphism corresponding to the isolated vertices in $y$.

Then if $\vert \ur{T_{u}} \vert > s$, 
$\psi^{\ur{T_{u}}}_{\ur{T_{\ell}}}$ acts on the space $T_{k,s}^{0}$ by $0$ and 
so $y_{inv} T_{k,s}^{0} =0$.

We can assume $\vert \ur{T_{u}} \vert \leqq s =\min(k,n)$. If $\vert \ur{T_{\ell}} \vert \ne \vert \ur{T_{u}} \vert$, the contractions or $\pr{}$ must appear in the upper row, so $y_{inv}\; T_{k,s}^{0} =0$ .(See \eqref{eqn psitorep}.)

We can assume $\vert \ur{T_{u}} \vert =\vert \ur{T_{\ell}} \vert \leqq s =\min(k,n)$. Then from Theorem \ref{thm psitorep}, $\psi^{\ur{T_{u}}}_{\ur{T_{\ell}}}$
 acts on this space by the alternating operators of degree $\vert \ur{T_{u}} \vert$.
Since the space $T_{k,s}^{0}$ is invariant under the action of 
$\frak{S}_k$, the lemma holds.
\end{proof}

Let us denote the set of the annihilators in $\mathbf{CP_k}$ of the module $T_{k,s}^{0}$ by $\frak{J}_{k,s}$. Then from the above, $\frak{J}_{k,s}$ becomes a two sided ideal in $\mathbf{CP_k}$ and $\mathbf{CP_k}/\frak{J}_{k,s}$ is the surjective image of the subalgebra $\R[\frak{S}_k]$.
Since $Pin(N)$ and $\mathbf{CP_k}$ are dual pair on $T_{k}$, on this direct summand $T_{k,s}^{0}$, the images of $\R[\frak{S}_k]$ and $Pin(2n)$ become a dual pair.

So we have the following theorem.

\begin{thm}\label{thm dualpin}
The symmetric group of degree $k$ and $Pin(N)$ become a dual pair 
on the subspace $T_{k,s}^{0}$ of $\Delta_{(\ep_1)} \bigotimes \bigotimes^k V$. 
Namely the space $T_{k,s}^{0}$ is decomposed into the direct sum of the tensor 
products of the irreducible representations with multiplicity free as follows.
\begin{equation}
T_{k,s}^{0} = \sum_{\substack{\lam : \text{partitions of size $k$} \\ \ell(\lam) \leqq s }} \lam_{\frak{S}_k} \textstyle{\bigotimes} [\Delta, \lam]_{Pin(N), (\ep_2)}.
\tag{\ref{thm dualpin}.1}\label{eqn dualpin}
\end{equation}
Here for $Pin(2n+1)$, $\ep_1 (-1)^k = \ep_2$
\end{thm}

\begin{rem}\label{rem dualspin}
For $Spin(2n+1)$, on the above $T_{k,s}^{0}$,  $\frak{S}_k$ and $Spin(2n+1)$ also act as a dual pair, since the center $z$ acts on it by the scalar $\ep_1 (-1)^k$. The decomposition is the same if we ignore $\ep_1$ and $\ep_2$.
\end{rem}

\begin{rem}\label{rem dualpin}
From the above, if $n \geqq k$, $\frak{J}_{k,k}$ is the linear space spanned by $y_{rt} \in \mathbf{GB_{k, rt}}$ 
with at most $k-1$ vertical edges and in this case, we have
$\mathbf{CP_k} = \R[\frak{S}_k] \bigoplus \frak{J}_{k,k}$.
\end{rem}

\begin{proof}
From the argument just before this theorem, 
$\R[\frak{S}_k]$ and $Pin(N)$ become a dual pair on the space $T_{k,s}^{0}$.

The Schur-Weyl reciprocity tells us that, 
as $\frak{S}_k \times GL(N)$ module,  $\bigotimes^k V$ is decomposed into 

\begin{equation}\label{eqn sw}
\textstyle{\bigotimes\limits^k V} =\sum_{\substack{ \lam : \text{partitions of size $k$}}} \lam_{\frak{S}_k} \textstyle{\bigotimes} \lam_{GL(N)}.
\end{equation}

For $\lam$ with $\ell(\lam) \leqq s \leqq n$, 
 the restriction rules (See the proposition 1.5.3 in \cite{kt1}.) tells us that we have 
\[
\lam_{GL(N)} \downarrow^{GL(N)}_{O(N)} = \sum_{\kappa, \mu} LR^{\lam}_{2\kappa, \mu} \mu_{O(N)}.
\]
Here $2\kappa$ denotes the even partition and $LR^{\lam}_{2\kappa, \mu}$ denotes the Littlewood-Richardson coefficient.
So  $\lam_{GL(2n)}\downarrow^{GL(2n)}_{O(2n)}$ contains $\lam_{O(2n)}$ with multiplicity one and $\lam_{O(2n)}$ never appears in the other components 
of the right-hand sides of \eqref{eqn sw}, 
since the size $\vert \lam \vert$ of $\lam$ is equal to $k$. 

From \eqref{eqn odddellam} in Theorem \ref{thm repodd} and \eqref{eqn evendelmun} in Theorem \ref{thm repeven}, $[\Delta, \lam]_{Pin(N)}$ occurs exactly once in $\Delta \bigotimes \lam_{O(2n)}$.

Since $\vert \lam \vert =k$, from \eqref{eqn odddelnat} and \eqref{eqn evendelnat} and \eqref {eqn evendelnnat}, $[\Delta, \lam]_{Pin(N)}$ never appears in 
$\Delta \bigotimes \bigotimes^u V$ for $u < k$. So this representation must be in the kernels of the contractions and $\pr{}$'s.

From the property of the Young symmetrizers,  all the alternating operators  $\Alt{\{i_1,i_2, \ldots, i_r\}}$ ($r > s $) of degree greater than $s$ 
annihilate the space 
${\sum\limits_{\ell(\lam) \leqq s } \lam_{\frak{S}_k} \bigotimes \lam_{GL(2n)}}$. 

So $\lam_{\frak{S}_k} \bigotimes [\Delta, \lam]_{Pin(N)}$ for partitions $\lam$ of size $k$ and of $\ell(\lam) \leqq s$ must be in the space 
$T_{k,s}^{0}$.

Since we already know that on the space $T_{k,s}^{0}$, $\frak{S}_k$ and $Pin(N)$ act as a dual pair and all the irreducible representations $\lam_{\frak{S}_k}$ of $\frak{S}_k$ with $\ell(\lam) \leqq s$ already appear, so we have the theorem.

\end{proof}

Next we move to the case for $Spin(2n)$.

For $Spin(2n)$ representations, we consider the endomorphism $(A \otimes \id) \in \End ( \Delta \bigotimes \bigotimes^k V)$. Here $A$ is the associator for  $\Delta$ and is given by the degree operator  $A [\ur{i_1, i_2, i_3, \ldots, i_r}] = (-1)^r  [\ur{i_1, i_2, i_3, \ldots, i_r}]$.

\begin{lem}\label{lem commute}
The endomorphism $A \otimes \id $ commutes with both of the actions of $\mathbf{CP_k}$ and $Spin(2n)$ on the space $\Delta \bigotimes \bigotimes^k V$. 

Moreover we have 
\[
\pr{p} \circ (A \otimes \id) = (-1)^p A \circ \pr{p}
\]
and 
\[
\inj{q} \circ A = (-1)^q (A \otimes \id) \inj{q}.
\]
\end{lem}

\begin{proof}

It is obvious that ($A \otimes \id $) commutes with the action of the Brauer algebra $\End_{O(2n)} (\bigotimes^k V)$ and the action of $Spin(2n)$, since the action of $Spin(2n)$ preserves the parity $r$, where $r$ is the number of the indices in $[\ur{i_1, i_2, i_3, \ldots, i_r}]$.

We note that from the parity argument, $\psi^{p}_{q}$ corresponding to the isolated points in any generalized Brauer diagrams $y_{inv} \in \mathbf{GB_{k, inv}}$ must satisfy  $p + q \equiv 0 \mod{2}$.

From the formula for $Pin(2n)$ in Theorem \ref{thm psi}, if $p + q \equiv 0 \mod{2}$,   we can easily check  $\psi^{p}_{q} \circ (A \otimes \id) = (A \otimes \id) \circ \psi^{p}_{q}$.

The last two formulas are straightforward.
\end{proof}

So we first decompose the space $T_{k} =\Delta \bigotimes \bigotimes^k V$ into the irreducible constituents under the action of $Pin(2n) \times \mathbf{CP_k}$. Since $Pin(2n)$ and $\mathbf{CP_k}$ are semisimple and act on the space as the dual pair, from the Wedderburn's theorem,  we have  
\[\textstyle{\Delta \bigotimes \bigotimes\limits^k V = \bigoplus\limits_{\substack{\ell(\lam) \leqq n \\ \vert \lam \vert \leqq k}} [\Delta, \lam]_{Pin(2n)} \bigotimes \lam_{\mathbf{CP_k}}.}
\]

Here $\lam_{\mathbf{CP_k}}$ is the irreducible representation of $\mathbf{CP_k}$. 
We consider the subspace $[\Delta, \lam]_{Pin(2n)} \bigotimes \lam_{\mathbf{CP_k}}$ as $Spin(2n) \times \mathbf{CP_k}$ module. Then we have 
\[\textstyle{[\Delta, \lam]_{Pin(2n)} \bigotimes \lam_{\mathbf{CP_k}}
= (1/2 + \lam)^+_{Spin(2n)} \bigotimes \lam_{\mathbf{CP_k}} \bigoplus (1/2 + \delta)^-_{Spin(2n)}\bigotimes \lam_{\mathbf{CP_k}}.}
\]

Since $A \otimes \id $ commutes with the action of $\mathbf{CP_k}$ and the space $[\Delta, \lam]_{Pin(2n)} \bigotimes \lam_{\mathbf{CP_k}}$ is the homogeneous component of the irreducible representation $\lam_{\mathbf{CP_k}}$ of $\mathbf{CP_k}$, we have 
\[\textstyle{
(A \otimes \id)([\Delta, \lam]_{Pin(2n)} \bigotimes \lam_{\mathbf{CP_k}})
=[\Delta, \lam]_{Pin(2n)} \bigotimes \lam_{\mathbf{CP_k}}.}
\]

Since $(A \otimes \id)^2 =\id$, $\pm 1$ eigenspaces of $(A \otimes \id)$ 
in $[\Delta, \lam] \bigotimes \lam_{\mathbf{CP_k}}$ become non-equivalent 
$Spin(2n) \times \mathbf{CP_k}$ modules.

Here $\Delta^+ \bigotimes \bigotimes^k V$ is $+1$ eigenspace and $\Delta^- \bigotimes \bigotimes^k V$ is $-1$ eigenspace of $A \otimes \id$ in 
$\Delta \bigotimes \bigotimes^k V$.

From Theorem \ref{thm repeven}, $(1/2 + \lam)^{\pm}_{Spin(2n)}$ appears in 
 $\Delta^{\pm} \bigotimes \bigotimes^k V$ if and only if $k - \vert \lam \vert \equiv 0 \mod{2}$. Also $(1/2 + \lam)^{\mp}_{Spin(2n)}$ appears in 
 $\Delta^{\pm} \bigotimes \bigotimes^k V$ if and only if $k - \vert \lam \vert \equiv 1 \mod{2}$. 

So  $(-1)^{k- \vert \lam \vert}$ eigenspace of  $(A \otimes \id)$ in $[\Delta, \lam]_{Pin(2n)} \bigotimes \lam_{\mathbf{CP_k}}$ is the irreducible representation $(1/2 + \lam)^+_{Spin(2n)} \bigotimes \lam_{\mathbf{CP_k}}$ and $(-1)^{k- \vert \lam \vert+1}$ eigenspace of  $(A \otimes \id)$ in $[\Delta, \lam]_{Pin(2n)} \bigotimes \lam_{\mathbf{CP_k}}$ is the irreducible representation $(1/2 + \lam)^-_{Spin(2n)} \bigotimes \lam_{\mathbf{CP_k}}$.

We define the extension of the algebra $\mathbf{CP_k}$ by $(A \otimes \id)$.

\begin{defn}\label{def cptilde}
Let $\mathbf{CS_k}$ be the subalgebra of 
$\Hom ( \Delta \bigotimes \bigotimes^k V,  \Delta \bigotimes \bigotimes^k V)$
 generated by $\mathbf{CP_k}$ and $(A \otimes \id)$.
\end{defn}

Then $\mathbf{CS_k}$ and $Spin(2n)$ act on the space 
$\Delta \bigotimes \bigotimes^k V$ as the dual pair.

Then Theorem \ref{thm dualpin} for $Spin(2n)$ is given as follows. 
Since we have Lemma \ref{lem commute}, $T_{k,s}^{0}$ in Definition \ref{def T0} is  $(A \otimes \id)$ stable and 
We decompose $T_{k,s}^{0}$ into the $\pm 1$ eigenspaces of $(A \otimes \id)$ 
and denote them by $T_{k,s}^{0, \pm}$.
Then from Theorem \ref{thm dualpin}, 
if we note that only $\lam$'s considering here satisfy $\vert \lam \vert =k$,
we have the following theorem.

\begin{thm}\label{thm dualspin}
Then the symmetric group of degree $k$ and $Spin(2n)$ becomes the dual pair 
on the subspaces $T_{k,s}^{0, +}$ and $T_{k,s}^{0, -}$ respectively.
Namely,  $T_{k,s}^{0, \pm}$ are decomposed into the direct sum of the tensor products of the irreducible representations with multiplicity free as follows.
\begin{equation}
\textstyle{
T_{k,s}^{0, \pm} = \sum\limits_{\substack{\lam : \text{partitions of size $k$} \\ \ell(\lam) \leqq s }} \lam_{\frak{S}_k} \bigotimes (1/2 + \lam)^{\pm}_{Spin(2n)}.}
\tag{\ref{thm dualspin}.1}\label{eqn dualspin}
\end{equation}

\end{thm}

\end{document}

%% file: spincentfig1.tex
\unitlength 0.1in
\begin{picture}(35.85,15.15)(3.30,-19.15)
%
\special{pn 20}%
\special{sh 1}%
\special{ar 395 400 10 10 0  6.28318530717959E+0000}%
\special{sh 1}%
\special{ar 395 720 10 10 0  6.28318530717959E+0000}%
\special{sh 1}%
\special{ar 715 400 10 10 0  6.28318530717959E+0000}%
\special{sh 1}%
\special{ar 715 720 10 10 0  6.28318530717959E+0000}%
\special{sh 1}%
\special{ar 715 720 10 10 0  6.28318530717959E+0000}%
%
\special{pn 20}%
\special{pa 395 400}%
\special{pa 395 720}%
\special{fp}%
%
\special{pn 20}%
\special{pa 715 400}%
\special{pa 715 720}%
\special{fp}%
\put(5.5500,-8.8000){\makebox(0,0){$y_1$}}%
%
\special{pn 20}%
\special{sh 1}%
\special{ar 1195 400 10 10 0  6.28318530717959E+0000}%
\special{sh 1}%
\special{ar 1195 720 10 10 0  6.28318530717959E+0000}%
\special{sh 1}%
\special{ar 1515 400 10 10 0  6.28318530717959E+0000}%
\special{sh 1}%
\special{ar 1515 720 10 10 0  6.28318530717959E+0000}%
\special{sh 1}%
\special{ar 1515 720 10 10 0  6.28318530717959E+0000}%
%
\special{pn 20}%
\special{pa 1195 400}%
\special{pa 1515 720}%
\special{fp}%
%
\special{pn 20}%
\special{pa 1515 400}%
\special{pa 1195 720}%
\special{fp}%
\put(13.5500,-8.8000){\makebox(0,0){$y_2$}}%
%
\special{pn 20}%
\special{sh 1}%
\special{ar 1995 400 10 10 0  6.28318530717959E+0000}%
\special{sh 1}%
\special{ar 1995 720 10 10 0  6.28318530717959E+0000}%
\special{sh 1}%
\special{ar 2315 400 10 10 0  6.28318530717959E+0000}%
\special{sh 1}%
\special{ar 2315 720 10 10 0  6.28318530717959E+0000}%
\special{sh 1}%
\special{ar 2315 720 10 10 0  6.28318530717959E+0000}%
%
\special{pn 20}%
\special{pa 1995 400}%
\special{pa 2315 400}%
\special{fp}%
%
\special{pn 20}%
\special{pa 1995 720}%
\special{pa 2315 720}%
\special{fp}%
\put(21.5500,-8.8000){\makebox(0,0){$y_3$}}%
%
\special{pn 20}%
\special{sh 1}%
\special{ar 2795 400 10 10 0  6.28318530717959E+0000}%
\special{sh 1}%
\special{ar 2795 720 10 10 0  6.28318530717959E+0000}%
\special{sh 1}%
\special{ar 3115 400 10 10 0  6.28318530717959E+0000}%
\special{sh 1}%
\special{ar 3115 720 10 10 0  6.28318530717959E+0000}%
\special{sh 1}%
\special{ar 3115 720 10 10 0  6.28318530717959E+0000}%
%
\special{pn 20}%
\special{pa 2795 400}%
\special{pa 2795 720}%
\special{fp}%
\put(29.5500,-8.8000){\makebox(0,0){$y_4$}}%
%
\special{pn 20}%
\special{sh 1}%
\special{ar 3595 400 10 10 0  6.28318530717959E+0000}%
\special{sh 1}%
\special{ar 3595 720 10 10 0  6.28318530717959E+0000}%
\special{sh 1}%
\special{ar 3915 400 10 10 0  6.28318530717959E+0000}%
\special{sh 1}%
\special{ar 3915 720 10 10 0  6.28318530717959E+0000}%
\special{sh 1}%
\special{ar 3915 720 10 10 0  6.28318530717959E+0000}%
%
\special{pn 20}%
\special{pa 3915 400}%
\special{pa 3915 720}%
\special{fp}%
\special{pa 3915 400}%
\special{pa 3915 720}%
\special{fp}%
\put(37.5500,-8.8000){\makebox(0,0){$y_5$}}%
%
\special{pn 20}%
\special{sh 1}%
\special{ar 395 1520 10 10 0  6.28318530717959E+0000}%
\special{sh 1}%
\special{ar 395 1840 10 10 0  6.28318530717959E+0000}%
\special{sh 1}%
\special{ar 715 1520 10 10 0  6.28318530717959E+0000}%
\special{sh 1}%
\special{ar 715 1840 10 10 0  6.28318530717959E+0000}%
\special{sh 1}%
\special{ar 707 1840 10 10 0  6.28318530717959E+0000}%
%
\special{pn 20}%
\special{pa 395 1520}%
\special{pa 715 1840}%
\special{fp}%
\put(5.5500,-20.0000){\makebox(0,0){$y_6$}}%
%
\special{pn 20}%
\special{sh 1}%
\special{ar 1195 1520 10 10 0  6.28318530717959E+0000}%
\special{sh 1}%
\special{ar 1195 1840 10 10 0  6.28318530717959E+0000}%
\special{sh 1}%
\special{ar 1515 1520 10 10 0  6.28318530717959E+0000}%
\special{sh 1}%
\special{ar 1515 1840 10 10 0  6.28318530717959E+0000}%
\special{sh 1}%
\special{ar 1515 1840 10 10 0  6.28318530717959E+0000}%
%
\special{pn 20}%
\special{pa 1515 1520}%
\special{pa 1195 1840}%
\special{fp}%
\put(13.5500,-20.0000){\makebox(0,0){$y_7$}}%
%
\special{pn 20}%
\special{sh 1}%
\special{ar 1995 1520 10 10 0  6.28318530717959E+0000}%
\special{sh 1}%
\special{ar 1995 1840 10 10 0  6.28318530717959E+0000}%
\special{sh 1}%
\special{ar 2315 1520 10 10 0  6.28318530717959E+0000}%
\special{sh 1}%
\special{ar 2315 1840 10 10 0  6.28318530717959E+0000}%
\special{sh 1}%
\special{ar 2315 1840 10 10 0  6.28318530717959E+0000}%
%
\special{pn 20}%
\special{pa 1995 1520}%
\special{pa 2315 1520}%
\special{fp}%
\put(21.5500,-20.0000){\makebox(0,0){$y_8$}}%
%
\special{pn 20}%
\special{sh 1}%
\special{ar 2795 1520 10 10 0  6.28318530717959E+0000}%
\special{sh 1}%
\special{ar 2795 1840 10 10 0  6.28318530717959E+0000}%
\special{sh 1}%
\special{ar 3115 1520 10 10 0  6.28318530717959E+0000}%
\special{sh 1}%
\special{ar 3115 1840 10 10 0  6.28318530717959E+0000}%
\special{sh 1}%
\special{ar 3115 1840 10 10 0  6.28318530717959E+0000}%
%
\special{pn 20}%
\special{pa 2795 1840}%
\special{pa 3115 1840}%
\special{fp}%
\put(29.5500,-20.0000){\makebox(0,0){$y_9$}}%
%
\special{pn 20}%
\special{sh 1}%
\special{ar 3595 1520 10 10 0  6.28318530717959E+0000}%
\special{sh 1}%
\special{ar 3595 1840 10 10 0  6.28318530717959E+0000}%
\special{sh 1}%
\special{ar 3915 1520 10 10 0  6.28318530717959E+0000}%
\special{sh 1}%
\special{ar 3915 1840 10 10 0  6.28318530717959E+0000}%
\special{sh 1}%
\special{ar 3915 1840 10 10 0  6.28318530717959E+0000}%
\put(37.5500,-20.0000){\makebox(0,0){$y_{10}$}}%
\end{picture}%

%% file: spincentfig2.tex
\unitlength 0.1in
\begin{picture}(21.22,9.57)(9.80,-11.47)
\put(9.9600,-3.8000){\makebox(0,0)[lb]{1}}%
\put(13.5400,-3.7400){\makebox(0,0)[lb]{2}}%
\put(17.7500,-3.6000){\makebox(0,0)[lb]{3}}%
\put(22.2800,-3.7400){\makebox(0,0)[lb]{4}}%
\put(26.3900,-3.8000){\makebox(0,0)[lb]{5}}%
\put(31.0200,-3.9300){\makebox(0,0)[lb]{6}}%
%
\special{pn 20}%
\special{sh 1}%
\special{ar 996 504 10 10 0  6.28318530717959E+0000}%
\special{sh 1}%
\special{ar 996 504 10 10 0  6.28318530717959E+0000}%
%
\special{pn 20}%
\special{sh 1}%
\special{ar 1412 504 10 10 0  6.28318530717959E+0000}%
\special{sh 1}%
\special{ar 1417 504 10 10 0  6.28318530717959E+0000}%
%
\special{pn 20}%
\special{sh 1}%
\special{ar 1828 504 10 10 0  6.28318530717959E+0000}%
\special{sh 1}%
\special{ar 1828 504 10 10 0  6.28318530717959E+0000}%
%
\special{pn 20}%
\special{sh 1}%
\special{ar 2244 504 10 10 0  6.28318530717959E+0000}%
\special{sh 1}%
\special{ar 2244 504 10 10 0  6.28318530717959E+0000}%
%
\special{pn 20}%
\special{sh 1}%
\special{ar 2660 504 10 10 0  6.28318530717959E+0000}%
\special{sh 1}%
\special{ar 2660 504 10 10 0  6.28318530717959E+0000}%
%
\special{pn 20}%
\special{sh 1}%
\special{ar 3076 504 10 10 0  6.28318530717959E+0000}%
\special{sh 1}%
\special{ar 3076 504 10 10 0  6.28318530717959E+0000}%
%
\special{pn 20}%
\special{pa 996 504}%
\special{pa 1412 1024}%
\special{fp}%
%
\special{pn 20}%
\special{pa 2660 504}%
\special{pa 1828 1024}%
\special{fp}%
%
\special{pn 20}%
\special{sh 1}%
\special{ar 1412 1024 10 10 0  6.28318530717959E+0000}%
\special{sh 1}%
\special{ar 1828 1024 10 10 0  6.28318530717959E+0000}%
\special{sh 1}%
\special{ar 2244 1024 10 10 0  6.28318530717959E+0000}%
\special{sh 1}%
\special{ar 2660 1024 10 10 0  6.28318530717959E+0000}%
\special{sh 1}%
\special{ar 2660 1024 10 10 0  6.28318530717959E+0000}%
%
\special{pn 20}%
\special{sh 1}%
\special{ar 996 1024 10 10 0  6.28318530717959E+0000}%
\special{sh 1}%
\special{ar 996 1024 10 10 0  6.28318530717959E+0000}%
\put(9.8000,-12.1200){\makebox(0,0)[lb]{1}}%
\put(13.9600,-11.9900){\makebox(0,0)[lb]{2}}%
\put(18.1700,-11.9900){\makebox(0,0)[lb]{3}}%
\put(21.4000,-12.2000){\makebox(0,0)[lb]{4}}%
\put(26.9600,-11.8000){\makebox(0,0)[lb]{5}}%
%
\special{pn 20}%
\special{ar 1830 516 421 130  3.2339656 6.1912466}%
%
\special{pn 20}%
\special{ar 2450 994 217 153  0.1437387 2.8653366}%
\end{picture}%

%% file: spincentfig3.tex
\unitlength 0.1in
\begin{picture}(49.75,8.00)(0.06,-10.00)
\put(5.0100,-3.6000){\makebox(0,0){$\emptyset$}}%
%
\special{pn 8}%
\special{pa 981 360}%
\special{pa 1141 360}%
\special{fp}%
\special{sh 1}%
\special{pa 1141 360}%
\special{pa 1074 340}%
\special{pa 1088 360}%
\special{pa 1074 380}%
\special{pa 1141 360}%
\special{fp}%
%
\special{pn 8}%
\special{pa 1301 200}%
\special{pa 1301 360}%
\special{fp}%
\special{pa 1301 360}%
\special{pa 1461 360}%
\special{fp}%
\special{pa 1461 360}%
\special{pa 1461 200}%
\special{fp}%
\special{pa 1461 200}%
\special{pa 1301 200}%
\special{fp}%
%
\special{pn 8}%
\special{pa 1621 360}%
\special{pa 1781 360}%
\special{fp}%
\special{sh 1}%
\special{pa 1781 360}%
\special{pa 1714 340}%
\special{pa 1728 360}%
\special{pa 1714 380}%
\special{pa 1781 360}%
\special{fp}%
%
\special{pn 8}%
\special{pa 1941 200}%
\special{pa 1941 360}%
\special{fp}%
\special{pa 1941 360}%
\special{pa 2101 360}%
\special{fp}%
\special{pa 2101 360}%
\special{pa 2101 200}%
\special{fp}%
\special{pa 2101 200}%
\special{pa 1941 200}%
\special{fp}%
%
\special{pn 8}%
\special{pa 1941 360}%
\special{pa 1941 520}%
\special{fp}%
\special{pa 1941 520}%
\special{pa 2101 520}%
\special{fp}%
\special{pa 2101 520}%
\special{pa 2101 360}%
\special{fp}%
\special{pa 2101 360}%
\special{pa 1941 360}%
\special{fp}%
%
\special{pn 8}%
\special{pa 2261 360}%
\special{pa 2421 360}%
\special{fp}%
\special{sh 1}%
\special{pa 2421 360}%
\special{pa 2354 340}%
\special{pa 2368 360}%
\special{pa 2354 380}%
\special{pa 2421 360}%
\special{fp}%
%
\special{pn 8}%
\special{pa 2581 200}%
\special{pa 2581 360}%
\special{fp}%
\special{pa 2581 360}%
\special{pa 2741 360}%
\special{fp}%
\special{pa 2741 360}%
\special{pa 2741 200}%
\special{fp}%
\special{pa 2741 200}%
\special{pa 2581 200}%
\special{fp}%
%
\special{pn 8}%
\special{pa 2581 360}%
\special{pa 2581 520}%
\special{fp}%
\special{pa 2581 520}%
\special{pa 2741 520}%
\special{fp}%
\special{pa 2741 520}%
\special{pa 2741 360}%
\special{fp}%
\special{pa 2741 360}%
\special{pa 2581 360}%
\special{fp}%
%
\special{pn 8}%
\special{pa 2901 360}%
\special{pa 3061 360}%
\special{fp}%
\special{sh 1}%
\special{pa 3061 360}%
\special{pa 2994 340}%
\special{pa 3008 360}%
\special{pa 2994 380}%
\special{pa 3061 360}%
\special{fp}%
%
\special{pn 8}%
\special{pa 3221 200}%
\special{pa 3221 360}%
\special{fp}%
\special{pa 3221 360}%
\special{pa 3381 360}%
\special{fp}%
\special{pa 3381 360}%
\special{pa 3381 200}%
\special{fp}%
\special{pa 3381 200}%
\special{pa 3221 200}%
\special{fp}%
%
\special{pn 8}%
\special{pa 3541 360}%
\special{pa 3701 360}%
\special{fp}%
\special{sh 1}%
\special{pa 3701 360}%
\special{pa 3634 340}%
\special{pa 3648 360}%
\special{pa 3634 380}%
\special{pa 3701 360}%
\special{fp}%
%
\special{pn 8}%
\special{pa 3861 200}%
\special{pa 3861 360}%
\special{fp}%
\special{pa 3861 360}%
\special{pa 4181 360}%
\special{fp}%
\special{pa 4181 360}%
\special{pa 4181 200}%
\special{fp}%
\special{pa 4181 200}%
\special{pa 3861 200}%
\special{fp}%
\special{pa 4021 200}%
\special{pa 4021 360}%
\special{fp}%
%
\special{pn 8}%
\special{pa 4341 360}%
\special{pa 4501 360}%
\special{fp}%
\special{sh 1}%
\special{pa 4501 360}%
\special{pa 4434 340}%
\special{pa 4448 360}%
\special{pa 4434 380}%
\special{pa 4501 360}%
\special{fp}%
%
\special{pn 8}%
\special{pa 4661 200}%
\special{pa 4661 360}%
\special{fp}%
\special{pa 4661 360}%
\special{pa 4981 360}%
\special{fp}%
\special{pa 4981 360}%
\special{pa 4981 200}%
\special{fp}%
\special{pa 4981 200}%
\special{pa 4661 200}%
\special{fp}%
\special{pa 4821 200}%
\special{pa 4821 360}%
\special{fp}%
%
\special{pn 8}%
\special{pa 4661 360}%
\special{pa 4661 520}%
\special{fp}%
\special{pa 4661 520}%
\special{pa 4821 520}%
\special{fp}%
\special{pa 4821 360}%
\special{pa 4821 520}%
\special{fp}%
%
\special{pn 8}%
\special{pa 181 840}%
\special{pa 341 840}%
\special{fp}%
\special{sh 1}%
\special{pa 341 840}%
\special{pa 274 820}%
\special{pa 288 840}%
\special{pa 274 860}%
\special{pa 341 840}%
\special{fp}%
%
\special{pn 8}%
\special{pa 501 680}%
\special{pa 501 840}%
\special{fp}%
\special{pa 501 840}%
\special{pa 821 840}%
\special{fp}%
\special{pa 821 840}%
\special{pa 821 680}%
\special{fp}%
\special{pa 821 680}%
\special{pa 501 680}%
\special{fp}%
\special{pa 661 680}%
\special{pa 661 840}%
\special{fp}%
%
\special{pn 8}%
\special{pa 501 840}%
\special{pa 501 1000}%
\special{fp}%
\special{pa 501 1000}%
\special{pa 661 1000}%
\special{fp}%
\special{pa 661 840}%
\special{pa 661 1000}%
\special{fp}%
%
\special{pn 8}%
\special{pa 981 840}%
\special{pa 1141 840}%
\special{fp}%
\special{sh 1}%
\special{pa 1141 840}%
\special{pa 1074 820}%
\special{pa 1088 840}%
\special{pa 1074 860}%
\special{pa 1141 840}%
\special{fp}%
%
\special{pn 8}%
\special{pa 1301 680}%
\special{pa 1301 840}%
\special{fp}%
\special{pa 1301 840}%
\special{pa 1621 840}%
\special{fp}%
\special{pa 1621 840}%
\special{pa 1621 680}%
\special{fp}%
\special{pa 1621 680}%
\special{pa 1301 680}%
\special{fp}%
\special{pa 1461 680}%
\special{pa 1461 840}%
\special{fp}%
%
\special{pn 8}%
\special{pa 1301 840}%
\special{pa 1301 1000}%
\special{fp}%
\special{pa 1301 1000}%
\special{pa 1461 1000}%
\special{fp}%
\special{pa 1461 840}%
\special{pa 1461 1000}%
\special{fp}%
%
\special{pn 8}%
\special{pa 1621 840}%
\special{pa 1621 1000}%
\special{fp}%
\special{pa 1621 1000}%
\special{pa 1461 1000}%
\special{fp}%
\special{pa 1461 1000}%
\special{pa 1461 1000}%
\special{fp}%
\end{picture}%

%% file: spincentfig4.tex
\unitlength 0.1in
\begin{picture}(39.50,18.82)(6.00,-24.82)
%
\special{pn 20}%
\special{sh 1}%
\special{ar 600 600 10 10 0  6.28318530717959E+0000}%
\special{sh 1}%
\special{ar 880 600 10 10 0  6.28318530717959E+0000}%
\special{sh 1}%
\special{ar 600 880 10 10 0  6.28318530717959E+0000}%
\special{sh 1}%
\special{ar 880 880 10 10 0  6.28318530717959E+0000}%
\special{sh 1}%
\special{ar 880 880 10 10 0  6.28318530717959E+0000}%
\put(10.4800,-9.8500){\makebox(0,0){$inv$}}%
\put(13.0000,-7.4000){\makebox(0,0){$=$}}%
%
\special{pn 20}%
\special{sh 1}%
\special{ar 1580 600 10 10 0  6.28318530717959E+0000}%
\special{sh 1}%
\special{ar 1860 600 10 10 0  6.28318530717959E+0000}%
\special{sh 1}%
\special{ar 1580 880 10 10 0  6.28318530717959E+0000}%
\special{sh 1}%
\special{ar 1860 880 10 10 0  6.28318530717959E+0000}%
\special{sh 1}%
\special{ar 1860 880 10 10 0  6.28318530717959E+0000}%
\put(20.1400,-10.1300){\makebox(0,0){$rt$}}%
\put(22.8000,-7.4000){\makebox(0,0){$-$}}%
%
\special{pn 20}%
\special{sh 1}%
\special{ar 2560 600 10 10 0  6.28318530717959E+0000}%
\special{sh 1}%
\special{ar 2840 600 10 10 0  6.28318530717959E+0000}%
\special{sh 1}%
\special{ar 2560 880 10 10 0  6.28318530717959E+0000}%
\special{sh 1}%
\special{ar 2840 880 10 10 0  6.28318530717959E+0000}%
\special{sh 1}%
\special{ar 2840 880 10 10 0  6.28318530717959E+0000}%
%
\special{pn 20}%
\special{pa 2560 600}%
\special{pa 2560 880}%
\special{fp}%
\put(29.8000,-9.9200){\makebox(0,0){$rt$}}%
\put(31.2000,-7.4000){\makebox(0,0){$+$}}%
%
\special{pn 20}%
\special{sh 1}%
\special{ar 3400 600 10 10 0  6.28318530717959E+0000}%
\special{sh 1}%
\special{ar 3400 880 10 10 0  6.28318530717959E+0000}%
\special{sh 1}%
\special{ar 3680 600 10 10 0  6.28318530717959E+0000}%
\special{sh 1}%
\special{ar 3680 880 10 10 0  6.28318530717959E+0000}%
\special{sh 1}%
\special{ar 3680 880 10 10 0  6.28318530717959E+0000}%
%
\special{pn 20}%
\special{pa 3400 600}%
\special{pa 3680 880}%
\special{fp}%
\put(38.7600,-9.7800){\makebox(0,0){$rt$}}%
\put(40.3700,-7.5400){\makebox(0,0){$+$}}%
%
\special{pn 20}%
\special{sh 1}%
\special{ar 4240 600 10 10 0  6.28318530717959E+0000}%
\special{sh 1}%
\special{ar 4240 880 10 10 0  6.28318530717959E+0000}%
\special{sh 1}%
\special{ar 4520 600 10 10 0  6.28318530717959E+0000}%
\special{sh 1}%
\special{ar 4520 880 10 10 0  6.28318530717959E+0000}%
\special{sh 1}%
\special{ar 4520 880 10 10 0  6.28318530717959E+0000}%
%
\special{pn 20}%
\special{pa 4520 600}%
\special{pa 4240 880}%
\special{fp}%
\put(47.3000,-9.9200){\makebox(0,0){$rt$}}%
\put(15.1000,-14.3300){\makebox(0,0){$-$}}%
%
\special{pn 20}%
\special{sh 1}%
\special{ar 1720 1300 10 10 0  6.28318530717959E+0000}%
\special{sh 1}%
\special{ar 1720 1580 10 10 0  6.28318530717959E+0000}%
\special{sh 1}%
\special{ar 2000 1300 10 10 0  6.28318530717959E+0000}%
\special{sh 1}%
\special{ar 2000 1580 10 10 0  6.28318530717959E+0000}%
\special{sh 1}%
\special{ar 2000 1573 10 10 0  6.28318530717959E+0000}%
%
\special{pn 20}%
\special{pa 2000 1300}%
\special{pa 2000 1580}%
\special{fp}%
\put(21.7500,-17.3400){\makebox(0,0){$rt$}}%
\put(24.2000,-14.4000){\makebox(0,0){$+$}}%
%
\special{pn 20}%
\special{sh 1}%
\special{ar 2700 1300 10 10 0  6.28318530717959E+0000}%
\special{sh 1}%
\special{ar 2700 1580 10 10 0  6.28318530717959E+0000}%
\special{sh 1}%
\special{ar 2980 1300 10 10 0  6.28318530717959E+0000}%
\special{sh 1}%
\special{ar 2980 1580 10 10 0  6.28318530717959E+0000}%
\special{sh 1}%
\special{ar 2980 1580 10 10 0  6.28318530717959E+0000}%
%
\special{pn 20}%
\special{pa 2700 1300}%
\special{pa 2700 1580}%
\special{fp}%
\special{pa 2980 1300}%
\special{pa 2980 1580}%
\special{fp}%
\put(31.6200,-17.0600){\makebox(0,0){$rt$}}%
\put(33.6500,-14.4000){\makebox(0,0){$-$}}%
%
\special{pn 20}%
\special{sh 1}%
\special{ar 3680 1300 10 10 0  6.28318530717959E+0000}%
\special{sh 1}%
\special{ar 3680 1580 10 10 0  6.28318530717959E+0000}%
\special{sh 1}%
\special{ar 3960 1300 10 10 0  6.28318530717959E+0000}%
\special{sh 1}%
\special{ar 3960 1580 10 10 0  6.28318530717959E+0000}%
\special{sh 1}%
\special{ar 3960 1580 10 10 0  6.28318530717959E+0000}%
%
\special{pn 20}%
\special{pa 3680 1300}%
\special{pa 3960 1580}%
\special{fp}%
\special{pa 3960 1300}%
\special{pa 3680 1580}%
\special{fp}%
\put(41.6300,-16.8500){\makebox(0,0){$rt$}}%
%
\special{pn 20}%
\special{sh 1}%
\special{ar 600 2140 10 10 0  6.28318530717959E+0000}%
\special{sh 1}%
\special{ar 600 2420 10 10 0  6.28318530717959E+0000}%
\special{sh 1}%
\special{ar 880 2140 10 10 0  6.28318530717959E+0000}%
\special{sh 1}%
\special{ar 880 2420 10 10 0  6.28318530717959E+0000}%
\special{sh 1}%
\special{ar 880 2420 10 10 0  6.28318530717959E+0000}%
%
\special{pn 20}%
\special{pa 600 2140}%
\special{pa 600 2420}%
\special{fp}%
\put(10.4800,-25.6700){\makebox(0,0){$inv$}}%
\put(13.1400,-23.0100){\makebox(0,0){$=$}}%
%
\special{pn 20}%
\special{sh 1}%
\special{ar 1860 2140 10 10 0  6.28318530717959E+0000}%
\special{sh 1}%
\special{ar 2140 2140 10 10 0  6.28318530717959E+0000}%
\special{sh 1}%
\special{ar 1860 2420 10 10 0  6.28318530717959E+0000}%
\special{sh 1}%
\special{ar 2140 2420 10 10 0  6.28318530717959E+0000}%
\special{sh 1}%
\special{ar 2140 2420 10 10 0  6.28318530717959E+0000}%
%
\special{pn 20}%
\special{pa 1860 2140}%
\special{pa 1860 2420}%
\special{fp}%
\put(22.8000,-25.3200){\makebox(0,0){$rt$}}%
\put(15.7300,-22.8700){\makebox(0,0){$-$}}%
\put(25.1400,-22.8000){\makebox(0,0){$+$}}%
%
\special{pn 20}%
\special{sh 1}%
\special{ar 2794 2140 10 10 0  6.28318530717959E+0000}%
\special{sh 1}%
\special{ar 2794 2420 10 10 0  6.28318530717959E+0000}%
\special{sh 1}%
\special{ar 3074 2140 10 10 0  6.28318530717959E+0000}%
\special{sh 1}%
\special{ar 3074 2420 10 10 0  6.28318530717959E+0000}%
\special{sh 1}%
\special{ar 3074 2420 10 10 0  6.28318530717959E+0000}%
%
\special{pn 20}%
\special{pa 2794 2140}%
\special{pa 2794 2420}%
\special{fp}%
\special{pa 3074 2140}%
\special{pa 3074 2420}%
\special{fp}%
\put(32.5600,-25.4600){\makebox(0,0){$rt$}}%
\end{picture}%

%% file: spincentfig5.tex
\unitlength 0.1in
\begin{picture}(38.42,6.23)(0.05,-9.23)
%
\special{pn 20}%
\special{sh 1}%
\special{ar 900 300 10 10 0  6.28318530717959E+0000}%
\special{sh 1}%
\special{ar 1180 300 10 10 0  6.28318530717959E+0000}%
\special{sh 1}%
\special{ar 900 580 10 10 0  6.28318530717959E+0000}%
\special{sh 1}%
\special{ar 1180 580 10 10 0  6.28318530717959E+0000}%
\special{sh 1}%
\special{ar 1180 580 10 10 0  6.28318530717959E+0000}%
%
\special{pn 20}%
\special{pa 900 300}%
\special{pa 1180 300}%
\special{fp}%
%
\special{pn 20}%
\special{sh 1}%
\special{ar 900 629 10 10 0  6.28318530717959E+0000}%
\special{sh 1}%
\special{ar 1180 629 10 10 0  6.28318530717959E+0000}%
\special{sh 1}%
\special{ar 900 909 10 10 0  6.28318530717959E+0000}%
\special{sh 1}%
\special{ar 1180 909 10 10 0  6.28318530717959E+0000}%
\special{sh 1}%
\special{ar 1180 909 10 10 0  6.28318530717959E+0000}%
%
\special{pn 20}%
\special{pa 1190 643}%
\special{pa 1190 923}%
\special{fp}%
\put(14.6000,-6.2200){\makebox(0,0){$ = $}}%
\put(19.9200,-6.2200){\makebox(0,0){$(X - 1)$}}%
%
\special{pn 20}%
\special{sh 1}%
\special{ar 2321 475 10 10 0  6.28318530717959E+0000}%
\special{sh 1}%
\special{ar 2601 475 10 10 0  6.28318530717959E+0000}%
\special{sh 1}%
\special{ar 2321 755 10 10 0  6.28318530717959E+0000}%
\special{sh 1}%
\special{ar 2601 755 10 10 0  6.28318530717959E+0000}%
\special{sh 1}%
\special{ar 2601 755 10 10 0  6.28318530717959E+0000}%
%
\special{pn 20}%
\special{pa 2321 475}%
\special{pa 2601 475}%
\special{fp}%
\put(27.7600,-6.1500){\makebox(0,0){$+$}}%
\put(32.2000,-6.1500){\makebox(0,0){$(X - 1)$}}%
%
\special{pn 20}%
\special{sh 1}%
\special{ar 3550 468 10 10 0  6.28318530717959E+0000}%
\special{sh 1}%
\special{ar 3830 468 10 10 0  6.28318530717959E+0000}%
\special{sh 1}%
\special{ar 3550 748 10 10 0  6.28318530717959E+0000}%
\special{sh 1}%
\special{ar 3830 748 10 10 0  6.28318530717959E+0000}%
\special{sh 1}%
\special{ar 3830 748 10 10 0  6.28318530717959E+0000}%
%
\special{pn 20}%
\special{pa 3550 468}%
\special{pa 3830 468}%
\special{fp}%
%
\special{pn 20}%
\special{pa 3567 741}%
\special{pa 3847 741}%
\special{fp}%
\put(5.0000,-6.1000){\makebox(0,0){$y_5 y_8 =$}}%
\end{picture}%

%% file: spincentfig6.tex
\unitlength 0.1in
\begin{picture}(46.38,28.08)(1.20,-29.08)
%
\special{pn 20}%
\special{sh 1}%
\special{ar 127 2097 10 10 0  6.28318530717959E+0000}%
\special{sh 1}%
\special{ar 267 2097 10 10 0  6.28318530717959E+0000}%
\special{sh 1}%
\special{ar 407 2097 10 10 0  6.28318530717959E+0000}%
\special{sh 1}%
\special{ar 547 2097 10 10 0  6.28318530717959E+0000}%
\special{sh 1}%
\special{ar 687 2097 10 10 0  6.28318530717959E+0000}%
\special{sh 1}%
\special{ar 827 2097 10 10 0  6.28318530717959E+0000}%
\special{sh 1}%
\special{ar 967 2097 10 10 0  6.28318530717959E+0000}%
\special{sh 1}%
\special{ar 967 2097 10 10 0  6.28318530717959E+0000}%
%
\special{pn 20}%
\special{sh 1}%
\special{ar 127 1663 10 10 0  6.28318530717959E+0000}%
\special{sh 1}%
\special{ar 267 1663 10 10 0  6.28318530717959E+0000}%
\special{sh 1}%
\special{ar 407 1663 10 10 0  6.28318530717959E+0000}%
\special{sh 1}%
\special{ar 547 1663 10 10 0  6.28318530717959E+0000}%
\special{sh 1}%
\special{ar 687 1663 10 10 0  6.28318530717959E+0000}%
\special{sh 1}%
\special{ar 827 1663 10 10 0  6.28318530717959E+0000}%
\special{sh 1}%
\special{ar 967 1663 10 10 0  6.28318530717959E+0000}%
\special{sh 1}%
\special{ar 967 1663 10 10 0  6.28318530717959E+0000}%
%
\special{pn 20}%
\special{sh 1}%
\special{ar 1387 1663 10 10 0  6.28318530717959E+0000}%
\special{sh 1}%
\special{ar 1527 1663 10 10 0  6.28318530717959E+0000}%
\special{sh 1}%
\special{ar 1667 1663 10 10 0  6.28318530717959E+0000}%
\special{sh 1}%
\special{ar 1807 1663 10 10 0  6.28318530717959E+0000}%
\special{sh 1}%
\special{ar 1947 1663 10 10 0  6.28318530717959E+0000}%
\special{sh 1}%
\special{ar 2087 1663 10 10 0  6.28318530717959E+0000}%
\special{sh 1}%
\special{ar 2227 1663 10 10 0  6.28318530717959E+0000}%
\special{sh 1}%
\special{ar 2227 1663 10 10 0  6.28318530717959E+0000}%
%
\special{pn 20}%
\special{sh 1}%
\special{ar 2654 1663 10 10 0  6.28318530717959E+0000}%
\special{sh 1}%
\special{ar 2794 1663 10 10 0  6.28318530717959E+0000}%
\special{sh 1}%
\special{ar 2934 1663 10 10 0  6.28318530717959E+0000}%
\special{sh 1}%
\special{ar 3074 1663 10 10 0  6.28318530717959E+0000}%
\special{sh 1}%
\special{ar 3214 1663 10 10 0  6.28318530717959E+0000}%
\special{sh 1}%
\special{ar 3354 1663 10 10 0  6.28318530717959E+0000}%
\special{sh 1}%
\special{ar 3494 1663 10 10 0  6.28318530717959E+0000}%
\special{sh 1}%
\special{ar 3494 1663 10 10 0  6.28318530717959E+0000}%
%
\special{pn 20}%
\special{sh 1}%
\special{ar 3914 1663 10 10 0  6.28318530717959E+0000}%
\special{sh 1}%
\special{ar 4054 1663 10 10 0  6.28318530717959E+0000}%
\special{sh 1}%
\special{ar 4194 1663 10 10 0  6.28318530717959E+0000}%
\special{sh 1}%
\special{ar 4334 1663 10 10 0  6.28318530717959E+0000}%
\special{sh 1}%
\special{ar 4474 1663 10 10 0  6.28318530717959E+0000}%
\special{sh 1}%
\special{ar 4614 1663 10 10 0  6.28318530717959E+0000}%
\special{sh 1}%
\special{ar 4754 1663 10 10 0  6.28318530717959E+0000}%
\special{sh 1}%
\special{ar 4754 1663 10 10 0  6.28318530717959E+0000}%
%
\special{pn 20}%
\special{ar 3914 1663 0 0  3.1415927 6.2831853}%
%
\special{pn 20}%
\special{ar 1450 1649 70 77  3.1415927 6.2831853}%
%
\special{pn 20}%
\special{ar 764 1663 214 77  3.1415927 6.2831853}%
%
\special{pn 20}%
\special{ar 2024 1670 214 77  3.1415927 6.2831853}%
%
\special{pn 20}%
\special{ar 3277 1663 213 77  3.1415927 6.2831853}%
%
\special{pn 20}%
\special{ar 4544 1663 214 77  3.1415927 6.2831853}%
%
\special{pn 20}%
\special{ar 2790 1656 137 77  3.1415927 6.2831853}%
%
\special{pn 20}%
\special{ar 477 2097 350 133  3.1415927 6.2831853}%
%
\special{pn 20}%
\special{ar 4190 1663 284 133  3.1415927 6.2831853}%
%
\special{pn 20}%
\special{ar 757 2097 213 77  3.1415927 6.2831853}%
\put(5.6800,-18.5900){\makebox(0,0){$y_1$}}%
\put(18.4200,-18.5900){\makebox(0,0){$y_2$}}%
\put(30.3200,-18.7300){\makebox(0,0){$y_3$}}%
\put(30.3200,-18.7300){\makebox(0,0){$y_3$}}%
\put(43.6900,-18.7300){\makebox(0,0){$y_4$}}%
\put(5.6100,-23.0700){\makebox(0,0){$y_5$}}%
%
\special{pn 20}%
\special{sh 1}%
\special{ar 1387 2090 10 10 0  6.28318530717959E+0000}%
\special{sh 1}%
\special{ar 1527 2090 10 10 0  6.28318530717959E+0000}%
\special{sh 1}%
\special{ar 1667 2090 10 10 0  6.28318530717959E+0000}%
\special{sh 1}%
\special{ar 1807 2090 10 10 0  6.28318530717959E+0000}%
\special{sh 1}%
\special{ar 1947 2090 10 10 0  6.28318530717959E+0000}%
\special{sh 1}%
\special{ar 2087 2090 10 10 0  6.28318530717959E+0000}%
\special{sh 1}%
\special{ar 2227 2090 10 10 0  6.28318530717959E+0000}%
\special{sh 1}%
\special{ar 2227 2090 10 10 0  6.28318530717959E+0000}%
\put(18.2800,-23.2100){\makebox(0,0){$z_1$}}%
%
\special{pn 20}%
\special{ar 1810 2091 144 55  0.1382430 3.0301680}%
%
\special{pn 20}%
\special{sh 1}%
\special{ar 2640 2104 10 10 0  6.28318530717959E+0000}%
\special{sh 1}%
\special{ar 2780 2104 10 10 0  6.28318530717959E+0000}%
\special{sh 1}%
\special{ar 2920 2104 10 10 0  6.28318530717959E+0000}%
\special{sh 1}%
\special{ar 3060 2104 10 10 0  6.28318530717959E+0000}%
\special{sh 1}%
\special{ar 3200 2104 10 10 0  6.28318530717959E+0000}%
\special{sh 1}%
\special{ar 3340 2104 10 10 0  6.28318530717959E+0000}%
\special{sh 1}%
\special{ar 3480 2104 10 10 0  6.28318530717959E+0000}%
\special{sh 1}%
\special{ar 3480 2104 10 10 0  6.28318530717959E+0000}%
\put(30.7400,-23.9100){\makebox(0,0){$z_2$}}%
%
\special{pn 20}%
\special{ar 3064 2105 143 55  0.1394401 3.0309354}%
%
\special{pn 20}%
\special{ar 2850 2110 203 120  6.2339640 6.2831853}%
\special{ar 2850 2110 203 120  0.0000000 3.1864057}%
%
\special{pn 20}%
\special{sh 1}%
\special{ar 3893 2076 10 10 0  6.28318530717959E+0000}%
\special{sh 1}%
\special{ar 4033 2076 10 10 0  6.28318530717959E+0000}%
\special{sh 1}%
\special{ar 4173 2076 10 10 0  6.28318530717959E+0000}%
\special{sh 1}%
\special{ar 4313 2076 10 10 0  6.28318530717959E+0000}%
\special{sh 1}%
\special{ar 4453 2076 10 10 0  6.28318530717959E+0000}%
\special{sh 1}%
\special{ar 4593 2076 10 10 0  6.28318530717959E+0000}%
\special{sh 1}%
\special{ar 4733 2076 10 10 0  6.28318530717959E+0000}%
\special{sh 1}%
\special{ar 4733 2076 10 10 0  6.28318530717959E+0000}%
\put(43.4100,-23.6300){\makebox(0,0){$z_3$}}%
%
\special{pn 20}%
\special{ar 4316 2077 144 55  0.1382430 3.0301680}%
%
\special{pn 20}%
\special{ar 4180 2083 144 55  0.1382430 3.0309354}%
%
\special{pn 20}%
\special{sh 1}%
\special{ar 127 2671 10 10 0  6.28318530717959E+0000}%
\special{sh 1}%
\special{ar 267 2671 10 10 0  6.28318530717959E+0000}%
\special{sh 1}%
\special{ar 407 2671 10 10 0  6.28318530717959E+0000}%
\special{sh 1}%
\special{ar 547 2671 10 10 0  6.28318530717959E+0000}%
\special{sh 1}%
\special{ar 687 2671 10 10 0  6.28318530717959E+0000}%
\special{sh 1}%
\special{ar 827 2671 10 10 0  6.28318530717959E+0000}%
\special{sh 1}%
\special{ar 967 2671 10 10 0  6.28318530717959E+0000}%
\special{sh 1}%
\special{ar 967 2671 10 10 0  6.28318530717959E+0000}%
\put(5.6800,-29.0200){\makebox(0,0){$z_4$}}%
%
\special{pn 20}%
\special{ar 550 2672 144 55  0.1382430 3.0301680}%
%
\special{pn 20}%
\special{ar 694 2664 144 55  0.1382430 3.0309354}%
%
\special{pn 20}%
\special{sh 1}%
\special{ar 1373 2664 10 10 0  6.28318530717959E+0000}%
\special{sh 1}%
\special{ar 1513 2664 10 10 0  6.28318530717959E+0000}%
\special{sh 1}%
\special{ar 1653 2664 10 10 0  6.28318530717959E+0000}%
\special{sh 1}%
\special{ar 1793 2664 10 10 0  6.28318530717959E+0000}%
\special{sh 1}%
\special{ar 1933 2664 10 10 0  6.28318530717959E+0000}%
\special{sh 1}%
\special{ar 2073 2664 10 10 0  6.28318530717959E+0000}%
\special{sh 1}%
\special{ar 2213 2664 10 10 0  6.28318530717959E+0000}%
\special{sh 1}%
\special{ar 2213 2664 10 10 0  6.28318530717959E+0000}%
\put(18.0000,-29.9300){\makebox(0,0){$z_5$}}%
%
\special{pn 20}%
\special{ar 1796 2665 144 55  0.1382430 3.0301680}%
%
\special{pn 20}%
\special{ar 2010 2685 203 120  6.2339640 6.2831853}%
\special{ar 2010 2685 203 120  0.0000000 3.1864057}%
%
\special{pn 20}%
\special{sh 1}%
\special{ar 319 170 10 10 0  6.28318530717959E+0000}%
\special{sh 1}%
\special{ar 459 170 10 10 0  6.28318530717959E+0000}%
\special{sh 1}%
\special{ar 599 170 10 10 0  6.28318530717959E+0000}%
\special{sh 1}%
\special{ar 739 170 10 10 0  6.28318530717959E+0000}%
\special{sh 1}%
\special{ar 879 170 10 10 0  6.28318530717959E+0000}%
\special{sh 1}%
\special{ar 1019 170 10 10 0  6.28318530717959E+0000}%
\special{sh 1}%
\special{ar 1159 170 10 10 0  6.28318530717959E+0000}%
\special{sh 1}%
\special{ar 1159 170 10 10 0  6.28318530717959E+0000}%
%
\special{pn 20}%
\special{sh 1}%
\special{ar 319 310 10 10 0  6.28318530717959E+0000}%
\special{sh 1}%
\special{ar 459 310 10 10 0  6.28318530717959E+0000}%
\special{sh 1}%
\special{ar 599 310 10 10 0  6.28318530717959E+0000}%
\special{sh 1}%
\special{ar 739 310 10 10 0  6.28318530717959E+0000}%
\special{sh 1}%
\special{ar 879 310 10 10 0  6.28318530717959E+0000}%
\special{sh 1}%
\special{ar 1019 310 10 10 0  6.28318530717959E+0000}%
\special{sh 1}%
\special{ar 1159 310 10 10 0  6.28318530717959E+0000}%
\special{sh 1}%
\special{ar 1159 310 10 10 0  6.28318530717959E+0000}%
%
\special{pn 20}%
\special{pa 319 170}%
\special{pa 599 310}%
\special{fp}%
%
\special{pn 20}%
\special{ar 946 163 206 63  3.1415927 6.2831853}%
%
\special{pn 20}%
\special{ar 1020 310 139 63  6.2831853 6.2831853}%
\special{ar 1020 310 139 63  0.0000000 3.1415927}%
%
\special{pn 20}%
\special{sh 1}%
\special{ar 326 527 10 10 0  6.28318530717959E+0000}%
\special{sh 1}%
\special{ar 466 527 10 10 0  6.28318530717959E+0000}%
\special{sh 1}%
\special{ar 606 527 10 10 0  6.28318530717959E+0000}%
\special{sh 1}%
\special{ar 746 527 10 10 0  6.28318530717959E+0000}%
\special{sh 1}%
\special{ar 886 527 10 10 0  6.28318530717959E+0000}%
\special{sh 1}%
\special{ar 1026 527 10 10 0  6.28318530717959E+0000}%
\special{sh 1}%
\special{ar 1166 527 10 10 0  6.28318530717959E+0000}%
\special{sh 1}%
\special{ar 1166 527 10 10 0  6.28318530717959E+0000}%
%
\special{pn 20}%
\special{sh 1}%
\special{ar 326 660 10 10 0  6.28318530717959E+0000}%
\special{sh 1}%
\special{ar 466 660 10 10 0  6.28318530717959E+0000}%
\special{sh 1}%
\special{ar 606 660 10 10 0  6.28318530717959E+0000}%
\special{sh 1}%
\special{ar 746 660 10 10 0  6.28318530717959E+0000}%
\special{sh 1}%
\special{ar 886 660 10 10 0  6.28318530717959E+0000}%
\special{sh 1}%
\special{ar 1026 660 10 10 0  6.28318530717959E+0000}%
\special{sh 1}%
\special{ar 1166 660 10 10 0  6.28318530717959E+0000}%
\special{sh 1}%
\special{ar 1166 660 10 10 0  6.28318530717959E+0000}%
%
\special{pn 20}%
\special{pa 886 520}%
\special{pa 746 660}%
\special{fp}%
%
\special{pn 20}%
\special{ar 816 522 210 65  2.9617392 6.1929548}%
%
\special{pn 20}%
\special{ar 747 664 139 59  6.0872582 6.2831853}%
\special{ar 747 664 139 59  0.0000000 3.1415927}%
\put(12.1000,-4.3000){\makebox(0,0)[lb]{$ = 3(X-2)(X-3)$}}%
%
\special{pn 20}%
\special{sh 1}%
\special{ar 2680 260 10 10 0  6.28318530717959E+0000}%
\special{sh 1}%
\special{ar 2820 260 10 10 0  6.28318530717959E+0000}%
\special{sh 1}%
\special{ar 2960 260 10 10 0  6.28318530717959E+0000}%
\special{sh 1}%
\special{ar 3100 260 10 10 0  6.28318530717959E+0000}%
\special{sh 1}%
\special{ar 3240 260 10 10 0  6.28318530717959E+0000}%
\special{sh 1}%
\special{ar 3380 260 10 10 0  6.28318530717959E+0000}%
\special{sh 1}%
\special{ar 3520 260 10 10 0  6.28318530717959E+0000}%
\special{sh 1}%
\special{ar 3520 260 10 10 0  6.28318530717959E+0000}%
%
\special{pn 20}%
\special{sh 1}%
\special{ar 2680 400 10 10 0  6.28318530717959E+0000}%
\special{sh 1}%
\special{ar 2820 400 10 10 0  6.28318530717959E+0000}%
\special{sh 1}%
\special{ar 2960 400 10 10 0  6.28318530717959E+0000}%
\special{sh 1}%
\special{ar 3100 400 10 10 0  6.28318530717959E+0000}%
\special{sh 1}%
\special{ar 3240 400 10 10 0  6.28318530717959E+0000}%
\special{sh 1}%
\special{ar 3380 400 10 10 0  6.28318530717959E+0000}%
\special{sh 1}%
\special{ar 3520 400 10 10 0  6.28318530717959E+0000}%
\special{sh 1}%
\special{ar 3520 400 10 10 0  6.28318530717959E+0000}%
%
\special{pn 20}%
\special{ar 3307 253 206 63  3.1415927 6.2831853}%
%
\special{pn 20}%
\special{pa 2680 253}%
\special{pa 3107 400}%
\special{fp}%
%
\special{pn 20}%
\special{ar 3095 415 145 118  6.2831853 6.2831853}%
\special{ar 3095 415 145 118  0.0000000 3.3435268}%
\put(17.2000,-10.0000){\makebox(0,0)[lb]{$- (X-2)(X-3)(X-4) {\ds \sum_{i=1}^5 \sum_{j=1}^5 (-1)^{i+j}z_iy_j}$}}%
\put(1.2000,-13.2000){\makebox(0,0)[lb]{\text{Here $y_j$ denotes the upper row and $z_i$ denotes the lower row given as follows.}}}%
\end{picture}%